\documentclass{amsart}
\usepackage{amssymb,stmaryrd}
\usepackage{amsfonts}
\usepackage{amstext}
\usepackage{algorithmic}
\usepackage{algorithm}
\usepackage{graphicx}
\usepackage{epstopdf}
\usepackage[all]{xy}
\usepackage{MnSymbol}
\usepackage{color,float}

\numberwithin{equation}{section}
\newtheorem{theorem}{Theorem}[section]
\newtheorem{lemma}[theorem]{Lemma}
\newtheorem{proposition}[theorem]{Proposition}

\theoremstyle{definition}

\theoremstyle{remark}
\newtheorem{remark}[theorem]{\bf{Remark}}

\newcommand{\C}{{\mathbb{C}}}
\newcommand{\Z}{{\mathbb{Z}}}

\newcommand{\F}{{\mathbb{F}}}

\renewcommand{\>}{{\rangle}}

\newcommand{\wedgeq}{{\wedge\kern-5pt\cdot}}
\newcommand{\cg}{{g}}

\newcommand{\tens}{\otimes}

\newcommand{\id}{{\rm id}}

\newcommand{\extd}{{\rm d}}

\newcommand{\eps}{\epsilon}

\newcommand{\Vol}{{\rm Vol}}
\renewcommand{\cg}{{\bf g}}

\begin{document}
\title[Digital finite quantum Riemannian geometries]{Digital finite quantum
Riemannian geometries}
\keywords{noncommutative geometry, quantum groups, quantum gravity, finite fields,digital electronics, quantum computing}
\subjclass[2000]{Primary 81R50, 58B32, 83C57}
\author{Shahn Majid \& Anna Pacho{\l }}
\address{Queen Mary University of London, Mile End Rd, London E1 4NS, UK}
\email{s.majid@qmul.ac.uk, a.pachol@qmul.ac.uk}

\begin{abstract} We study bimodule quantum Riemannian geometries over the field $\Bbb F_2$ of two elements as the extreme case of a finite-field adaptation of noncommutative-geometric methods for physics. We classify all parallelisable such geometries for coordinate algebras up to vector space dimension $n\le 3$, finding a rich moduli of examples for $n=3$ and top form degree 2, including many that are not flat. Their coordinate algebras are commutative but their differentials are not. We also study the quantum Laplacian $\Delta=(\ ,\ )\nabla\extd$ 
on our models and characterise when it has a massive eigenvector.  \end{abstract}

\maketitle

\section{Introduction}

Quantum or noncommutative geometry has been extensively developed since the 1980s and from several different motivations. One is the widely accepted view that momentum space could be curved and conversely that
spacetime could be noncommutative or `quantum' due to
Planck scale corrections \cite{Ma:pla} and that this could in principle be measurable by secondary effects such as in \cite{AmeMa,KohBra}.   There is also evidence for such a {\em quantum spacetime hypothesis} in 3D quantum gravity \cite{FL,FreMa}, where the theory is better understood albeit topological, while in 4D the hypothesis gives a route into elements of effective quantum gravity without knowing the full theory. There are also plenty of other potential applications  of quantum geometry, such to the geometry of other quantum integrable systems where quantum groups play a central role,  and in principle to actual quantum systems where there can be phase spaces that also have a metric structure (as in the Kahler manifold case \cite{Brod}) which we can now follow into the quantum algebra of observables. So far, this last possibility has been little studied but there are instances, such as the quantum Hall effect \cite{Bel} and more recently the fractional quantum Hall effect \cite{Lut} where noncommutative geometry is thought to be relevant. Mathematically speaking, the most well-known approaches are the one of Connes, coming out of operator algebras as a `spectral triple' generalisation of the Dirac operator \cite{Con}, and a more algebraic `quantum groups' approach coming out of experience with quantum groups but not limited to them, see e.g.  \cite{Ma:ltcc, BegMa1,BegMa2, Ma:gra,MaTao1}. The starting point here is a bimodule $\Omega^1$ of quantum differential forms over an algebra $A$, a quantum metric in $\Omega^1\tens_A\Omega^1$ and a `bimodule connection'  \cite{DVM,Mou} for the formulation of a quantum Levi-Civita connection.  

In this paper we explore a natural `spin-off' direction initiated in \cite{BasMa,MaPac} in which we use this second  approach to  noncommutative geometry to replace the field $\C$ of complex numbers that we work over by the field ${\mathbb{F}}_2=\{0,1\}$ of two elements. Note that algebraic geometers have long formulated the geometry of commutative algebras over general fields, but such methods become very involved over the field $\F_{p^d}$ of prime power $p^d$ elements, arguably due to incompatibility of the positive characteristic with usual commutative differentials as explained in \cite{BasMa}. By contrast, we offer a very different and more  calculable approach to geometry over $\F_{p^d}$ of interest {\em even} when the coordinate algebra is commutative, made possible because in our case differentials need not commute with functions. Indeed, in this paper we only look at `coordinate algebras' $A$ of vector space dimension $\le 3$, which limits us mostly to commutative $A$, but with noncommuting differentials.  Restricting to finite dimensional $A$ corresponds in classical geometry to  functions on a finite set, so in physics this would be like making a finite model of space or spacetime. Also note that being commutative over $\F_{p^d}$ does {\em not} mean that an algebra is actually that of functions on some space. In short, much of the character and richness of noncommutative geometry carries over to this apparently simple setting over a finite field like $\F_{p^d}$. Making discrete geometry models in physics is of course an established way to tame infinities that might be encountered, and this too is included in quantum geometry \cite{Ma:gra,Ma:sq},  but now we go further and in an orthogonal direction of making the field finite. For example, if `path integrals' are done with function values in a finite field then integration of the function value at each point becomes a finite sum (such possible applications remain to be explored as one should first formulate path integrals in way that does not need $A$ to be an actual algebra of functions). 

The first step is to have a sense of what this moduli space of $\F_{p^d}$-geometries looks like, which  is the aim of present paper for $A$ of small dimension and in the `digital case' of  $\F_2$-geometries. Remarkably, the answer turns out to be quite rich and to 
include many nonflat ones (all with commutative $A$).  Another immediate application of having a repertoire of digital quantum geometries is that they can be used to test ideas and conjectures in the general theory if we expect them to hold for any field, even if we are mainly interested in the theory over $\C$. Indeed, the nonlinear nature of the quantum Levi-Civita condition in quantum Riemannian geometry makes it very hard to solve by analytic means for a general quantum metric, with the result that to date only the
square graph and the integer lattice $\Z$ were fully solved for general metrics \cite{Ma:sq,Ma:haw}. It also means that the quantum Levi-Civita connection may not exist, and even if it does, it may not be unique. This is a phenomenon that deserves more study and we will see it in our models. Having so few general examples to work with also makes it hard to further develop the theory in a convincing way (notably the correct notion of stress energy and conservation laws are poorly understood).  It is also possible, geometry being ubiquitous in science and engineering, that there could be applications of  $\F_{p^d}$ and $\F_2$ quantum geometries in their own right. One of these could be to transfer geometric ideas into digital electronics as explained in \cite{MaPac}. Why exactly one would want to do this remains to be seen, but one area of application could be to build \emph{digital
quantum computing} gates as analogues of what we may wish to build in an actual quantum computer and with potentially some of the benefits and as `training wheels' for the real thing.  

The paper begins in Section~\ref{secpre} with some preliminary definitions from the constructive `bimodule' 
approach to noncommutative Riemannian geometry but written entirely in
tensor-calculus terms, which is needed for computer implementation. The paper is a sequel to \cite
{MaPac} where we classified digital geometries on $A={\mathbb{F}}
_2[x_1,\cdots,x_n]$ with differential calculus defined by commutative $n$
-dimensional algebras $(V,\circ)$ and we already know all possible such
algebras over ${\mathbb{F}}_2$ up to $n\le 4$ from that work. The difference
now is that $A=(V,\circ)$ itself is our coordinate algebra on which we do
the noncommutative differential geometry. We also do not need $A$ to be commutative, but there are no
noncommutative unital algebras as $n=2$ and only one over $\F_2$ at  $n=3$, and this does not appear to admit an interesting quantum geometry, so in practice  $A$ will mainly be commutative. The classification is then done by computer
methods using Mathematica to try all possible values for the Christoffel symbols, 
with some work by hand as a check of the implementation. Results for $n=2$ are at the start
of  Section~\ref{secclass} and the more interesting results for $n=3$ are in Sections~\ref{secDmodel} -- \ref{secFmodel}. 
After this, Section~\ref{seclap} 
computes the Laplacian for the main examples and Section~\ref{secricci} studies
the Ricci tensor and scalar. A summary of the results is provided in the concluding Section~\ref{secfin}.

\section{Preliminaries}

\label{secpre}

This is a short account of the bimodule approach in \cite{Ma:ltcc,BegMa1,BegMa2}
but in an explicit form with structure constants and over a general field $k$. We will then look for solutions over ${\mathbb{F}}_2$ in later sections.
Compared to \cite{MaPac}, we denote the finite-dimensional algebra by $A$, its product by
omission and the identity element by 1 (in the previous work they were denoted $\circ
$, $e$).

(i) Let $\{x^\mu\}$ be a basis of our algebra $A$ with $e=x^0=1$ and $
\mu =0,\cdots ,n-1$, where $n=\mathrm{dim}\ A$. We write structure constants by 
\begin{equation}
x^{\mu } x^{\nu }=V^{\mu \nu }{}_{\rho }x^{\rho },\quad V^{\mu \nu }{}_\rho\in k.
\end{equation}
We sum over repeated indices. For a unital commutative associative algebra we of course need 
\begin{equation}
V^{0\mu}{}_\nu=\delta^\mu_\nu,\quad V^{\mu \nu }{}_\rho=V^{\nu
\mu}{}_\rho,\quad V^{\rho \nu }{}_{\lambda }V^{\lambda
\mu }{}_{\gamma }=V^{\nu \mu }{}_{\lambda }V^{\rho\lambda}{}_{\gamma }. 
\end{equation}

Next we define the differential structure by specifying a space of 1-forms $
\Omega^1$ which we assume has a basis $\{\omega^i\}$, $i=1,\cdots,m$ over $A$
, where necessarily $m\le n-1$ is the dimension of the calculus over $A$.
The case $m=n-1$ is the `universal calculus' and any other is a quotient of
this. Our assumption is that $\Omega^1=A.\{\omega^i\}$ as a free left module by
the product in $A$ and we also require a right action of $A$ which we
specify by structure constants, along with structure constants for the
exterior differential $\extd:A\to \Omega^1$, by 
\begin{equation}
\omega ^{i}.x^{\mu }=a^{i\mu }{}_{\nu j}x^{\nu }.\omega ^{j},\quad \extd
x^{\mu }=d^{\mu }{}_{\nu i}x^{\nu }.\omega ^{i},\quad a^{i\mu }{}_{\nu j},\
d^{\mu }{}_{\nu i}\in k.
\end{equation}
Such a calculus is called `left parallelisable'. In the nicest case, the commutation rules for moving
algebra generators to the left can be inverted, so we can equally take $
\Omega^1=\{\omega^i\}.A$. The two actions are required to obey the axioms of a
bimodule $a((\extd b)c)=(a\extd b)c$ (the left and right actions
commute with each other) and $\extd$ is required to obey the Leibniz rule $\extd (ab)=(\extd
a)b+a\extd b$, which become
\begin{equation}
a^{i\mu }{}_{\sigma j}a^{j\nu }{}_{\tau k}V^{\sigma \tau }{}_{\eta
}=a^{i\rho }{}_{\eta k}V^{\mu \nu }{}_{\rho },\quad V^{\mu \nu }{}_{\rho
}d^{\rho }{}_{\sigma j}=d^{\mu }{}_{\lambda i}a^{i\nu }{}_{\gamma
j}V^{\lambda \gamma }{}_{\sigma }+d^{\nu }{}_{\delta j}V^{\mu \delta
}{}_{\sigma }.
\end{equation}
We also want $\Omega^1$ to be spanned by things of the form $a\extd b$
(the so-called surjectivity axiom) and optionally we ask for the calculus to be
connected in the sense that only the constant function is killed by $\mathrm{
d}$. These translate respectively as 
\begin{equation}
B^{\mu \nu }{}_{\eta i}:=V^{\mu \rho }{}_{\eta }d^{\nu }{}_{\rho i},\quad 
\mathrm{rank}\ nm
\end{equation}
\begin{equation}
d^{\mu }{}_{\nu i},\quad \mathrm{1-dimensional\ null\ space}.
\end{equation}

Note that for any unital algebra, we have $\extd1=0$ and $\omega
^{i}.1=\omega ^{i}=1.\omega ^{i}$ which also implies that $\extd x^{\mu
}.1=1.\extd x^{\mu }=\extd x^{\mu }$ and fixes some of the structure
constants. Also note that there is a standard notion of isomorphism $(\Omega^1,\extd)\cong (\Omega^{1'},\extd')$ between differential structures in their classification over a fixed algebra, namely invertible bimodule maps $\phi:\Omega^1\to\Omega^{1'}$ such that $\phi(\extd a)=\extd' \phi(a)$ for all $a\in A$. All calculi of maximal dimension $m=n-1$ are isomorphic (to the universal calculus) so in this case there is one distinct calculus. But for $n=3,m=1$ we will point out the one case where there are nontrivial isomorphisms among our solutions. A weaker equivalence is if there is an algebra automorphism $A\to A$ compatible with a map between potentially different differential calculi before and after, which we will again see some  examples of but which we do not study systematically. 

(ii) Next we define a metric as an invertible element $\cg\in
\Omega^1\otimes_A\Omega^1$ in the sense that there exists a bimodule `inner product' map $(\ ,\ ):\Omega^1\otimes_A\Omega^1\to A$ such that 
$((\eta,\ )\otimes \mathrm{id})\cg=\eta=(\mathrm{id}\otimes(\ ,\eta))\cg$ for all $\eta\in\Omega^1$. 
Note that invertibility in this sense forces $\cg$ to be
central (to commute with functions)\cite{BegMa2}. To explain these ideas in terms of structure constants, it is useful to  use two different expansions of the metric,
\begin{equation}
\cg=g_{\mu ij}x^{\mu }.\omega ^{i}\otimes \omega ^{j}=\omega^i.\tilde g_{\mu
ij}x^\mu\otimes\omega^j;\quad g_{\mu ij}=a^{k\nu}{}_{\mu i}\tilde g_{\nu kj};\quad g_{\mu i j }, \tilde g_{\mu i j }\in
k. 
\end{equation}
Then invertibility of the metric amounts to existence of an algebra-valued matrix $g^{ij}=(\omega^i,
\omega^j)\in A$ which is inverse to $\tilde g_{ij}=\tilde g_{\mu ij}x^\mu\in
A$ in the sense that $g^{ij}\tilde g_{jk}=\delta_{ik}=\tilde g_{ij}g^{jk}$. Meanwhile, centrality of the metric implied by invertibility appears as 
\begin{equation}\label{centralg} g_{\mu ij}a_{\quad \lambda k}^{j\nu }a_{\quad \gamma m}^{i\lambda
}V_{\quad \sigma }^{\mu \gamma }=g_{\rho mk}V_{\quad \sigma }^{\nu \rho }.\end{equation}
In practice, it is easier to impose centrality first and then check for invertibility amongst the results. Also, given a quantum metric, one has a natural `metric quantum dimension'
\begin{equation}\label{qdim} \underline{\rm dim}=(\ , \ )(\cg)\in k.\end{equation}

Next, by a `left connection' on $\Omega ^{1}$ we mean $\nabla :\Omega
^{1}\rightarrow \Omega ^{1}\otimes _{A}\Omega ^{1}$ such that 
\begin{equation*}
\nabla (a\omega )=a(\nabla \omega )+\extd a\otimes \omega,\quad \forall
a\in A,\ \omega \in \Omega ^{1}.
\end{equation*}
The reader can think of a covariant derivative $\nabla_v:\Omega^1\to
\Omega^1 $ along all vector fields $v$ but we specify all of these together
by an extra left-most copy of $\Omega^1$ waiting to be evaluated against any
vector field. This is a normal approach in noncommutative
geometry and allows us to dispense with vector fields entirely. By `bimodule
connection' we mean a left connection such that in addition there is a 
bimodule map $\sigma :\Omega ^{1}\otimes _{A}\Omega ^{1}\rightarrow \Omega
^{1}\otimes _{A}\Omega ^{1}$ such that
\begin{equation*}
\nabla (\omega a)=(\nabla \omega )a+\sigma (\omega \otimes \extd
a),\quad \forall a\in A,\ \omega \in \Omega ^{1}. 
\end{equation*}
If a left connection admits such a $\sigma $ then the latter is unique,
hence this is a property of $\nabla $ and not further data. For this reason, we will not list the $\sigma$ explicitly in the tables (they are not particularly illuminating) but they should always be understood for a quantum Levi-Civita connection, and can be computed from 
\begin{equation}\label{sigma} \sigma(\omega^i\otimes \extd x^\mu)=\extd x^\mu\otimes \omega^i +\nabla[\omega^i,x^\mu]-[\nabla \omega^i,x^\mu].\end{equation}
Next, any bimodule connection canonically extends to $\nabla \cg:=(\nabla\otimes\mathrm{id})\cg+(\sigma\otimes
\mathrm{id})(\mathrm{id}\otimes\nabla)\cg$ and we say in this case that $\nabla$ is
{\em metric-compatible} if $\nabla \cg=0$. We write structure constants for the
connection as 
\begin{equation}
\nabla \omega ^{i}=\Gamma^i{}_{\nu km}x^{\nu }\omega ^{k}\otimes \omega
^{m},\quad \sigma \left( \omega ^{i}\otimes \omega ^{j}\right)
=\sigma^{ij}{}_{\mu km}x^{\mu }\omega ^{k}\otimes \omega ^{m}. 
\end{equation}

\begin{lemma}
\label{lemma2.1} In terms of our structure constants, the full left
connection is 
\begin{equation*}
\nabla (x^{\mu }.\omega ^{j})=\left( V^{\mu \rho }{}_{\nu }\Gamma^j{}_{\rho
km }+d^\mu{}_{\nu k}\delta _{m}^{j}\right) x^{\nu }.\omega ^{k}\otimes
\omega ^{m}. 
\end{equation*}
A bimodule connection has $\sigma$ obeying 
\begin{equation*}
\Gamma^i{}_{\nu st}a^{t\mu }{}_{\lambda m}a^{s\lambda}{}_{\gamma k}V^{\nu
\gamma }{}_{\rho }+d^\mu{}_{\lambda j}a^{i\lambda }{}_{\gamma
r}\sigma^{rj}{}_{\beta km}V^{\gamma \beta }{}_{\rho }=a^{i\mu}{}_{\nu
j}\left( V^{\nu \alpha }{}_{\rho }\Gamma^j{}_{\alpha km }+d^\nu{}_{\rho
k}\delta _{m}^{j}\right)
\end{equation*}
\begin{equation*}
\left(a^{j\mu}{}_{\nu k}a^{i\nu}{}_{\rho m}\sigma^{mk}{}_{\sigma
st}-\sigma^{ij}{}_{\rho km}a^{m\mu}{}_{\nu t}a^{k\nu}{}_{\sigma
s}\right)V^{\rho\sigma}{}_\tau=0.
\end{equation*}
for the connection to obey the right Leibniz rule and for $\sigma$ to be a bimodule map, respectively.
Such a connection is metric compatible if 
\begin{equation*}
g_{\mu mn}d^\mu{}_{\nu k}+g_{\mu in}\Gamma^i{}_{\rho km}V^{\mu\rho}{}_{\nu
}+g_{\mu ij}\Gamma^j{}_{\beta sn}a^{i\beta }{}_{\alpha
t}\sigma^{ts}{}_{\sigma km }V^{\mu \alpha }{}_{\rho }V^{\rho \sigma }{}_{\nu
} =0.
\end{equation*}
\end{lemma}

\proof Here by the connection derivation rule, $\nabla (x^{\mu }.\omega
^{j})=x^{\mu }(\nabla \omega ^{j})+\extd x^{\mu }\otimes \omega
^{j}=x^{\mu }(\Gamma ^{j}{}_{\rho km}x^{\rho }\omega ^{k}\otimes \omega
^{m})+d^{\mu }{}_{\nu k}x^{\nu }.\omega ^{k}\otimes \omega ^{j}=V^{\mu \rho
}{}_{\nu }x^{\nu }\Gamma ^{j}{}_{\rho km}\omega ^{k}\otimes \omega
^{m}+d^{\mu }{}_{\nu k}x^{\nu }.\omega ^{k}\otimes \delta _{m}^{j}\omega
^{m}$, giving the expression stated. For a bimodule connection, we
additionally require 
\begin{align*}
\nabla (\omega ^{i}x^{\mu })& =(\nabla \omega ^{i})x^{\mu }+\sigma (\omega
^{i}\otimes \extd x^{\mu })=(\Gamma ^{i}{}_{\nu st}x^{\nu }\omega
^{s}\otimes \omega ^{t}).x^{\mu }+d^{\mu }{}_{\lambda j}\sigma (\omega
^{i}.x^{\lambda }\otimes \omega ^{j}) \\
& =\Gamma ^{i}{}_{\nu st}a^{t\mu }{}_{\lambda r}x^{\nu }\omega
^{s}.x^{\lambda }\otimes \omega ^{r}+d^{\mu }{}_{\lambda j}\sigma
(a^{i\lambda }{}_{\gamma r}x^{\gamma }\omega ^{r}\otimes \omega ^{j}) \\
& =\Gamma ^{i}{}_{\nu st}a^{t\mu }{}_{\lambda r}a^{s\lambda }{}_{\gamma
j}x^{\nu }x^{\gamma }\omega ^{j}\otimes \omega ^{r}+d^{\mu }{}_{\lambda
j}a^{i\lambda }{}_{\gamma r}x^{\gamma }\sigma (\omega ^{r}\otimes \omega
^{j}) \\
& =\Gamma ^{i}{}_{st\nu }a^{t\mu }{}_{\lambda m}a^{s\lambda }{}_{\gamma
k}V^{\nu \gamma }{}_{\rho }x^{\rho }\omega ^{k}\otimes \omega ^{m}+d^{\mu
}{}_{\lambda j}a^{i\lambda }{}_{\gamma r}\sigma ^{rj}{}_{\beta km}V^{\gamma
\beta }{}_{\rho }x^{\rho }\omega ^{k}\otimes \omega ^{m} \\
& =\left( \Gamma ^{i}{}_{\nu st}a^{t\mu }{}_{\lambda m}a^{s\lambda
}{}_{\gamma k}V^{\nu \gamma }{}_{\rho }+d^{\mu }{}_{\lambda j}a^{i\lambda
}{}_{\gamma r}\sigma ^{rj}{}_{\beta km}V^{\gamma \beta }{}_{\rho }\right)
x^{\rho }\omega ^{k}\otimes \omega ^{m} \\
& =\nabla (a^{i\mu }{}_{\nu j}x^{\nu }\omega ^{j})=a^{i\mu }{}_{\nu j}\left(
V^{\nu \alpha }{}_{\rho }\Gamma ^{j}{}_{\alpha km}+d^{\nu }{}_{\rho k}\delta
_{m}^{j}\right) x^{\rho }.\omega ^{k}\otimes \omega ^{m}
\end{align*}
where the last line uses the calculus commutation relations and our previous
formula for the connection. Comparing gives the first condition stated
which, if $\sigma $ exists, characterises it. We also need $\sigma $ to be
well-defined as a bimodule map. This comes down to equality of the
expressions 
\begin{align*}
\sigma (\omega ^{i}\otimes \omega ^{j}x^{\mu })& =a^{j\mu }{}_{\nu k}\sigma
(\omega ^{i}x^{\nu }\otimes \omega ^{k}) \\
& =a^{j\mu }{}_{\nu k}a^{i\nu }{}_{\rho m}x^{\rho }\sigma (\omega
^{m}\otimes \omega ^{k})=a^{j\mu }{}_{\nu k}a^{i\nu }{}_{\rho m}x^{\rho
}\sigma ^{mk}{}_{\sigma st}x^{\sigma }\omega ^{s}\otimes \omega ^{t} \\
& =a^{j\mu }{}_{\nu k}a^{i\nu }{}_{\rho m}\sigma ^{mk}{}_{\sigma st}V^{\rho
\sigma }{}_{\tau }x^{\tau }\omega ^{s}\otimes \omega ^{t} \\
\sigma (\omega ^{i}\otimes \omega ^{j})x^{\mu }& =\sigma ^{ij}{}_{\rho
km}x^{\rho }\omega ^{k}\otimes \omega ^{m}x^{\mu } \\
& =\sigma ^{ij}{}_{\rho km}x^{\rho }\omega ^{k}\otimes a^{m\mu }{}_{\nu
t}x^{\nu }\omega ^{t}=\sigma ^{ij}{}_{\rho km}x^{\rho }a^{m\mu }{}_{\nu
t}a^{k\nu }{}_{\sigma s}x^{\sigma }\omega ^{s}\otimes \omega ^{t} \\
& =\sigma ^{ij}{}_{\rho km}a^{m\mu }{}_{\nu t}a^{k\nu }{}_{\sigma s}V^{\rho
\sigma }{}_{\tau }x^{\tau }\omega ^{s}\otimes \omega ^{t}
\end{align*}
which is the second condition stated. Once we have a bimodule connection,
the metric compatibility makes sense and reads 
\begin{align*}
\nabla \cg& =\nabla (g_{\mu ij}x^{\mu }\omega ^{i}\otimes \omega ^{j}) \\
& =g_{\mu ij}d^{\mu }{}_{\nu k}x^{\nu }\omega ^{k}\otimes \omega ^{i}\otimes
\omega ^{j}+g_{\mu ij}x^{\mu }\Gamma ^{i}{}_{\rho km}x^{\rho }\omega
^{k}\otimes \omega ^{m}\otimes \omega ^{j}+g_{\mu ij}x^{\mu }\sigma (\omega
^{i}\otimes \Gamma ^{j}{}_{\beta sn}x^{\beta }\omega ^{s})\otimes \omega ^{n}
\end{align*}
and the last term here computes further to 
\begin{equation*}
g_{\mu ij}x^{\mu }\Gamma ^{j}{}_{\beta sn}a^{i\beta }{}_{\alpha t}x^{\alpha
}\sigma (\omega ^{t}\otimes \omega ^{s})\otimes \omega ^{n}=g_{\mu ij}x^{\mu
}\Gamma ^{j}{}_{\beta sn}a^{i\beta }{}_{\alpha t}x^{\alpha }\sigma
^{ts}{}_{\sigma km}x^{\sigma }\omega ^{k}\otimes \omega ^{m}\otimes \omega
^{n}.
\end{equation*}
We then replace products giving application of $V$ and relabel so that all
terms are multiples of $x^{\nu }\omega ^{k}\otimes \omega ^{m}\otimes \omega
^{n}$. This then gives the condition stated for $\nabla \cg=0$. \endproof

Note that $\sigma$ in principle depends on the connection via (\ref{sigma}), making the metric compatibility condition quadratic in $\Gamma$. In general, however, it is more realistic to find the possible $\sigma$ 
first as an auxiliary variable.  For $m=1$, we drop the form indices so that 
\begin{equation*}
\nabla \omega =\Gamma _{\nu }x^{\nu }\omega \otimes \omega ,\quad \sigma
\left( \omega \otimes \omega \right) =\sigma _{\beta }x^{\beta }\omega
\otimes \omega ,\quad \cg=g_{\gamma }x^{\gamma }\omega \otimes \omega
\end{equation*}
and the first stated bimodule connection condition and the metric condition
become respectively 
\begin{equation*}
\Gamma _{\nu }a^{\mu }{}_{\lambda }a^{\lambda }{}_{\gamma }V^{\nu \gamma
}{}_{\rho }+d^{\mu }{}_{\lambda }a^{\lambda }{}_{\gamma }\sigma _{\beta
}V^{\gamma \beta }{}_{\rho }=a^{\mu }{}_{\nu }\left( V^{\nu \alpha }{}_{\rho
}\Gamma {}_{\alpha }+d^{\nu }{}_{\rho }\right)
\end{equation*}
\begin{equation*}
g_{\mu }d^{\mu }{}_{\nu }+g_{\mu }\Gamma {}_{\rho }V^{\mu \rho }{}_{\nu
}+g_{\mu }\Gamma _{\beta }a^{\beta }{}_{\alpha }\sigma _{\sigma }V^{\mu
\alpha }{}_{\rho }V^{\rho \sigma }{}_{\nu }=0,
\end{equation*}
while the second stated bimodule connection condition drops out provided the algebra is
commutative (which we have not assumed in proving the lemma but which in practice will be our
case of interest). 

(iii) Next we suppose that $A=\Omega ^{0}$ and $\Omega ^{1}$ are part of a
differential graded or `exterior algebra' $\Omega =\oplus _{i}\Omega ^{i}$
with a wedge product $\wedge $ and with $\extd$ extending so as to obey $
\extd^{2}=0 $ and the graded-Leibniz rule. This requires a little more
data to specify at least $\Omega ^{2}$. We will limit attention to the case
where $\Omega ^{2}=A.\Vol$ is a 1-dimensional free module over $A$, with basis $
\Vol$ say (and we assume that this is also true from the right). For $m=1
$ it is natural to take $\Vol=\omega \wedge \omega $ but we could take any
invertible function of this, while more generally it depends on
the calculus if $\Vol$ exists and how it looks. If it does, we define 
\begin{equation*}
\omega ^{i}\wedge \omega ^{j}=\epsilon ^{ij}{}_{\mu }x^{\mu }\Vol
=\epsilon ^{ij}\Vol,\quad \epsilon ^{ij}{}_{\mu }\in k,\ \epsilon
^{ij}\in A.
\end{equation*}
We want an associative wedge product including the action by elements of $A$, in which case
centrality of the volume form comes down to certain commutation relations between $\omega^i\wedge
\omega^j $, namely
\begin{equation}\label{epsrel}
\epsilon ^{ij} x^{\mu }=a^{j\mu }{}_{\nu k}a^{i\nu }{}_{\rho m}x^{\rho
}\epsilon ^{mk}.
\end{equation}
This property holds for some examples for $m=1$ presented in Sec. 3.2, while for $m=2$ in Sec. 3.3--3.5 it is imposed as a desirable condition as we focus on better behaved calculi. In general, one can suppose that $\epsilon ^{ij}$ is invertible with $
\epsilon _{ij}\epsilon ^{jk}=\delta _{ik}=\epsilon ^{ij}\epsilon _{jk}$. 
Finally, we want $\extd$ to extend with $\extd^{2}=0$, which we
write as 
\begin{equation}  \label{dext1}
\extd\omega ^{i}=\tau^i{}_{\mu }x^{\mu }\Vol=\tau ^{i}\Vol;\quad d^{\mu }{}_{\nu i}x^{\nu }\tau ^{i}+d^{\mu }{}_{\nu i}d^{\nu }{}_{\rho
j}x^\rho\epsilon ^{ji}=0.
\end{equation}
The extension of $\extd$ to general 1-forms is by the Leibniz rule and
this has to be consistent with the bimodule commutation relations, which is
\begin{equation}  \label{dext2}
\tau ^{i}a^{j_{0}\mu }{}_{\nu k}a^{i_{0}\nu }{}_{\rho m}x^{\rho }\epsilon
^{mk}-d^{\mu }{}_{\nu j}a^{i\nu }{}_{\rho k}x^{\rho }\epsilon ^{kj}=a^{i\mu
}{}_{\nu j}d^{\nu }{}_{\rho k}x^{\rho }\epsilon ^{jk}+a^{i\mu }{}_{\nu
j}x^{\nu }\tau ^{i}.
\end{equation}
Note that $\tau^i,\eps^{ij}$ are elements of $A$ so there are $V$ products if we wish to write equations
involving these more explicitly. If $\omega^i=\extd x^i$ then $\tau^i=0$ and in this case
if $\eps^{ij}$ is a solution then so is any invertible function times it, which corresponds to the same $\Omega^2$ with a different volume form. One could have other choices e.g. $\Omega ^{2}=0$ for which such
restrictions and our further conditions below involving $\wedge $ would be
empty.

For $m=1$, we drop the indices so that $\omega\wedge\omega=\eps\Vol$ where $\eps\in A$ is invertible. Then
(\ref{epsrel})--(\ref{dext2}) become
\begin{equation}\label{epsrel1}
\eps x^{\mu }=a^{\mu }{}_{\nu }a^{\nu }{}_{\rho }x^{\rho
}\eps 
\end{equation}
\begin{equation}  \label{dext1m1}
d^{\mu }{}_{\nu}x^{\nu }\tau+d^{\mu }{}_{\nu}d^{\nu }{}_{\rho}x^\rho\epsilon=0
\end{equation}
\begin{equation}  \label{dext2m1}
\tau a^{\mu }{}_{\nu }a^{\nu }{}_{\rho }x^{\rho }\eps-d^{\mu }{}_{\nu }a^{\nu
}{}_{\rho }x^{\rho }\eps=a^{\mu }{}_{\nu }d^{\nu }{}_{\rho }x^{\rho }\eps+a^{\mu
}{}_{\nu }x^{\nu }\tau 
\end{equation}
with $\eps=1$ the canonical choice.

(iv) Once we have specified at least $\Omega^2$, we normally ask for our metric to be
`quantum symmetric' in the sense 
\begin{equation}  \label{qs}
\wedge (\cg)=0;\quad g_{\mu ij}V^{\mu\nu}{}_\rho \epsilon^{ij}{}_\nu=0
\end{equation}
or in $A$-valued notation $g_{ij}\epsilon^{ij}=0$. Also, for any left connection, we define its torsion by $T_{\nabla
}:=\wedge\nabla -\extd$. This commutes with left multiplication by $A$, so for it to 
 vanish it is enough to look on basis elements.  A \emph{quantum Levi-Civita connection} (QLC) is then defined as one where $
\nabla \cg=T_{\nabla }=0$, where the latter torsion freeness in our terms is 
\begin{equation}  \label{Tor}
\Gamma ^{i}{}_{\mu jk}x^{\mu }V^{\mu\nu}{}_\rho \epsilon ^{jk}{}_\nu=\tau ^{i}{}_\rho.
\end{equation}
There is also a weaker notion of \emph{weak quantum
Levi-Civita connection} (WQLC) where we have a left connection which is
torsion free and also \emph{cotorsion free} in the sense 
$
coT_{\nabla }:=\left( \extd\otimes \mathrm{id}-\mathrm{id}\wedge \nabla
\right) \cg=0$, which in our terms is  
\begin{equation}  \label{coT}
g_{\gamma im}d^{\gamma }{}_{\lambda k}\epsilon^{ki}{}_\sigma 
V^{\lambda \sigma}{}_\theta+g_{\gamma im}\tau^{i}{}_\lambda V^{\gamma \lambda }{}_\theta =g_{\gamma ij}V^{\gamma \sigma}{}_
\alpha a^{i\lambda }{}_{\sigma r}\Gamma^{j}{}_{\lambda
km}\epsilon^{rk}{}_\beta V^{\alpha \beta }{}_\theta. 
\end{equation}
This can be useful when $\nabla $ is not a bimodule connection, e.g. when one
does not exist, or when it is computationally too hard to search at first for
all bimodule connections. In this case one can first try to classify all
WQLCs and then see which of them are QLCs. 

For $m=1$ with $\Omega^2$ one-dimensional as above, (\ref{qs})  cannot hold
and (\ref{Tor})--(\ref{coT}) become
\begin{equation}  \label{Tom1}
\Gamma _{\mu }\eps_\nu V^{\mu\nu}{}_\rho=\tau _{\rho }\quad \Leftrightarrow \quad \Gamma\eps=\tau
\end{equation}
\begin{equation}  \label{coTm1}
g_{\gamma}d^{\gamma }{}_{\lambda}\epsilon_\sigma 
V^{\lambda \sigma}{}_\theta+g_{\gamma}\tau_\lambda V^{\gamma \lambda }{}_\theta =g_{\gamma}V^{\gamma \sigma}{}_
\alpha a^{\lambda }{}_{\sigma }\Gamma_{\lambda}\epsilon_\beta V^{\alpha \beta }{}_\theta\quad \Leftrightarrow \quad g'\eps+g\tau=g\bar\Gamma\eps
\end{equation}
where if $\tau=0$ then (\ref{Tom1}) cannot hold unless $\Gamma=0$, while (\ref{coTm1}) becomes $g'=g\bar\Gamma$.
Here we adopted a shorthand $f^{\prime }=f_{\mu }d^{\mu }{}_{\nu }x^{\nu }$
and $\bar f=f_{\mu }a^{\mu }{}_{\nu }x^{\nu }$ for any $f=f_{\mu }x^{\mu }\in
A $, so that $\extd f=f^{\prime }\omega$ and $\omega f=\bar f\omega $.

(v) Once we have constructed our geometries, we will be interested in the
geometric Laplacian and the curvature. These are given respectively by 
\begin{equation*}
\Delta =(\ ,\ )\nabla \extd:A\rightarrow A,\quad R_{\nabla }=(\extd
\otimes \mathrm{id}-\mathrm{id}\wedge \nabla )\nabla :\Omega ^{1}\rightarrow
\Omega ^{2}\otimes _{A}\Omega ^{1}.
\end{equation*}
In terms of components, these are 
\begin{equation}  \label{laplacian}
\Delta x^{\mu }=g_{\tau
}{}^{km} d^{\mu }{}_{\alpha i}(V^{\alpha \rho }{}_{\sigma }\Gamma
^{i}{}_{\rho km}+d^{\alpha }{}_{\sigma k}\delta_{m}^{i})V^{\tau \sigma }{}_{\nu }\,x^{\nu }
\end{equation}
\begin{equation*}  \label{curvm2}
R_{\nabla }\omega ^{i}=\rho ^{i}{}_{j}{}_{\mu }x^{\mu }\Vol\otimes
\omega ^{j}=\rho ^{i}{}_{j}\Vol\otimes \omega ^{j};
\end{equation*}
\begin{equation}  \label{rhom2}
\quad \rho ^i{}_{j\beta}=\Gamma^i{}_{\mu kj}d^{\mu }{}_{\nu n}V^{\nu \alpha }{}_\beta 
\epsilon^{nk}{}_\alpha +\Gamma^i{}_{\mu kj}\tau^k{}_\alpha V^{\mu \alpha }{}_\beta
-\Gamma^i{}_{\mu km}\Gamma^m{}_{\lambda pj} a^{k\lambda }{}_{\alpha s}V^{\mu \alpha }{}_{\sigma }\epsilon^{sp}{}_\theta V^{\sigma \theta }{}_\beta .
\end{equation}

When $m=1$, these become
\[\Delta x^{\mu }=\bar g^{-1}_{\tau
} d^{\mu }{}_{\alpha}(V^{\alpha \rho }{}_{\sigma }\Gamma_{\rho}+d^{\alpha }{}_{\sigma})V^{\tau \sigma }{}_{\nu }\,x^{\nu }
;\quad  \Delta f=\bar g^{-1}(f''+f'\Gamma)\]
for all $f\in A$, where $(\omega,\omega)=\bar g^{-1}$ is the inverse metric, and 
$R_{\nabla }\omega =\rho 
\Vol\otimes \omega$ with
\begin{equation*} \rho_\beta=\Gamma_{\mu}d^{\mu }{}_{\nu}\epsilon_\alpha V^{\nu \alpha }{}_\beta 
 +\Gamma_{\mu}\tau_\alpha V^{\mu \alpha }{}_\beta 
-\Gamma_{\mu}\Gamma_{\lambda} a^{\lambda }{}_{\alpha }V^{\mu \alpha }{}_{\sigma }\epsilon_\theta V^{\sigma \theta }{}_\beta;\quad \rho =\Gamma^{\prime }\eps+\Gamma\tau -\Gamma\bar \Gamma \eps. 
\end{equation*}

It should be stressed that these definitions are not ad-hoc; they are
part of a general noncommutative or `quantum' Riemannian geometry that
applies to a large range of examples, including $q$-deformation ones and graph
geometries.

(vi) The Ricci tensor is less well understood but the proposal in \cite{BegMa2} is to define it with respect to a `lifting' bimodule map $i:\Omega^2\to \Omega^1\tens_A\Omega^1$ 
such that $\wedge\circ i=\id$. If we assume that $\Omega^2$ is one-dimensional with central basis $\Vol$
then we write $i(\Vol)=I_{ij}\omega^i\otimes\omega^j$ for some central element of $\Omega^1\otimes_A\Omega^1$ such that $\wedge(I)=\Vol$. The latter is explicitly,
\begin{equation}\label{wedgeI} I_{ij}\eps^{ij}=1;\quad  I_{ij\mu}\eps^{ij}{}_\nu V^{\mu\nu}{}_\alpha=\delta_{\alpha,0}.\end{equation}
Note that such $I$ is generally not unique, namely we can add
any functional multiple $\gamma \cg$ for $\gamma\in A$ if $\cg$ is central and quantum
symmetric (which may not always be the case). Then 
\[
\mathrm{Ricci}=g_{ij}((\omega^i ,\ )\tens\id)(i\tens\id)R_\nabla \omega^j=g_{ij}(\omega^i ,\rho^j{}_k I_{mn}\omega^m)\omega^n\otimes
\omega^k;\]
\begin{equation}\label{ricciV}
\mathrm{Ricci}_{ij}=g_{\alpha mn}\rho ^n{}_{j\beta} I_{\mu li}V^{\beta\mu}{}_{\nu}a^{m\nu}{}_{\gamma k}g_{\eta}{}^{kl}V^{\alpha \gamma}{}_{\delta}V^{\delta\eta}{}_{\zeta} x^{\zeta} 
\end{equation}
where we used the bimodule commutation relations and $g^{ij}=(\omega^i,
\omega^j)$ as inverse to $\tilde g_{ij}$ according to the general analysis
above, and $\rho ^n{}_{j\beta}$ is as in \eqref{rhom2}. The main idea then is to use the freedom in $I$ to adjust 
\textrm{Ricci} to have the same quantum symmetry and centrality as $\cg$. This may or may not be possible, and if it is it may not be uniquely so. The Ricci scalar is defined as
\begin{equation}\label{S} S=(\ ,\ ){\rm Ricci}={\rm Ricci}_{ij}\tilde g^{ij}.\end{equation}

If $m=1$ and $\Omega^2$ is 1-dimensional then we have $i(\Vol)=\eps^{-1} \omega\tens\omega $ as the unique
choice and\begin{equation*}
\mathrm{Ricci}=g \bar\rho \bar\eps^{-1} \bar g^{-1}\omega\tens\omega. \end{equation*}
This has the same (not quantum symmetric) form as the metric in this case. 
\section{Geometries with $n\le 3$}

\label{secclass}

Here we look at the noncommutative Riemannian geometry of commutative
unital algebras up to dimension $n\le 3$. For $n=1$, there is just $k.1$ for
any field $k$, $\extd 1=0$ and no metric. For $n=2$ we do the
calculation by hand as a check of some computer methods and for $n=3$ we
then proceed entirely by computer (using Mathematica and R) for $m=1$ and
partly by hand and computer in nontrivial cases with $m=2$.

\subsection{Classification of geometries for $n=2$}

\label{Sec3.1}

For $n=2$ over any field $k$ and up to normalisation of the element $x\ne 1$, we have (i) $x^2=\lambda $ for $\lambda\in k$ (some of these could be the
same up to normalisation) or (ii) $x^2=\lambda +x$ where $\lambda\in k$.
Note that if $\lambda=0$ in case (i) then $(1+x)^2=1+2x=-1+2(1+x)$ which
is in case (ii) after normalisation if 2 is invertible in the field and in
case (i) with $\lambda=1$ if $2$ is not invertible in the field. Similarly if $\lambda=1$
over ${\mathbb{F}}_2$ then this is equivalent to $\lambda=0$ by a change in
variable. Over ${\mathbb{F}}_2$ this means three possibilities, but for
future reference we work with general $k$ for as long as we can. In each
case there is a unique non-zero differential calculus $\Omega^1$, namely the
universal calculus with dimension $m=1$ (any other is a quotient).

In each case we first describe this universal calculus to degree 2 under the
assumption that $\omega=\extd x$ is our basis of $\Omega^1$ from the
left or right. The universal calculus for $n=2$ must have $m=1$ (so no indices related to that) and
we must have the form
\begin{equation*}
d^\mu{}_\nu=
\begin{pmatrix}
0 & 0 \\ 
1 & 0
\end{pmatrix}
,\quad \tau=0.
\end{equation*}
The universal calculus is given by applying $\extd$ to the relations and is always connected, so $f^{\prime }=0$ implies that $f$ is a
multiple of 1, where we recall our notation $\extd f=f^{\prime }\omega$.

(i) (This includes algebra $\rm{A}$ in the table). We have $\extd x.x+ x.
\extd x=0$ so $\{\omega,x\}=0$. Applying $\extd$ again we have $
-\omega\wedge\omega+\omega\wedge\omega=0$ hence there are no relations in
degree 2 for the universal calculus. In terms of structure constants, this is 
\begin{equation*}
a^\mu{}_\nu=
\begin{pmatrix}
1 & 0 \\ 
0 & -1
\end{pmatrix}
,\quad \tau=0
\end{equation*}
which one can check is the most general solution for our structure constants 
$V^{11}{}_0=\lambda$ and $V^{11}{}_1=0$ given that we fixed the form of $d$.

We take $\Omega^2=A.\Vol$ where  $\Vol=\omega\wedge\omega$ is
central. In this case a general central non-zero metric is $
\cg=g\omega\otimes\omega$ for any invertible $g\in A$, and can never be
quantum symmetric. A general connection is $\nabla \omega=\Gamma
\omega\otimes\omega$ for $\Gamma\in A$ and can never be torsion free unless $
\Gamma=0$, in which case $\nabla(x\omega)=\omega\otimes\omega$. This is a
bimodule connection if $\nabla(\omega x)=\sigma(\omega\otimes\omega)$ which
given the commutation relations means $\sigma(\omega\otimes\omega)=-\omega
\otimes\omega$. One can see this also from the bimodule conditions in Lemma~
\ref{lemma2.1}. Then metric compatibility 
\begin{equation*}
0=\nabla \cg=g^{\prime }\omega\otimes\omega\otimes\omega
\end{equation*}
needs $g$ to be a multiple of 1. Hence up to normalisation only $
\cg=\omega\otimes\omega$ admits a QLC. If we look only for a WQLC then we have 
\begin{equation*}
(\extd \otimes\mathrm{id}-\mathrm{id}\wedge\nabla)(g\omega\otimes
\omega)=g^{\prime }\omega\wedge\omega\otimes\omega=0
\end{equation*}
giving the same conclusion of $\cg=\omega\otimes\omega$ for a WQLC. If we
relaxed $T_\nabla=0$ then $coT_\nabla=0$ alone is $g^{\prime }=g\bar\Gamma$
where $\bar{\ }:A\to A$ is the automorphism $\bar x=-x$. Over ${\mathbb{F}}_2
$ this admits $g=1, \Gamma=0$ again and also $g=1+x=\Gamma$.

{\em Alternatively, we can quotient to $\Omega^2=0$}. In this case any connection
is (trivially) a WQLC for any metric and any metric is quantum symmetric. In
this case we do not need to set $\Gamma=0$ and in this case the condition
for $\sigma$ by the same steps as above is $\nabla(x\omega)=\omega\otimes
\omega+x\Gamma\omega\otimes\omega=-\nabla(\omega
x)=-\Gamma\omega\otimes\omega x-\sigma(\omega\otimes\omega)$ which now gives 
$\sigma(\omega\otimes\omega)=-(1+2x\Gamma)\omega\otimes\omega$. This time
the metric compatibility condition is 
\begin{equation*}
0=\nabla \cg=g^{\prime
}\omega\otimes\omega\otimes\omega+g\Gamma\omega\otimes\omega\otimes\omega+
\sigma(g\omega\otimes \Gamma\omega)\otimes\omega
\end{equation*}
which simplifies to 
\begin{equation*}
g^{\prime }+g\Gamma-g\bar \Gamma(1+2x\Gamma)=0
\end{equation*}
in the algebra. Over ${\mathbb{F}}_2$, this simplifies further to $g^{\prime 2}=0$
which for $\lambda=0$ has only $g=1$ and $\Gamma=0,x$ as solutions. This is
also in the row for algebra A in Table~1 and applies for any $\Omega^2$.

(ii) (This includes algebras $\rm{B,C}$ in the table.) Now we have $\extd
x.x+ x.\extd x=\extd x$ so $\{\omega,x\}=\omega$. Applying $
\extd$, we again have $-\omega\wedge\omega+\omega\wedge\omega=0$ hence
there are no relations in degree 2 for the universal calculus. In terms of
structure constants, this is 
\begin{equation*}
a^\mu{}_\nu=
\begin{pmatrix}
1 & 0 \\ 
1 & -1
\end{pmatrix}
,\quad \tau=0
\end{equation*}
which one can check is the most general solution for our structure constants 
$V^{11}{}_0=\lambda$ and $V^{11}{}_1=1$ given that we fixed the form of $\extd$.

We again can take $\Omega^2=A.\Vol$ where $\Vol
=\omega\wedge\omega$ now obeys $\Vol x=-\omega\wedge x\omega+\mathrm{
Vol}=x\omega\wedge \omega - \Vol+\Vol=x\Vol$, so this
is again central. Similarly, a general central non-zero metric is again of
the form $\cg=g\omega\otimes\omega$ for any invertible $g\in A$ and can never
be quantum symmetric. As before, a general connection is $\nabla
\omega=\Gamma \omega\otimes\omega$ for $\Gamma\in A$ and can never be
torsion free unless $\Gamma=0$, in which case a bimodule connection requires 
$\nabla(x\omega)=\omega\otimes\omega=-\nabla(\omega
x)+\nabla\omega=-\sigma(\omega\otimes\omega)$ which gives $
\sigma(\omega\otimes\omega)=-\omega\otimes\omega$ as before. Then metric
compatibility needs $g$ constant hence up to normalisation only $
\cg=\omega\otimes\omega$ admits a QLC just as before. Also as before, only $
\cg=\omega\otimes\omega$  to have a WQLC (these steps are identical the the
previous case). If we drop $T_\nabla=0$ and ask only for $coT_\nabla=0$ then
we need $g^{\prime }=g\bar\Gamma$ again, where now $\bar{\ }:A\rightarrow A$
is the automorphism $\bar{x}=1-x$. Over ${\mathbb{F}}_2$ with $\lambda=0$, 
this has only $g=1$ and $\Gamma=0$ again as solutions. With $\lambda=1$ we
have this and also $\Gamma=g$ for $g=x,1+x$ as shown in the table for algebras $\rm{
B,C}$.

{\em Alternatively, we can quotient to $\Omega ^{2}=0$}. As before, in this case any
connection is (trivially) a WQLC for any metric and any metric is quantum
symmetric. The condition for $\sigma $ by the same steps as above is $\nabla
(x\omega )=\omega \otimes \omega +x\Gamma \omega \otimes \omega =-\nabla
(\omega x)+\nabla \omega =-\Gamma \omega \otimes \omega x-\sigma (\omega
\otimes \omega )+\Gamma \omega \otimes \omega $ which now gives $\sigma
(\omega \otimes \omega )=-(1+(2x-1)\Gamma )\omega \otimes \omega $. This
time the metric compatibility condition is 
\begin{equation*}
0=\nabla \cg=g^{\prime }\omega \otimes \omega \otimes \omega +g\Gamma \omega
\otimes \omega \otimes \omega +\sigma (g\omega \otimes \Gamma \omega
)\otimes \omega
\end{equation*}
which simplifies to 
\begin{equation*}
g^{\prime }+g\Gamma -g\bar{\Gamma}(1+(2x-1)\Gamma )=0
\end{equation*}
in the algebra. Over ${\mathbb{F}}_{2}$, this simplifies to 
\begin{equation*}
g^{\prime }=g(\Gamma+\bar\Gamma+\Gamma\bar\Gamma)
\end{equation*}
with only $g=1$ and $\Gamma=0$ as solutions. This is shown also in the table
for rows $\rm{B,C}$ as it applies for any $\Omega^2$. In geometric terms, the algebra
A is the group algebra ${\mathbb{F}}_2{\mathbb{Z}}_2$ with $z=1+x$ obeying $
z^2=1$ and the algebra B is the function algebra ${\mathbb{F}}_2({\mathbb{Z}}_2)$
with $x,1+x$ the delta-functions at the two points (and the calculi are
covariant with respect to the Hopf algebra structure). The algebra C is ${
\mathbb{F}}_2[x]/(x^2+x+1)$ which is isomorphic to the field ${\mathbb{F}}_4$
as an extension of ${\mathbb{F}}_2$.

\begin{remark} Over $\F_2$, one can check for all the algebras A,B,C  that the only solutions to (\ref{epsrel1})-(\ref{dext2m1}) have
$\tau=0$, so $\omega=\extd x$ as we assumed in our analysis, and that $\eps$ is unique up to multiplication by an invertible functions so that without loss of generality we can take $\eps=1$, or $\Vol=\omega\wedge\omega$, as we also assumed. 
\end{remark}

\begin{table}[tbp]
{\footnotesize {\ 
\begin{tabular}{|l|l|l|l|l|}
\hline
& $
\begin{array}{c}
\text{{Relations}} \\ 
\end{array}
$ & $
\begin{array}{c}
\extd1=0 \\ 
\extd x=\omega
\end{array}
$ & $
\begin{array}{c}
\text{metrics} \\ 
\end{array}
$ & Connections (only $T_\nabla=0$: $\nabla\omega=0$) \\ \hline
$
\begin{array}{c}
\text{A.} \\ 
\\ 
{\mathbb{F}}_2{\mathbb{Z}}_2
\end{array}
$ & $x^{2}=0$ & $\quad \omega .x=x.\omega $ & ${\quad }g_{A.1}=\omega
\otimes \omega $ & $
\begin{array}{l}
g_{A.1}\text{\ compatible (no other $coT_\nabla=0$): } \\ 
\quad \nabla_{A.1.1}\omega=0 \\ 
\quad \nabla_{A.1.2}\omega=x\omega\otimes\omega
\end{array}
$ \\ 
&  &  & $g_{A.2}=\left( 1+x\right) \omega \otimes \omega $ & $
\begin{array}{l}
\text{no $g_{A.2}$ compatible but } coT_\nabla=0: \\ 
\quad \nabla _{A.2}\omega =\left( 1+x\right) \omega \otimes \omega
\end{array}
$ \\ \hline
$
\begin{array}{c}
\text{B.} \\ 
{\mathbb{F}}_2({\mathbb{Z}}_2)
\end{array}
$ & $x^{2}=x,$ & $\omega .x=\omega +x.\omega $ & $
\begin{array}{c}
g_{B}=\omega \otimes \omega
\end{array}
$ & $
\begin{array}{c}
\text{metric compatible (no other $coT_\nabla=0$): } \\ 
\nabla _{B}\omega =0
\end{array}
$ \\ \hline
$
\begin{array}{c}
\text{C.} \\ 
\\ 
{\mathbb{F}}_4
\end{array}
$ & $x^{2}=1+x$ & $\omega .x=\omega +x.\omega $ & $
\begin{array}{c}
g_{C.1}=\omega \otimes \omega \\ 
g_{C.2}=x\omega\otimes\omega \\ 
g_{C.3}=(1+x)\omega\otimes\omega
\end{array}
$ & $
\begin{array}{l}
g_{C.1}\text{\ compatible (no other $coT_\nabla=0$):\ } \\ 
\quad \nabla_{C.1.1}\omega=0 \\ 
\quad \nabla_{C.1.2}\omega=x\omega\otimes\omega \\ 
\quad \nabla_{C.1.3}\omega=(1+x)\omega\otimes\omega \\ 
\text{no $g_{C.2}$ compatible but $coT_\nabla=0$:} \\ 
\quad \nabla _{C.2}\omega= x\omega \otimes \omega \\ 
\text{no $g_{C.3}$ compatible but $coT_\nabla=0$:} \\ 
\quad \nabla _{C.3}\omega=(1+ x)\omega \otimes \omega \\ 
\end{array}
$ \\ \hline
\end{tabular}
} }
\caption{Classification for $n=2,m=1$ and $\Omega^2$ one dimensional. The
central metrics are never quantum symmetric and only $\cg=\protect\omega\otimes
\protect\omega$ admits a QLC and it is $\protect\nabla\protect\omega=0$.
(For $\Omega^2=0$ all metrics are quantum symmetric, all connections are
WQLC and the metric compatible ones in the table are the QLCs.)\label{t1}}
\end{table}

\subsection{Classification geometries for $n=3$ and $m=1$}

In this case there are many more algebras and we restrict to ${\mathbb{F}}_2$
. Then by the results in \cite{MaPac} there are 6 unital commutative
algebras A -- F up to isomorphism and we consider each in turn followed by a further noncommutative one G which turns
up by the same method when we drop commutativity. For differential structures, we
first list the universal one with $m=2$ (the geometry of which we
consider later) and then use a computer to find all possible 1-dimensional
quotients. The algebra relations and the universal calculus relations by
applying $\extd$ to them are

A) $x^{2}=0=y^{2},\ xy=0,$

\begin{equation*}
\extd x.y=x\extd y,\quad \extd y.x=y\extd x,\quad \lbrack 
\extd x,x]=[\extd y,y]=0;
\end{equation*}

B) (this is ${\mathbb{F}}_2({\mathbb{Z}}_3)$ or functions on a triangle) $
x^{2}=x,\ y^{2}=y,\ xy=0$

\begin{equation*}
\extd x.y=x\extd y,\quad \extd y.x=y\extd x,\quad \lbrack 
\extd x,x]=\extd x,\quad \lbrack \extd y,y]=\extd y;
\end{equation*}

C) $x^{2}=x,\ y^{2}=xy=0$

\begin{equation*}
\extd x.y=x\extd y,\quad \extd y.x=y\extd x,\quad \lbrack 
\extd x,x]=\extd x,\quad \lbrack \extd y,y]=0;
\end{equation*}

D) (this is ${\mathbb{F}}_2{\mathbb{Z}}_3$) $x^{2}=y,\ y^{2}=x,\ xy=x+y$

\begin{equation*}
\extd x.y=x\extd y+\extd x+\mathrm{d }y,\quad \extd y.x=y
\extd x+\extd x+\extd y,\quad \lbrack \extd x,x]=\extd
y,\quad \lbrack \extd y,y]=\extd x;
\end{equation*}

E) (this is an anyonic line ${\mathbb{F}}_2[x]/(x^2)$) $x^{2}=y,\ y^{2}=xy=0$

\begin{equation*}
\extd x.y=x\extd y,\quad \extd y.x=y\extd x,\quad \lbrack 
\extd x,x]=\extd y,\quad \lbrack \extd y,y]=0;
\end{equation*}

F) (this is the field extension ${\mathbb{F}}_8=\mathbb{F}_2[y]/(y^3+y^2+1)$)\  $y^{2}=x,\ xy=1+x$  (and $x^{2}=1+x+y$ implied)

\begin{equation*}
\extd x.y=x\extd y+\extd x,\quad \extd y.x=y\extd x+
\extd x,\quad \lbrack \extd x,x]= \extd x+\extd y,\quad
\lbrack \extd y,y]=\extd x;
\end{equation*}

G) (this is noncommutative) $\F_2\<x,y\>$ modulo the ideal generated by the relations $x^2=x$, $y^2=0$, $xy=y$, $yx=0$

\begin{equation*}
\extd x.x= (1+x)\extd x,\quad \extd x.y=(1+x)\extd y,\quad \extd y.x= y\extd x,\quad \extd y.y=y\extd y.
\end{equation*}

For $m=1$, we have in each case to add a relation to the ones coming from the
universal calculus. The results are as follows.

(i) For algebras A, D, and E there are no solutions for 1-dimensional
differential calculi $\Omega^1$ of the left-parallelisable form assumed in our
general analysis.

(ii) Algebra F has potentially 14 differential calculi $\Omega^1$ by solving our equations:
\begin{equation*}
\begin{array}{lll}
{\rm F.1} & \quad\extd x=\omega ,\ \extd y=y.\omega,\ & \
\omega .x=\omega +x.\omega +y.\omega ,\ \omega .y=x.\omega ,\ ~ \\ 
{\rm F.2} & \quad\extd x=y.\omega ,\ \extd y=x.\omega,\  & \ \mathrm{
bimodule\ relations\ as\ in\ } F.1 \\ 
{\rm F.3} & \quad\extd x=x.\omega ,\ \extd y=\omega +x.\omega,\  & \ 
\mathrm{bimodule\ relations\ as\ in\ } F.1 \\ 
{\rm F.4} & \quad \extd x=x.\omega +y.\omega ,\ \extd y=\omega,\   & \ \mathrm{
bimodule\ relations\ as\ in\ } F.1 \\ 
{\rm F.5} & \quad\extd x=\omega +x.\omega +y.\omega ,\ \extd y=\omega
+y.\omega,\  & \ \mathrm{bimodule\ relations\ as\ in\ } F.1 \\ 
{\rm F.6} & \quad\extd x=\omega +y.\omega ,\ \extd y=x.\omega
+y.\omega,\  & \ \mathrm{bimodule\ relations\ as\ in\ } F.1 \\ 
{\rm F.7} & \quad\extd x=\omega +x.\omega ,\ \extd y=\omega +x.\omega
+y.\omega,\  & \ \mathrm{bimodule\ relations\ as\ in\ } F.1 \\ 
{\rm F.8} & \quad\extd x=\omega+x.\omega ,\ \extd y=\omega,\  & \
\omega .x=y.\omega ,\ \omega .y=\omega +x.\omega +y.\omega,\  \\ 
{\rm F.9} & \quad\extd x=\omega +y.\omega ,\ \extd y=x.\omega,\ 
  & \ \mathrm{bimodule\ relations\ as\ in\ } F.8 \\ 
{\rm F.10} & \quad\extd x=x.\omega +y.\omega ,\ \extd y=\omega
+x.\omega,\  & \ \mathrm{bimodule\ relations\ as\ in\ } F.8 \\ 
{\rm F.11} & \quad\extd x=\omega +x.\omega +y.\omega ,\ \extd
y=y.\omega,\  & \ \mathrm{bimodule\ relations\ as\ in\ } F.8 \\ 
{\rm F.12} & \quad\extd x=y.\omega ,\ \extd y=\omega +y.\omega,\  & \ 
\mathrm{bimodule\ relations\ as\ in\ } F.8 \\ 
{\rm F.13} & \quad\extd x=x.\omega ,\ \extd y=x.\omega +y.\omega,\  & \ 
\mathrm{bimodule\ relations\ as\ in\ } F.8 \\ 
{\rm F.14} & \quad\extd x=\omega ,\ \extd y=\omega +x.\omega +y.\omega
,\  & \ \mathrm{bimodule\ relations\ as\ in\ } F.8
\end{array}
\end{equation*}

However, the first 7 are all isomorphic under a change of basis $\phi(\omega)=y^i\omega$ for some power $i$ and likewise for the last 7, so there are only two non-isomorphic calculi. Since $\Bbb F_8$ is a field, multiplication by any nonzero element is an isomorphism and there are exactly 7 in each isomorphism class. One the other hand, all the calculi have invertible $a^{\mu }{}_{\nu }$ matrices and none of them admits a central non-zero metric. One is also forced to $\Omega^2=0$.

(iii) Algebra B has 8 left-parallelisable differential calculi by solving our equations, and they are all distinct since the algebra has only 1 as invertible element. Of them,  only B.4 and B.8 have invertible $a^{\mu }{}_{\nu }$ needed for 
$\omega $ to be both a left and a right basis (our preferred case).  However,
it is exactly these more geometrical ones which admit no non-zero central
metric even when we relax invertibility. This is summarised in Table~2. None
of the non-invertible metrics admit a metric compatible connection
either. Also for these calculi, $\Omega^2$ is forced to have vector space dimension less than that of the algebra so we either have
to take $\Omega^2=0$ or it is not a free module with
a single basis element $\Vol$, which is an added complication as the analysis in Section~\ref{secpre} won't apply. For example, for calculus $\rm{B.1}$ one can apply $\extd$ to the first order relations to conclude that for any exterior algebra $\extd x\wedge\extd y=\extd y\wedge\extd x=\extd y\wedge\extd y=0$ giving only $\omega\wedge\omega, x\omega\wedge\omega$ as the 2-forms with $\Vol=\omega\wedge\omega$ as a generator but not a basis over the algebra (some products with it are zero). In some cases this generator is also not central. 

\begin{table}[tbp]
\begin{tabular}{|l|l|l|}
\hline
Family B & Differential calculi relations & Central metrics (non-invertible)
\\ \hline
$\text{B.1}$ & $
\begin{array}{c}
\extd x=\omega ,\ \extd y=y.\omega \\ 
\omega .x=\omega +x.\omega ,\ \omega .y=0 \\ 
\end{array}
$ & $
\begin{array}{c}
g_{B.1}=\left( 1+y\right) \omega \otimes \omega \\ 
g_{B.1.1}=xg_{B.1}=x\omega \otimes \omega {\ } \\ 
g_{B.1.2}=\text{ }(1+x)g_{B.1}=\left( 1+x+y\right) \omega \otimes \omega \\ 
\end{array}
$ \\ \hline
$\text{B.2}$ & $
\begin{array}{c}
\hline
\extd x=\omega +y.\omega ,\ \extd y=y.\omega \\ 
\omega .x=\omega +x.\omega +y.\omega ,\ \omega .y=0 \\ 
\end{array}
$ & $
\begin{array}{c}
\text{metrics the same as in B.1}
\end{array}
$ \\ \hline
B.3 & $
\begin{array}{c}
\extd x=\omega ,\ \extd y=x.\omega +y.\omega ,\  \\ 
\omega .x=\omega +x.\omega ,\ \omega .y=x.\omega \\ 
\end{array}
$ & $
\begin{array}{c}
g_{B.3}=\left( x+y\right) \omega \otimes \omega \\ 
g_{B.3.1}=xg_{B.3}=x\omega \otimes \omega \\ 
g_{B.3.2}=yg_{B.3}=y\omega \otimes \omega
\end{array}
{\ }$ \\ \hline
B.4 & $
\begin{array}{c}
\extd x=\omega +y.\omega ,\ \extd y=x.\omega +y.\omega ,\  \\ 
\omega .x=\omega +x.\omega +y.\omega ,\ \omega .y=x.\omega \\ 
\end{array}
$ & $\text{{no central metrics}}$ \\ \hline
$\text{B.5}$ & $\text{equivalent to B.1 under }x\leftrightarrow y$ &  \\ 
\hline
$\text{B.6}$ & $\text{equivalent to B.3 under }x\leftrightarrow y$ &  \\ 
\hline
$\text{B.7}$ & $\text{{equivalent to B.2} under }x\leftrightarrow y$ &  \\ 
\hline
$\text{B.8}$ & $\text{{equivalent to B.4} under }x\leftrightarrow y$ &  \\ 
\hline
\end{tabular}
\caption{1-dimensional differential calculi for the algebra B with the
corresponding central metrics, but none invertible. None of them admit a
metric compatible connection (or a cotorsion free or torsion free one unless $\Omega^2=0$).\label{t2}}
\end{table}

(iv) The similar  computer results for the algebra C at $m=1$ are summarised in Table~3 and again they are all distinct as only 1 is invertible in the algebra. Again we see that there are some
non-invertible central metrics. 
Also, in all cases $a^{\mu }{}_{\nu }$ is not invertible so these are not fully parallelisable 
in our sense and we have similar issues that we should take $\Omega^2=0$ or it is not a free module over the algebra even though we can take $\Vol=\omega\wedge\omega$ as the generator, in some cases not central. On the other hand, some of the metrics do admit compatible connections. They
all have torsion in the case of $\Omega^2\ne 0$ with $T_{\nabla }\omega=\Gamma\Vol$ for any connection
$\Gamma=n_{0}+n_{1}x+n_{2}y$, say. Note that because the torsion is not zero, being metric compatible does not imply cotorsion free and in fact for $g_{C.1.2}$ the cotorsion free connections are a subset of the metric compatible ones as shown in the table. The same for the $\rm{C.2}$ calculus. The curvatures for the four calculi are:

 C.1.: $R_{\nabla }\omega=(n_0+n_1+n_0 n_{1}+(1+n_0+n_1)n_2 y)
\Vol\otimes \omega$. For the six metric compatible ones in the table the curvature is always $R_{\nabla }\omega=\Vol\otimes \omega$ 
except for $\Gamma=1+x+y$ which has a factor $1+y$ out front. 

C.2.: $R_{\nabla}(\omega)=(n_{0}+n_{1}+n_{0}n_{1}+(n_{1}+(1+n_{0}+n_1)n_{2})y)\Vol\otimes \omega $. For the 6 metric compatible ones we have the same options for the curvature, with the $1+y$ factor when $\Gamma=x,x+y,1+x$.                           

C.3.: $R_{\nabla }\omega=(n_{0}+n_{2}+(n_{1}+n_{0}n_{1}+n_{2})x)\Vol
\otimes \omega $

C.4.: $R_{\nabla }\omega=(n_{0}+n_{2}+(n_{1}+n_{0}n_{1}+n_{2})x+n_{2}y)\Vol\otimes \omega $.

\begin{table}[tbp]
{\footnotesize {\ 
\begin{tabular}{|l|l|l|l|}
\hline
Family C & $
\begin{array}{c}
\text{Differential } \\ 
\text{calculi relations}
\end{array}
$ & $\text{ Central metrics (non-invertible)}$ & Connections (with torsion) \\ \hline
C.1 & $
\begin{array}{c}
\extd x=\omega , \\ 
\extd y=y.\omega , \\ 
\omega .x=\omega +x.\omega , \\ 
\omega .y=0
\end{array}
$ & $
\begin{array}{c}
g_{C.1}=\left( x+y\right) \omega \otimes \omega \\ 
g_{C.1.1}=xg_{C.1}=x\omega \otimes \omega \\ 
g_{C.1.2}=\left( 1+x+y\right) g_{C.1}=y\omega \otimes \omega
\end{array}
$ & $
\begin{array}{l}
\text{only $g_{C.1.2}$ metric compatible}\\
\quad \nabla _{C.1.2.1}-\nabla_{C.1.2.6}: \ \Gamma\ne 0,y\\
\text{($coT_\nabla=0$ if also $\Gamma\ne 1+x,1+x+y$)}
\end{array}
$ \\ \hline
C.2 & $
\begin{array}{c}
\extd x=\omega +y.\omega , \\ 
\extd y=y.\omega , \\ 
\omega .x=\omega +x.\omega , \\ 
\omega .y=0
\end{array}
$ & metrics the same as for C.1 & $
\begin{array}{c}
\text{connections the same as for C.1}
\end{array}
$ \\ \hline
C.3 & $
\begin{array}{c}
\extd x=x.\omega , \\ 
\extd y=\omega +x.\omega , \\ 
\omega .x=0, \\ 
\omega .y=y.\omega
\end{array}
$ & $
\begin{array}{c}
g_{C.3}=\left( 1+x\right) \omega \otimes \omega \\ 
g_{C.3.1}=\left( 1+y\right) g_{C.3}\\ \quad\qquad\qquad=\left( 1+x+y\right) \omega \otimes \omega
\\ 
g_{C.3.1}=\left( x+y\right) g_{C.3}=y\omega \otimes \omega
\end{array}
$ & none metric comp. or $coT_\nabla=0$ \\ \hline
C.4 & $
\begin{array}{c}
\extd x=x.\omega , \\ 
\extd y=\omega +x.\omega +y.\omega , \\ 
\omega .x=0, \\ 
\omega .y=y.\omega
\end{array}
$ & metrics the same as for C.3 & none metric comp. or $coT_\nabla=0$  \\ \hline
\end{tabular}
}}
\caption{1-dimensional differential calculi for the algebra C with the
corresponding central but non-invertible metrics. None of the metrics are
quantum symmetric and none of the metric compatible connections have $T_
\protect\nabla=0$ for $\Omega^2\ne 0$. (For $\Omega^2=0$ all
metrics are quantum symmetric, all connections are WQLC and the metric
compatible ones in the table are the QLCs.)\label{t3}}
\end{table}

(v) For the noncommutative algebra G there are 48 parallelisable calculi with $m=1$ and only 1 is invertible in the algebra so they are distinct in the sense of a change of $\omega$. But none of them have a  central metric, so we do not discuss them individually. 

Summarising, for the $n=3,m=1$ case there are 6 possible commutative algebras and one noncommutative one. None of them admits a calculus $\Omega^1$ of dimension $m=1$ satisfying all our requirements. Some algebras (B and C) admit a central candidate for a metric but none of these are invertible.  These algebras also have issues with $\Omega^2\neq 0$. The noncommutative algebra G, for $m=1$, does not admit a suitable calculus having a central metric.

\subsection{Classification of $n=3,m=2$ geometries on the algebra D}\label{secDmodel}

We will now consider each of the 6 algebras above but with the $m=2$
case of the universal calculus $\Omega^1$. To keep things simple we
consider geometries with basis $\omega ^{1}=\extd
x,\omega ^{2}=\extd y$ for $\Omega^1$, so that $\tau^i:=\extd \omega^i=0$. The universal calculus at $
\Omega^2$ is normally too large to be reasonable for a geometry -- we will
need to quotient it to obtain something more `reasonable' such as $\Omega^2$
1-dimensional. This has to be searched for on a case-by-case basis for each
algebra. In this section we illustrate the method in detail on the algebra D and then for 
the other algebras we just list the results. 

In fact the algebra D is isomorphic to the group algebra of the group ${
\mathbb{Z}}_{3}$ since $z=1+x$ obeys $z^{2}=1+y$ and $z^{3}=1$ in the
algebra. The bimodule commutation relations in terms of these are
\begin{equation*}
\omega ^{1}.z=\omega ^{1}(x+1)=(x+1)\omega ^{1}+\omega ^{2}=z\omega
^{1}+\omega ^{2}
\end{equation*}
\begin{equation*}
\omega ^{2}.z=\omega ^{2}(x+1)=(y+1)\omega ^{1}+\omega ^{2}+\omega
^{2}=z^{2}\omega ^{1}.
\end{equation*}
From these it is easy to see that 
\begin{eqnarray}  \label{liftZ3}
\cg &=&\alpha z\omega ^{1}\otimes \omega ^{1}+\alpha \omega ^{1}\otimes \omega
^{2}+\beta z\omega ^{2}\otimes \omega ^{1}+\beta \omega ^{2}\otimes \omega
^{2} \\
&=&\alpha \left( 1+x\right) \omega ^{1}\otimes \omega ^{1}+\alpha \omega
^{1}\otimes \omega ^{2}+\beta \left( 1+x\right) \omega ^{2}\otimes \omega
^{1}+\beta \omega ^{2}\otimes \omega ^{2}  \notag
\end{eqnarray}
is the general form of a central element in the tensor square, for any two
functions $\alpha ,\beta $ in the algebra.

Next, applying $\extd$ to these gives no relations with the result that 
$\Omega^2$ for the universal calculus is 4-dimensional. For a natural
1-dimensional $\Omega^2$, we appeal to the group theory where we have
left-invariant 1-forms $e_\pm$ forming a Grassmann algebra. Making the
isomorphism formally (which one can do by making a field extension to adjoin
a cube root of 1) and transferring back, we are led to define 
\begin{equation*}
\omega^1\omega^2=\omega^2\omega^1=0,\quad( \omega^1)^2+z(\omega^2)^2=0
\end{equation*}
giving a 1-dimensional $\Omega^2$, which we take with basis $\Vol
=(\omega^1)^2=z(\omega^2)^2$ say. This is central and has 
\begin{equation*}
\omega^i\omega^j=\epsilon^{ij}\Vol,\quad\epsilon^{11}=1,\quad
\epsilon^{22}=z^2,\quad \epsilon^{12}=\epsilon^{21}=0.
\end{equation*}
One can check that one has a DGA with $\extd\omega^i=0$ (so $\tau^i=0$
).

Then $\wedge (\cg)=0$ requires $\alpha =z\beta $ so we have just a
1-functional parameter of non-zero central quantum symmetric metrics, 
\begin{eqnarray*}
\cg &=&\beta \left( z(z\omega ^{1}+\omega ^{2})\otimes \omega ^{1}+(z\omega
^{1}+\omega ^{2})\otimes \omega ^{2}\right) =\beta \left( (\omega
^{1}z^{2}+\omega ^{2}z)\otimes \omega ^{1}+\omega ^{1}z\otimes \omega
^{2}\right) \\
&=&\beta \left( 1+y\right) \omega ^{1}\otimes \omega ^{1}+\beta \left(
1+x\right) \omega ^{2}\otimes \omega ^{1}+\beta \left( 1+x\right) \omega
^{1}\otimes \omega ^{2}+\beta \omega ^{2}\otimes \omega ^{2}
\end{eqnarray*}
where the latter expression makes it clear that this is invertible at least
when $\beta =1$ because the `internal' coefficient matrix $\tilde g$ is invertible (and
for typical $\beta $ according to how the coefficients look when $\beta $ is
commuted to the middle).

We now use a computer to solve for QLCs with torsion depending on the choice
of $\Omega^2$, which in the 1-dimensional case comes down to the choice of
volume element $\Vol$. Fixing one of these, we  find four QLCs for each choice of $
\beta=1,z,z^2$ of invertible $\beta$, all with curvature except for one flat
one in each group.

1) ($\beta =1$) $g_{D.1}=z(z\omega ^{1}+\omega ^{2})\otimes \omega
^{1}+(z\omega ^{1}+\omega ^{2})\otimes \omega ^{2}$: 
\begin{align*}
\nabla _{D.1.1}\omega ^{1}& =z^{2}\omega ^{1}\otimes \omega^{1}+(1+z)(\omega ^{1}\otimes \omega ^{2}+\omega ^{2}\otimes \omega
^{1})+\omega ^{2}\otimes \omega ^{2} \\
\nabla _{D.1.1}\omega ^{2}& =z^{2}\omega ^{1}\otimes \omega ^{1}+z\omega^{1}\otimes \omega ^{2}+z^{2}\omega ^{2}\otimes \omega ^{1}+\omega^{2}\otimes \omega ^{2} \\
R_{\nabla _{D.1.1}}\omega ^{1}& =z^{2}\mathrm{Vol}\otimes \omega ^{1}+z\mathrm{Vol}\otimes \omega ^{2}~,\quad R_{\nabla _{D.1.1}}\omega^{2}=\mathrm{Vol}\otimes \omega ^{1}~;
\end{align*}
\begin{align*}
\nabla _{D.1.2}\omega ^{1}& =z^{2}\omega ^{1}\otimes \omega ^{1}+z(\omega
^{1}\otimes \omega ^{2}+\omega ^{2}\otimes \omega ^{1})+\omega ^{2}\otimes
\omega ^{2} \\
\nabla _{D.1.2}\omega ^{2}& =z^{2}\omega ^{2}\otimes \omega ^{1} \\
R_{\nabla _{D.1.2}}\omega ^{1}& =R_{\nabla _{D.1.2}}\omega ^{2}=0~;
\end{align*}
\begin{align*}
\nabla _{D.1.3}\omega ^{1}& =(z+z^{2})\omega ^{1}\otimes \omega
^{1}+(1+z)\omega ^{1}\otimes \omega ^{2}+z\omega ^{2}\otimes \omega^{1}+\left( 1+z^{2}\right) \omega ^{2}\otimes \omega ^{2} \\
\nabla _{D.1.3}\omega ^{2}& =z^{2}\omega ^{1}\otimes \omega ^{1}+\left(z+z^{2}\right) \omega ^{2}\otimes \omega ^{1}+\omega ^{2}\otimes \omega ^{2}
\\
R_{\nabla _{D.1.3}}\omega ^{1}& =z^{2}\mathrm{Vol}\otimes \omega ^{1}+z\mathrm{Vol}\otimes \omega ^{2}~\ ,\quad R_{\nabla _{D.1.3}}\omega^{2}=\mathrm{Vol}\otimes \omega ^{1}~;
\end{align*}
\begin{align*}
\nabla _{D.1.4}\omega ^{1}& =(z+z^{2})\omega ^{1}\otimes \omega ^{1}+z\omega^{1}\otimes \omega ^{2}+(1+z)\omega ^{2}\otimes \omega ^{1}+\left(1+z^{2}\right) \omega ^{2}\otimes \omega ^{2} \\
\nabla _{D.1.4}\omega ^{2}& =z\omega ^{1}\otimes \omega ^{2}+\left(z+z^{2}\right) \omega ^{2}\otimes \omega ^{1} \\
R_{\nabla _{D.1.4}}\omega ^{1}& =z^{2}\mathrm{Vol}\otimes \omega ^{1}+z\mathrm{Vol}\otimes \omega ^{2}~\ \ ,\quad R_{\nabla _{D.1.4}}\omega^{2}=\mathrm{Vol}\otimes \omega ^{1}~.
\end{align*}

\medskip 2) ($\beta =z$) $g_{D.2}=z^{2}(z\omega ^{1}+\omega ^{2})\otimes\omega ^{1}+z(z\omega ^{1}+\omega ^{2})\otimes \omega ^{2}$: 
\begin{align*}
\nabla _{D.2.1}\omega ^{1}& =\omega ^{1}\otimes \omega ^{1}+z\omega^{1}\otimes \omega ^{2}+z^{2}\omega ^{2}\otimes \omega ^{1}+z\omega^{2}\otimes \omega ^{2} \\
\nabla _{D.2.1}\omega ^{2}& =\omega ^{1}\otimes \omega ^{1}+(1+z^{2})\left(\omega ^{1}\otimes \omega ^{2}+\omega ^{2}\otimes \omega ^{1}\right)
+z\omega ^{2}\otimes \omega ^{2} \\
R_{\nabla _{D.2.1}}\omega ^{1}& =z^{2}\mathrm{Vol}\otimes \omega ^{2},\quad R_{\nabla _{D.2.1}}\omega ^{2}=z\mathrm{Vol}\otimes \omega^{1}+\mathrm{Vol}\otimes \omega ^{2}~;
\end{align*}%
\begin{align*}
\nabla _{D.2.2}\omega ^{1}& =\omega ^{1}\otimes \omega ^{1}+(z+z^{2})\omega^{1}\otimes \omega ^{2}+z\omega ^{2}\otimes \omega ^{2} \\
\nabla _{D.2.2}\omega ^{2}& =\left( 1+z\right) \omega ^{1}\otimes \omega^{1}+z^{2}\omega ^{1}\otimes \omega ^{2}+\left( 1+z^{2}\right) \omega^{2}\otimes \omega ^{1}+(z+z^{2})\omega ^{2}\otimes \omega ^{2} \\
R_{\nabla _{D.2.2}}\omega ^{1}& =z^{2}\mathrm{Vol}\otimes \omega ^{2}~\
,\quad R_{\nabla _{D.2.2}}\omega ^{2}=z\mathrm{Vol}\otimes \omega ^{1}+\mathrm{Vol}\otimes \omega ^{2}~;
\end{align*}%
\begin{align*}
\nabla _{D.2.3}\omega ^{1}& =(z+z^{2})\omega ^{1}\otimes \omega^{2}+z^{2}\omega ^{2}\otimes \omega ^{1} \\
\nabla _{D.2.3}\omega ^{2}& =\left( 1+z\right) \omega ^{1}\otimes \omega^{1}+(1+z^{2})\omega ^{1}\otimes \omega ^{2}+z^{2}\omega ^{2}\otimes \omega
^{1}+(z+z^{2})\omega ^{2}\otimes \omega ^{2} \\
R_{\nabla _{D.2.3}}\omega ^{1}& =z^{2}\mathrm{Vol}\otimes \omega ^{2},\quad
R_{\nabla _{D.2.3}}\omega ^{2}=z\mathrm{Vol}\otimes \omega ^{1}+\mathrm{Vol}
\otimes \omega ^{2}~;
\end{align*}%
\begin{align*}
\nabla _{D.2.4}\omega ^{1}& =z\omega ^{1}\otimes \omega ^{2} \\
\nabla _{D.2.4}\omega ^{2}& =\omega ^{1}\otimes \omega ^{1}+z^{2}(\omega^{1}\otimes \omega ^{2}+\omega ^{2}\otimes \omega ^{1})+z\omega ^{2}\otimes
\omega ^{2} \\
R_{\nabla _{D.2.4}}\omega ^{1}& =R_{\nabla _{D.2.4}}\omega ^{2}=0.~\
\end{align*}

\bigskip 3) ($\beta =z^{2}$) $g_{D.3}=(z\omega ^{1}+\omega ^{2})\otimes\omega ^{1}+z^{2}(z\omega ^{1}+\omega ^{2})\otimes \omega ^{2}$: 
\begin{equation*}
\nabla _{D.3.1}\omega ^{1}=z\omega ^{2}\otimes \omega ^{1},\quad \nabla_{D.3.1}\omega ^{2}=z^{2}\omega ^{1}\otimes \omega ^{2},\quad R_{\nabla_{D.3.1}}\omega ^{1}=R_{\nabla _{D.3.1}}\omega ^{2}=0;
\end{equation*}
\begin{align*}
\nabla _{D.3.2}\omega ^{1}& =z^{2}\omega ^{1}\otimes \omega ^{1}+\omega^{2}\otimes \omega ^{2},\quad \nabla _{D.3.2}\omega ^{2}=z^{2}\omega^{2}\otimes \omega ^{1} \\
R_{\nabla _{D.3.2}}\omega ^{1}& =z\mathrm{Vol}\otimes \omega ^{1},\quad
R_{\nabla _{D.3.2}}\omega ^{2}=z\mathrm{Vol}\otimes \omega ^{2};~\
\end{align*}
\begin{align*}
\nabla _{D.3.3}\omega ^{1}& =z\omega ^{1}\otimes \omega ^{2},\quad \nabla_{D.3.3}\omega ^{2}=\omega ^{1}\otimes \omega ^{1}+z\omega ^{2}\otimes\omega ^{2} \\
R_{\nabla _{D.3.3}}\omega ^{1}& =z\mathrm{Vol}\otimes \omega ^{1},\quad
R_{\nabla _{D.3.3}}\omega ^{2}=z\mathrm{Vol}\otimes \omega ^{2};~\ 
\end{align*}%
\begin{align*}
\nabla _{D.3.4}\omega ^{1}& =z^{2}\omega ^{1}\otimes \omega ^{1}+z\left(\omega ^{1}\otimes \omega ^{2}+\omega ^{2}\otimes \omega ^{1}\right) +\omega
^{2}\otimes \omega ^{2} \\
\nabla _{D.3.4}\omega ^{2}& =\omega ^{1}\otimes \omega ^{1}+z^{2}\left(\omega ^{1}\otimes \omega ^{2}+\omega ^{2}\otimes \omega ^{1}\right)
+z\omega ^{2}\otimes \omega ^{2} \\
R_{\nabla _{D.3.4}}\omega ^{1}& =z\mathrm{Vol}\otimes \omega ^{1},\quad
R_{\nabla _{D.3.4}}\omega ^{2}=z\mathrm{Vol}\otimes \omega ^{2}.~\ 
\end{align*}

For this family there are in fact 3 possible (but equivalent) solutions for $
\epsilon $ (that are invertible and satisfy  (\ref{epsrel}), 
\eqref{dext1},\eqref{dext2},
\eqref{qs} for a central volume form with metric quantum symmetric). They are multiples by the invertible functions of
\[ \epsilon=\left( 
\begin{array}{cc}
1 & 0 \\ 
0 & z^{2}
\end{array}
\right) ,\quad \epsilon^{-1}=\left( 
\begin{array}{cc}
1 & 0 \\ 
0 & z
\end{array}
\right)\] which is the one used above and corresponds to $\Vol=(\omega ^{1})^{2}.$ The other solutions
are $z,z^2$ times this which one can view as $z^2,z$ times the volume form for the same calculus.

\subsection{ Classification of $n=3,m=2$ geometries on the algebra B}\label{secBmodel}

In this section we keep the $m=2$ or universal calculus with $\omega ^{1}=
\extd x$ and $\omega ^{2}=\extd y$ and $\Omega^2$ again
one-dimensional so that we can use the same general set-up as above, but we
consider the other possible algebras from our list. First we find by computer that:

\begin{lemma}
For algebras {\rm{A,C,E,G}}  there are no invertible central metrics for the $m=2$ calculus. 
\end{lemma}
For example, for the noncommutative algebra G, one finds from the commutation relations that there is a unique central element 
\[ g_G= (1+x)\extd y\tens \extd y+ y\extd y\tens \extd x=\extd x . y\tens\extd y+ \extd y.y\tens\extd x\]
so that $\tilde g_G=y\sigma_1$ (a Pauli matrix), which is not invertible since $y$ is not. 

That leaves algebras B,F and we consider B in this section.  There is one invertible quantum symmetric central metric
\begin{equation*}
g_{B}=(1+y)\omega ^{1}\otimes \omega ^{1}+(1+x+y)(\omega ^{1}\otimes \omega
^{2}+\omega ^{2}\otimes \omega^{1})+(1+x)\omega ^{2}\otimes \omega ^{2}. 
\end{equation*}

There is only one solution for $\epsilon$ (which is invertible and satisfies 
 (\ref{epsrel}),\eqref{dext1},\eqref{dext2},\eqref{qs} for a central volume form with metric quantum symmetric), namely 
\begin{equation*}
\epsilon =\left( 
\begin{array}{cc}
1+x & x+y \\ 
x+y & 1+y
\end{array}
\right) ,\quad \quad \epsilon ^{-1}=\left( 
\begin{array}{cc}
1+y & x+y \\ 
x+y & 1+x
\end{array}
\right).
\end{equation*}

There are 4 QLCs for this family and the metric above, of which three are 
flat: 
\begin{align*}
\nabla _{B.1}\omega ^{1}& =(1+x+y)\omega ^{1}\otimes \omega ^{1}+(1+y)\omega
^{1}\otimes \omega ^{2}+\omega ^{2}\otimes \omega ^{2} \\
\nabla _{B.1}\omega ^{2}& =\omega ^{1}\otimes \omega ^{1}+\omega ^{1}\otimes
\omega ^{2}+(1+x+y\bigskip )\omega ^{2}\otimes \omega ^{1}+(1+y)\omega
^{2}\otimes \omega ^{2} \\
R_{\nabla _{B.1}}\omega ^{1}& = R_{\nabla _{B.1}}\omega ^{2}=0;
\end{align*}%
\begin{align*}
\nabla _{B.2}\omega ^{1}& =(1+x)\omega ^{1}\otimes \omega ^{1}+(1+x+y)\omega
^{1}\otimes \omega ^{2}+\omega ^{2}\otimes \omega ^{1}+\omega ^{2}\otimes
\omega ^{2} \\
\nabla _{B.2}\omega ^{2}& =\omega ^{1}\otimes \omega ^{1}+\left( 1+x\right)
\omega ^{2}\otimes \omega ^{1}+(1+x+y)\omega ^{2}\otimes \omega ^{2} \\
R_{\nabla _{B.2}}\omega ^{1}& = R_{\nabla _{B.2}}\omega ^{2}=0;
\end{align*}%
\begin{align*}
\nabla _{B.3}\omega ^{1}& =(1+x)\omega ^{1}\otimes \omega ^{1}+(1+y)\omega
^{1}\otimes \omega ^{2}+\left( x+y\right) \omega ^{2}\otimes \omega
^{1}+\left( 1+x\right) \omega ^{2}\otimes \omega ^{2} \\
\nabla _{B.3}\omega ^{2}& =\left( 1+y\right) \omega ^{1}\otimes \omega
^{1}+\left( x+y\right) \omega ^{1}\otimes \omega ^{2}+\left( 1+x\right)
\omega ^{2}\otimes \omega ^{1}+(1+y)\omega ^{2}\otimes \omega ^{2} \\
R_{\nabla _{B.3}}\omega ^{1}& =R_{\nabla _{B.3}}\omega ^{2}=0;
\end{align*}%
\begin{align*}
\nabla _{B.4}\omega ^{1}& =(1+x+y)\left( \omega ^{1}\otimes \omega
^{1}+\omega ^{1}\otimes \omega ^{2}+\omega ^{2}\otimes \omega ^{1}\right)
+\left( 1+x\right) \omega ^{2}\otimes \omega ^{2} \\
\nabla _{B.4}\omega ^{2}& =\left( 1+y\right) \omega ^{1}\otimes \omega
^{1}+\left( 1+x+y\right) \left( \omega ^{1}\otimes \omega ^{2}+\omega
^{2}\otimes \omega ^{1}+\omega ^{2}\otimes \omega ^{2}\right)  \\
R_{\nabla _{B.4}}\omega ^{1}& =\mathrm{Vol}\otimes \omega ^{1}~, 
\\
R_{\nabla _{B.4}}\omega ^{2}& =\mathrm{Vol}\otimes \omega ^{2}.
\end{align*}

In fact the algebra here is ${\mathbb{F}}_{2}({\mathbb{Z}}_{3})$ with $
x=\delta _{1}$, $y=\delta _{2}$ and $1+x+y=\delta _{0}$ for the three
delta-functions. We can identify left-invariant 1-forms with respect to the
group structure, 
\begin{equation*}
e^{1}=(x+1)\omega ^{1}+(x+y)\omega ^{2},\quad e^{2}=(x+y)\omega
^{1}+(y+1)\omega ^{2}
\end{equation*}
which provide a manageable route to solving these equations by hand, with
the same results as above and simple commutation relations $e^{1}\delta
_{2}=\delta _{1}e^{1}$, $e^{2}\delta _{1}=\delta _{2}e^{2}$ etc (the
standard triangle graph calculus) and Grassmann algebra $(e^{i})^{2}=0$, $
e^{1}e^{2}=e^{2}e^{1}=\Vol$. From this point of view the metric is the Euclidean
metric $g_{B}=e^{1}\otimes e^{2}+e^{2}\otimes e^{1}$ which over ${\mathbb{F}}
_{2}$ is unique as 1 is the only invertible function.

\subsection{ Classification of $n=3,m=2$ geometries on the algebra F}\label{secFmodel}

For the algebra F with its $m=2$ universal $\Omega^1$, there are 7 invertible quantum symmetric central metrics namely any nonzero multiple of any one of them, e.g.  
\begin{eqnarray*}
\cg&=&\beta \left( y^2\omega^1\tens\omega^1+\omega^1\tens \omega^2+\omega^2\tens\omega^1+(1+y)\omega^2\tens\omega^2\right)
\end{eqnarray*}
for any nonzero $\beta$ (necessarily invertible since the algebra here is a field). However, only four of them admit QLCs:

1) ($\beta=y^2$) $g_{F.1}=(1+y+y^2)\omega ^{1}\otimes \omega
^{1}+y^2(\omega ^{1}\otimes \omega ^{2}+\omega ^{2}\otimes \omega ^{1})+\omega
^{2}\otimes \omega ^{2}$

with 12 QLCs, five of which are flat:
\begin{eqnarray*}
\nabla _{F.1.1}\omega ^{1} &=&(1+y^{2})\left( \omega ^{2}\otimes \omega
^{1}+\omega ^{1}\otimes \omega ^{1}+\omega ^{2}\otimes \omega ^{2}\right) ,
\\
\nabla _{F.1.1}\omega ^{2} &=&0,\quad R_{\nabla _{F.1.1}}\omega
^{1}=R_{\nabla _{F.1.1}}\omega ^{2}=0;
\end{eqnarray*}%
\begin{eqnarray*}
\nabla _{F.1.2}\omega ^{1} &=&(1+y)\left( \omega ^{1}\otimes \omega
^{1}+\omega ^{1}\otimes \omega ^{2}+\omega ^{2}\otimes \omega ^{2}\right) ,
\\
\nabla _{F.1.2}\omega ^{2} &=&0,\quad R_{\nabla _{F.1.2}}\omega
^{1}=R_{\nabla _{F.1.2}}\omega ^{2}=0;
\end{eqnarray*}%
\begin{eqnarray*}
\nabla _{F.1.3}\omega ^{1} &=&(1+y+y^{2})\omega ^{1}\otimes \omega
^{1}+y^{2}(\omega ^{1}\otimes \omega ^{2}+\omega ^{2}\otimes \omega
^{1})+\omega ^{2}\otimes \omega ^{2}, \\
\nabla _{F.1.3}\omega ^{2} &=&y^{2}(\omega ^{1}\otimes \omega ^{1}+\omega
^{1}\otimes \omega ^{2})+y\omega ^{2}\otimes \omega ^{1}+\omega ^{2}\otimes
\omega ^{2}, \\
R_{\nabla _{F.1.3}}\omega ^{1} &=& R_{\nabla _{F.1.3}}\omega
^{2}=0;
\end{eqnarray*}%
\begin{eqnarray*}
\nabla _{F.1.4}\omega ^{1} &=&(1+y^{2})\omega ^{1}\otimes \omega ^{1}+\left(
1+y\right) \omega ^{1}\otimes \omega ^{2}+\left( 1+y^{2}\right) \omega
^{2}\otimes \omega ^{1}+y^{2}\omega ^{2}\otimes \omega ^{2}, \\
\nabla _{F.1.4}\omega ^{2} &=&(1+y+y^{2})\omega ^{1}\otimes \omega
^{2}+\left( y + y^2\right) \omega ^{2}\otimes \omega ^{2}, \\
R_{\nabla _{F.1.4}}\omega ^{1} &=&R_{\nabla _{F.1.4}}\omega ^{2}=0;
\end{eqnarray*}%
\begin{eqnarray*}
\nabla _{F.1.5}\omega ^{1} &=&\omega ^{1}\otimes \omega ^{1}+y\omega
^{1}\otimes \omega ^{2}+\left( 1+y+y^{2}\right) \omega ^{2}\otimes \omega
^{1}+(y + y^2)\omega ^{2}\otimes \omega ^{2}, \\
\nabla _{F.1.5}\omega ^{2} &=&(1+y)\omega ^{1}\otimes \omega ^{1}+\left(
1+y+y^{2}\right) \omega ^{1}\otimes \omega ^{2}+y\omega ^{2}\otimes \omega
^{1}+\omega ^{2}\otimes \omega ^{2}, \\
R_{\nabla _{F.1.5}}\omega ^{1} &=&y^2\mathrm{Vol}\otimes \omega ^{1}+\mathrm{%
Vol}\otimes \omega ^{2},\quad R_{\nabla _{F.1.5}}\omega ^{2}=\mathrm{Vol}%
\otimes \omega ^{1}+(y + y^2)\mathrm{Vol}\otimes \omega ^{2};
\end{eqnarray*}%
\begin{eqnarray*}
\nabla _{F.1.6}\omega ^{1} &=&\left( 1+y\right) \omega ^{1}\otimes \omega
^{1}+y\omega ^{2}\otimes \omega ^{1}+(1+y+y^{2})\omega ^{2}\otimes \omega
^{2}, \\
\nabla _{F.1.6}\omega ^{2} &=&\left( y+y^{2}\right) \left( \omega
^{1}\otimes \omega ^{2}+\omega ^{2}\otimes \omega ^{1}\right) , \\
R_{\nabla _{F.1.6}}\omega ^{1} &=&y^2\mathrm{Vol}\otimes \omega ^{1},\quad
R_{\nabla _{F.1.6}}\omega ^{2}=y\mathrm{Vol}\otimes \omega ^{1}+(y + y^2)\mathrm{Vol}\otimes \omega ^{2};
\end{eqnarray*}%
\begin{eqnarray*}
\nabla _{F.1.7}\omega ^{1} &=&y\omega ^{1}\otimes \omega
^{1}+(1+y+y^{2})\omega ^{2}\otimes \omega ^{2}, \\
\nabla _{F.1.7}\omega ^{2} &=&(1+y^{2})\omega ^{1}\otimes \omega ^{1}+\left(
1+y\right) \omega ^{1}\otimes \omega ^{2}+\omega ^{2}\otimes \omega
^{1}+\left( 1+y\right) \omega ^{2}\otimes \omega ^{2}, \\
R_{\nabla _{F.1.7}}\omega ^{1} &=&(y + y^2)\mathrm{Vol}\otimes \omega ^{1},\quad
R_{\nabla _{F.1.7}}\omega ^{2}=(1+y^2)\mathrm{Vol}\otimes \omega ^{1}+y^2\mathrm{Vol}\otimes \omega ^{2};
\end{eqnarray*}%
\begin{eqnarray*}
\nabla _{F.1.8}\omega ^{1} &=&y\omega ^{1}\otimes \omega
^{2}+(1+y+y^{2})\omega ^{2}\otimes \omega ^{1}+y\omega ^{2}\otimes \omega
^{2}, \\
\nabla _{F.1.8}\omega ^{2} &=&(y+y^{2})\omega ^{1}\otimes \omega
^{1}+y^{2}\omega ^{1}\otimes \omega ^{2}+\left( 1+y\right) \omega
^{2}\otimes \omega ^{1}+y\omega ^{2}\otimes \omega ^{2}, \\
R_{\nabla _{F.1.8}}\omega ^{1} &=&\left( 1+y^2\right) \mathrm{Vol}\otimes
\omega ^{1}+y\mathrm{Vol}\otimes \omega ^{2},\quad R_{\nabla _{F.1.8}}\omega
^{2}=\mathrm{Vol}\otimes \omega ^{1}+(1+y + y^2)\mathrm{Vol}\otimes \omega^{2};
\end{eqnarray*}%
\begin{eqnarray*}
\nabla _{F.1.9}\omega ^{1} &=&y^{2}\omega ^{1}\otimes \omega ^{1}+\left(
1+y^{2}\right) \omega ^{1}\otimes \omega ^{2}+\left( 1+y\right) \omega
^{2}\otimes \omega ^{1}+\omega ^{2}\otimes \omega ^{2}, \\
\nabla _{F.1.9}\omega ^{2} &=&y\omega ^{1}\otimes \omega ^{1}+\left(
1+y+y^{2}\right) \omega ^{1}\otimes \omega ^{2}+\omega ^{2}\otimes \omega
^{1}+\left( 1+y^{2}\right) \omega ^{2}\otimes \omega ^{2}, \\
R_{\nabla _{F.1.9}}\omega ^{1} &=&\left( 1+y\right) \mathrm{Vol}\otimes
\omega ^{1}+\mathrm{Vol}\otimes \omega ^{2},\quad R_{\nabla _{F.1.9}}\omega
^{2}=\mathrm{Vol}\otimes \omega ^{2};
\end{eqnarray*}%
\begin{eqnarray*}
\nabla _{F.1.10}\omega ^{1} &=&(1+y+y^{2})\omega ^{1}\otimes \omega
^{2}+\left( 1+y\right) \omega ^{2}\otimes \omega ^{1}+\left(
1+y+y^{2}\right) \omega ^{2}\otimes \omega ^{2}, \\
\nabla _{F.1.10}\omega ^{2} &=&\omega ^{1}\otimes \omega ^{1}+\omega
^{1}\otimes \omega ^{2}+\left( 1+y^{2}\right) \omega ^{2}\otimes \omega
^{1}+y\omega ^{2}\otimes \omega ^{2}, \\
R_{\nabla _{F.1.10}}\omega ^{1} &=&y^2\mathrm{Vol}\otimes \omega ^{1}+\left(
1+y + y^2\right) \mathrm{Vol}\otimes \omega ^{2},\quad R_{\nabla
_{F.1.10}}\omega ^{2}=y^2\mathrm{Vol}\otimes \omega ^{1}+\left( 1+y^2\right) 
\mathrm{Vol}\otimes \omega ^{2};
\end{eqnarray*}%
\begin{eqnarray*}
\nabla _{F.1.11}\omega ^{1} &=&y\left( \omega ^{1}\otimes \omega ^{1}+\omega
^{2}\otimes \omega ^{1}+\omega ^{2}\otimes \omega ^{2}\right) , \\
\nabla _{F.1.11}\omega ^{2} &=&\left( 1+y^{2}\right) \omega ^{1}\otimes
\omega ^{1}+y\omega ^{1}\otimes \omega ^{2}+\left( 1+y+y^{2}\right) \omega
^{2}\otimes \omega ^{1}+\left( 1+y^{2}\right) \omega ^{2}\otimes \omega ^{2},
\\
R_{\nabla _{F.1.11}}\omega ^{1} &=&R_{\nabla _{F.1.11}}\omega ^{2}=0;
\end{eqnarray*}%
\begin{eqnarray*}
\nabla _{F.1.12}\omega ^{1} &=&y\omega ^{1}\otimes \omega ^{1}+\left(
1+y^{2}\right) \omega ^{1}\otimes \omega ^{2}+y^{2}\omega ^{2}\otimes \omega
^{1}+\left( 1+y\right) \omega ^{2}\otimes \omega ^{2}, \\
\nabla _{F.1.12}\omega ^{2} &=&y\omega ^{1}\otimes \omega ^{1}+\omega
^{1}\otimes \omega ^{2}+\left( 1+y+y^{2}\right) \omega ^{2}\otimes \omega
^{1}+\left( 1+y^{2}\right) \omega ^{2}\otimes \omega ^{2}, \\
R_{\nabla _{F.1.12}}\omega ^{1} &=&y\mathrm{Vol}\otimes \omega ^{2},\quad
R_{\nabla _{F.1.12}}\omega ^{2}=\left( 1+y + y^2\right) \mathrm{Vol}\otimes
\omega ^{1}+y\mathrm{Vol}\otimes \omega ^{2}.
\end{eqnarray*}

2) ($\beta=1+y^2$) $g_{F.2}=(1+y)\omega ^{1}\otimes \omega
^{1}+\left( 1+y^2\right) (\omega ^{1}\otimes \omega ^{2}+\omega ^{2}\otimes
\omega ^{1})+y\omega ^{2}\otimes \omega ^{2}$

with 12 QLCs, five of which are flat: 
\begin{eqnarray*}
\nabla _{F.2.1}\omega ^{1} &=&0,\quad \nabla _{F.2.1}\omega ^{2}=\left(
1+y^{2}\right) (\omega ^{1}\otimes \omega ^{1}+\omega ^{1}\otimes \omega
^{2}+\omega ^{2}\otimes \omega ^{2}), \\
R_{\nabla _{F.2.1}}\omega ^{1} &=&R_{\nabla _{F.2.1}}\omega ^{2}=0;
\end{eqnarray*}%
\begin{eqnarray*}
\nabla _{F.2.2}\omega ^{1} &=&y\omega ^{1}\otimes \omega ^{1}+\omega
^{1}\otimes \omega ^{2}+\left( 1+y+y^{2}\right) \omega ^{2}\otimes \omega
^{1}+\left( 1+y^{2}\right) \omega ^{2}\otimes \omega ^{2}, \\
\nabla _{F.2.2}\omega ^{2} &=&(y+y^{2})\omega ^{1}\otimes \omega
^{1}+y(\omega ^{1}\otimes \omega ^{2}+\omega ^{2}\otimes \omega
^{1})+y^{2}\omega ^{2}\otimes \omega ^{2}, \\
R_{\nabla _{F.2.2}}\omega ^{1} &=&\left( y + y^2\right) \mathrm{Vol}\otimes
\omega ^{1}+(1+y)\mathrm{Vol}\otimes \omega ^{2},\quad R_{\nabla
_{F.2.2}}\omega ^{2}=\mathrm{Vol}\otimes \omega ^{1}+\left( 1+y + y^2\right) 
\mathrm{Vol}\otimes \omega ^{2};
\end{eqnarray*}%
\begin{eqnarray*}
\nabla _{F.2.3}\omega ^{1} &=&y\omega ^{1}\otimes \omega ^{1}+\left(
1+y+y^{2}\right) \omega ^{1}\otimes \omega ^{2}+\omega ^{2}\otimes \omega
^{2}, \\
\nabla _{F.2.3}\omega ^{2} &=&\left( 1+y+y^{2}\right) \omega ^{1}\otimes
\omega ^{1}+y(\omega ^{1}\otimes \omega ^{2}+\omega ^{2}\otimes \omega
^{1})+\omega ^{2}\otimes \omega ^{2}, \\
R_{\nabla _{F.2.3}}\omega ^{1} &=&\left( 1+y\right) \mathrm{Vol}\otimes
\omega ^{1}+y^2\mathrm{Vol}\otimes \omega ^{2},\quad R_{\nabla _{F.2.3}}\omega
^{2}=\mathrm{Vol}\otimes \omega ^{1};
\end{eqnarray*}%
\begin{eqnarray*}
\nabla _{F.2.4}\omega ^{1} &=&0,\quad \nabla _{F.2.4}\omega ^{2}=\left(
y+y^{2}\right) (\omega ^{1}\otimes \omega ^{1}+\omega ^{2}\otimes \omega
^{1}+\omega ^{2}\otimes \omega ^{2}), \\
R_{\nabla _{F.2.4}}\omega ^{1} &=&R_{\nabla _{F.2.4}}\omega ^{2}=0;
\end{eqnarray*}%
\begin{eqnarray*}
\nabla _{F.2.5}\omega ^{1} &=&(y+y^{2})\omega ^{1}\otimes \omega
^{1}+y^{2}\omega ^{1}\otimes \omega ^{2}+y\omega ^{2}\otimes \omega
^{1}+\left( y+y^{2}\right) \omega ^{2}\otimes \omega ^{2}, \\
\nabla _{F.2.5}\omega ^{2} &=&y\omega ^{1}\otimes \omega ^{1}+y^{2}\omega
^{1}\otimes \omega ^{2}+\left( y+y^{2}\right) \omega ^{2}\otimes \omega
^{1}+y\omega ^{2}\otimes \omega ^{2}, \\
R_{\nabla _{F.2.5}}\omega ^{1} &=&R_{\nabla _{F.2.5}}\omega ^{2}=0;
\end{eqnarray*}%
\begin{eqnarray*}
\nabla _{F.2.6}\omega ^{1} &=&(1+y)\omega ^{1}\otimes \omega
^{2}+(1+y)\omega ^{2}\otimes \omega ^{1}, \\
\nabla _{F.2.6}\omega ^{2} &=&\left( 1+y\right) \omega ^{1}\otimes \omega
^{1}+y\omega ^{1}\otimes \omega ^{2}+\left( y+y^{2}\right) \omega
^{2}\otimes \omega ^{1}+\left( 1+y^{2}\right) \omega ^{2}\otimes \omega ^{2},
\\
R_{\nabla _{F.2.6}}\omega ^{1} &=&y^2\mathrm{Vol}\otimes \omega ^{1}+\mathrm{Vol}\otimes \omega ^{2},\quad R_{\nabla _{F.2.6}}\omega ^{2}=\left(
1+y^2\right) \mathrm{Vol}\otimes \omega ^{2};
\end{eqnarray*}%
\begin{eqnarray*}
\nabla _{F.2.7}\omega ^{1} &=&\omega ^{1}\otimes \omega ^{1}+y\omega
^{2}\otimes \omega ^{1}, \\
\nabla _{F.2.7}\omega ^{2} &=&\left( 1+y\right) \omega ^{1}\otimes \omega
^{1}+\omega ^{2}\otimes \omega ^{1}+\left( y+y^{2}\right) \omega ^{2}\otimes
\omega ^{2}, \\
R_{\nabla _{F.2.7}}\omega ^{1} &=&R_{\nabla _{F.2.7}}\omega ^{2}=0;
\end{eqnarray*}%
\begin{eqnarray*}
\nabla _{F.2.8}\omega ^{1} &=&\omega ^{1}\otimes \omega ^{1}+(y+y^{2})\omega
^{1}\otimes \omega ^{2}+(1+y^{2})\omega ^{2}\otimes \omega ^{1}+y^{2}\omega
^{2}\otimes \omega ^{2}, \\
\nabla _{F.2.8}\omega ^{2} &=&\omega ^{1}\otimes \omega ^{1}+\left(
1+y^{2}\right) \omega ^{1}\otimes \omega ^{2}+\left( 1+y+y^{2}\right) \omega
^{2}\otimes \omega ^{1}, \\
R_{\nabla _{F.2.8}}\omega ^{1} &=&\left( 1+y + y^2\right) \mathrm{Vol}\otimes
\omega ^{1}+\left( y + y^2\right) \mathrm{Vol}\otimes \omega ^{2},\quad
R_{\nabla _{F.2.8}}\omega ^{2}=(y + y^2)(\mathrm{Vol}\otimes \omega ^{1}+\mathrm{Vol}\otimes \omega ^{2});
\end{eqnarray*}%
\begin{eqnarray*}
\nabla _{F.2.9}\omega ^{1} &=&\omega ^{1}\otimes \omega ^{1}+(1+y^{2})\omega
^{1}\otimes \omega ^{2}+(y+y^{2})\omega ^{2}\otimes \omega ^{1}+y^{2}\omega
^{2}\otimes \omega ^{2}, \\
\nabla _{F.2.9}\omega ^{2} &=&y^{2}\omega ^{1}\otimes \omega ^{1}+\left(
1+y+y^{2}\right) \left( \omega ^{1}\otimes \omega ^{2}+\omega ^{2}\otimes
\omega ^{1}\right) +\left( 1+y\right) \omega ^{2}\otimes \omega ^{2}, \\
R_{\nabla _{F.2.9}}\omega ^{1} &=&\left( 1+y\right) \mathrm{Vol}\otimes
\omega ^{1},\quad R_{\nabla _{F.2.9}}\omega ^{2}=y(\mathrm{Vol}\otimes
\omega ^{1}+\mathrm{Vol}\otimes \omega ^{2});
\end{eqnarray*}%
\begin{eqnarray*}
\nabla _{F.2.10}\omega ^{1} &=&\left( 1+y+y^{2}\right) \left( \omega
^{1}\otimes \omega ^{1}+\omega ^{1}\otimes \omega ^{2}\right) +(1+y)\omega
^{2}\otimes \omega ^{1}+\left( y+y^{2}\right) \omega ^{2}\otimes \omega ^{2},
\\
\nabla _{F.2.10}\omega ^{2} &=&\omega ^{1}\otimes \omega ^{1}+\left(
1+y\right) \omega ^{1}\otimes \omega ^{2}+\left( 1+y+y^{2}\right) \omega
^{2}\otimes \omega ^{1}+y\omega ^{2}\otimes \omega ^{2}, \\
R_{\nabla _{F.2.10}}\omega ^{1} &=&\left( 1+y\right) \mathrm{Vol}\otimes
\omega ^{1}+\left( 1+y + y^2\right) \mathrm{Vol}\otimes \omega ^{2},\quad
R_{\nabla _{F.2.10}}\omega ^{2}=\left( y + y^2\right) \mathrm{Vol}\otimes \omega
^{1}+y\mathrm{Vol}\otimes \omega ^{2};
\end{eqnarray*}%
\begin{eqnarray*}
\nabla _{F.2.11}\omega ^{1} &=&\left( 1+y^{2}\right) \omega ^{1}\otimes
\omega ^{1}+\left( y+y^{2}\right) \omega ^{1}\otimes \omega ^{2}+y^{2}\omega
^{2}\otimes \omega ^{1}+\left( 1+y\right) \omega ^{2}\otimes \omega ^{2}, \\
\nabla _{F.2.11}\omega ^{2} &=&y\omega ^{1}\otimes \omega ^{1}+y^{2}\omega
^{1}\otimes \omega ^{2}+\left( 1+y^{2}\right) \omega ^{2}\otimes \omega
^{1}+\left( 1+y\right) \omega ^{2}\otimes \omega ^{2}, \\
R_{\nabla _{F.2.11}}\omega ^{1} &=&\left( 1+y^2\right) \mathrm{Vol}\otimes
\omega ^{1}+(1+y)\mathrm{Vol}\otimes \omega ^{2},\quad R_{\nabla
_{F.2.11}}\omega ^{2}=\left( 1+y\right) \mathrm{Vol}\otimes \omega ^{1}+y^2\mathrm{Vol}\otimes \omega ^{2};
\end{eqnarray*}%
\begin{eqnarray*}
\nabla _{F.2.12}\omega ^{1} &=&\left( 1+y^{2}\right) \omega ^{1}\otimes
\omega ^{1}+\left( 1+y\right) \omega ^{1}\otimes \omega ^{2}+\left(
1+y+y^{2}\right) \omega ^{2}\otimes \omega ^{2}, \\
\nabla _{F.2.12}\omega ^{2} &=&(y+y^{2})\left( \omega ^{1}\otimes \omega
^{1}+\omega ^{1}\otimes \omega ^{2}\right) +\left( 1+y^{2}\right) \left(
\omega ^{2}\otimes \omega ^{1}+\omega ^{2}\otimes \omega ^{2}\right) , \\
R_{\nabla _{F.2.12}}\omega ^{1} &=&R_{\nabla _{F.2.12}}\omega ^{2}=0.
\end{eqnarray*}

3) ($\beta=1+y$) $g_{F.3}=\omega ^{1}\otimes \omega
^{1}+\left( 1+y\right) (\omega ^{1}\otimes \omega ^{2}+\omega ^{2}\otimes
\omega ^{1})+\left( 1+y^2\right) \omega ^{2}\otimes \omega ^{2}$

with 12 QLCs, five of which are flat:
\begin{eqnarray*}
\nabla _{F.3.1}\omega ^{1} &=&y^{2}\omega ^{1}\otimes \omega ^{1}+\left(
1+y+y^{2}\right) (\omega ^{1}\otimes \omega ^{2}+\omega ^{2}\otimes \omega
^{1})+\left( 1+y\right) \omega ^{2}\otimes \omega ^{2}, \\
\quad \nabla _{F.3.1}\omega ^{2} &=&\left( 1+y^{2}\right) \omega ^{1}\otimes
\omega ^{1}+\omega ^{1}\otimes \omega ^{2}+y\omega ^{2}\otimes \omega
^{1}+\left( 1+y+y^{2}\right) \omega ^{2}\otimes \omega ^{2}, \\
R_{\nabla _{F.3.1}}\omega ^{1} &=&\left( 1+y^2\right) \mathrm{Vol}\otimes
\omega ^{1}+(1+y)\mathrm{Vol}\otimes \omega ^{2},\quad R_{\nabla
_{F.3.1}}\omega ^{2}=y\mathrm{Vol}\otimes \omega ^{1}+(y + y^2)\mathrm{Vol}\otimes \omega ^{2};
\end{eqnarray*}%
\begin{eqnarray*}
\nabla _{F.3.2}\omega ^{1} &=&\left( 1+y\right) \omega ^{1}\otimes \omega
^{1}+\left( y+y^{2}\right) (\omega ^{1}\otimes \omega ^{2}+\omega
^{2}\otimes \omega ^{1})+y\omega ^{2}\otimes \omega ^{2}, \\
\quad \nabla _{F.3.2}\omega ^{2} &=&y^{2}\omega ^{1}\otimes \omega
^{1}+y\omega ^{1}\otimes \omega ^{2}+\left( y+y^{2}\right) \omega
^{2}\otimes \omega ^{2}, \\
R_{\nabla _{F.3.2}}\omega ^{1} &=&\left( 1+y\right) (\mathrm{Vol}\otimes
\omega ^{1}+\mathrm{Vol}\otimes \omega ^{2}),\quad R_{\nabla _{F.3.2}}\omega
^{2}=\left( 1+y^2\right) \mathrm{Vol}\otimes \omega ^{1}+\mathrm{Vol}\otimes
\omega ^{2};
\end{eqnarray*}%
\begin{eqnarray*}
\nabla _{F.3.3}\omega ^{1} &=&y\omega ^{1}\otimes \omega ^{1}+\left(
1+y^{2}\right) (\omega ^{1}\otimes \omega ^{2}+\omega ^{2}\otimes \omega
^{1})+(1+y+y^{2})\omega ^{2}\otimes \omega ^{2}, \\
\quad \nabla _{F.3.3}\omega ^{2} &=&\omega ^{1}\otimes \omega
^{1}+(1+y+y^{2})\omega ^{1}\otimes \omega ^{2}+y\omega ^{2}\otimes \omega
^{1}+\left( y+y^{2}\right) \omega ^{2}\otimes \omega ^{2}, \\
R_{\nabla _{F.3.3}}\omega ^{1} &=&\mathrm{Vol}\otimes \omega ^{2},\quad
R_{\nabla _{F.3.3}}\omega ^{2}=y\mathrm{Vol}\otimes \omega ^{1}+\left(
1+y\right) \mathrm{Vol}\otimes \omega ^{2};
\end{eqnarray*}%
\begin{eqnarray*}
\nabla _{F.3.4}\omega ^{1} &=&\left( y+y^{2}\right) \left( \omega
^{1}\otimes \omega ^{1}+\omega ^{1}\otimes \omega ^{2}+\omega ^{2}\otimes
\omega ^{2}\right) =\nabla _{F.3.4}\omega ^{2}, \\
R_{\nabla _{F.3.4}}\omega ^{1} &=&R_{\nabla _{F.3.4}}\omega ^{2}=0;
\end{eqnarray*}%
\begin{eqnarray*}
\nabla _{F.3.5}\omega ^{1} &=&y\omega ^{1}\otimes \omega ^{1}+\left(
1+y\right) \omega ^{1}\otimes \omega ^{2}+y^{2}\omega ^{2}\otimes \omega
^{1}+\omega ^{2}\otimes \omega ^{2}, \\
\quad \nabla _{F.3.5}\omega ^{2} &=&y^{2}\left( \omega ^{1}\otimes \omega
^{1}+\omega ^{1}\otimes \omega ^{2}+\omega ^{2}\otimes \omega ^{1}\right)
+\left( 1+y\right) \omega ^{2}\otimes \omega ^{2}, \\
R_{\nabla _{F.3.5}}\omega ^{1} &=&y^2\mathrm{Vol}\otimes \omega ^{1}+(1+y^2)\mathrm{Vol}\otimes \omega ^{2},\quad R_{\nabla _{F.3.5}}\omega ^{2}=\left(
1+y + y^2\right) \mathrm{Vol}\otimes \omega ^{2};
\end{eqnarray*}%
\begin{eqnarray*}
\nabla _{F.3.6}\omega ^{1} &=&y^{2}\omega ^{1}\otimes \omega ^{1}+\left(
y+y^{2}\right) \omega ^{1}\otimes \omega ^{2}+\omega ^{2}\otimes \omega
^{1}+\left( 1+y^{2}\right) \omega ^{2}\otimes \omega ^{2}, \\
\quad \nabla _{F.3.6}\omega ^{2} &=&y^{2}\omega ^{1}\otimes \omega
^{1}+\left( 1+y\right) \omega ^{1}\otimes \omega ^{2}+y^{2}\omega
^{2}\otimes \omega ^{1}+\left( 1+y^{2}\right) \omega ^{2}\otimes \omega ^{2},
\\
R_{\nabla _{F.3.6}}\omega ^{1} &=&(1+y + y^2)\mathrm{Vol}\otimes \omega
^{1},\quad R_{\nabla _{F.3.6}}\omega ^{2}=y^2\mathrm{Vol}\otimes \omega
^{2};
\end{eqnarray*}%
\begin{eqnarray*}
\nabla _{F.3.7}\omega ^{1} &=&\left( 1+y+y^{2}\right) \omega ^{1}\otimes
\omega ^{1}+y^{2}\omega ^{1}\otimes \omega ^{2}+\left( 1+y^{2}\right) \left(
\omega ^{2}\otimes \omega ^{1}+\omega ^{2}\otimes \omega ^{2}\right) , \\
\quad \nabla _{F.3.7}\omega ^{2} &=&y\omega ^{1}\otimes \omega ^{1}+\omega
^{1}\otimes \omega ^{2}+y^{2}\left( \omega ^{2}\otimes \omega ^{1}+\omega
^{2}\otimes \omega ^{2}\right) , \\
R_{\nabla _{F.3.7}}\omega ^{1} &=&R_{\nabla _{F.3.7}}\omega ^{2}=0;
\end{eqnarray*}%
\begin{eqnarray*}
\nabla _{F.3.8}\omega ^{1} &=&y\omega ^{1}\otimes \omega ^{2}+\left(
1+y\right) \omega ^{2}\otimes \omega ^{1}+y^{2}\omega ^{2}\otimes \omega
^{2}, \\
\quad \nabla _{F.3.8}\omega ^{2} &=&(1+y+y^{2})\omega ^{1}\otimes \omega
^{1}+\left( 1+y\right) \omega ^{1}\otimes \omega ^{2}+\left( y+y^{2}\right)
\omega ^{2}\otimes \omega ^{1}+y\omega ^{2}\otimes \omega ^{2}, \\
R_{\nabla _{F.3.8}}\omega ^{1} &=&(y + y^2)\mathrm{Vol}\otimes \omega ^{1}+y\mathrm{Vol}\otimes \omega ^{2},\quad R_{\nabla _{F.3.8}}\omega ^{2}=\left(
1+y\right) \mathrm{Vol}\otimes \omega ^{1}+\left( 1+y^2\right) \mathrm{Vol}\otimes \omega ^{2};
\end{eqnarray*}%
\begin{eqnarray*}
\nabla _{F.3.9}\omega ^{1} &=&\left( 1+y^{2}\right) \omega ^{1}\otimes
\omega ^{2}+\left( 1+y\right) \omega ^{2}\otimes \omega ^{2}, \\
\quad \nabla _{F.3.9}\omega ^{2} &=&\omega ^{1}\otimes \omega ^{1}+\left(
1+y\right) \left( \omega ^{1}\otimes \omega ^{2}+\omega ^{2}\otimes \omega
^{1}\right) +\left( 1+y^{2}\right) \omega ^{2}\otimes \omega ^{2}, \\
R_{\nabla _{F.3.9}}\omega ^{1} &=&R_{\nabla _{F.3.9}}\omega ^{2}=0;
\end{eqnarray*}%
\begin{eqnarray*}
\nabla _{F.3.10}\omega ^{1} &=&\left( 1+y\right) \left( \omega ^{1}\otimes
\omega ^{1}+\omega ^{2}\otimes \omega ^{1}+\omega ^{2}\otimes \omega
^{2}\right) =\nabla _{F.3.10}\omega ^{2}, \\
R_{\nabla _{F.3.10}}\omega ^{1} &=&R_{\nabla _{F.3.10}}\omega ^{2}=0;
\end{eqnarray*}%
\begin{eqnarray*}
\nabla _{F.3.11}\omega ^{1} &=&\left( 1+y^{2}\right) \omega ^{1}\otimes
\omega ^{1}+\left( 1+y\right) \omega ^{1}\otimes \omega ^{2}+\left(
1+y+y^{2}\right) \omega ^{2}\otimes \omega ^{1}+\omega ^{2}\otimes \omega
^{2}, \\
\quad \nabla _{F.3.11}\omega ^{2} &=&\left( y+y^{2}\right) \omega
^{1}\otimes \omega ^{1}+y^{2}\omega ^{1}\otimes \omega ^{2}+\left(
1+y^{2}\right) \omega ^{2}\otimes \omega ^{1}, \\
R_{\nabla _{F.3.11}}\omega ^{1} &=&\left( 1+y^2\right) \mathrm{Vol}\otimes
\omega ^{2},\quad R_{\nabla _{F.3.11}}\omega ^{2}=\left( 1+y\right) (\mathrm{Vol}\otimes \omega ^{1}+\mathrm{Vol}\otimes \omega ^{2});
\end{eqnarray*}%
\begin{eqnarray*}
\nabla _{F.3.12}\omega ^{1} &=&y^{2}\omega ^{1}\otimes \omega ^{1}+\left(
1+y+y^{2}\right) \omega ^{1}\otimes \omega ^{2}+\left( 1+y\right) \omega
^{2}\otimes \omega ^{1}+y^{2}\omega ^{2}\otimes \omega ^{2}, \\
\quad \nabla _{F.3.12}\omega ^{2} &=&\left( 1+y+y^{2}\right) \left( \omega
^{1}\otimes \omega ^{1}+\omega ^{2}\otimes \omega ^{1}+\omega ^{2}\otimes
\omega ^{2}\right) , \\
R_{\nabla _{F.3.12}}\omega ^{1} &=&R_{\nabla _{F.3.12}}\omega ^{2}=0.
\end{eqnarray*}

4) ($\beta=y$) $g_{F.4}=\left( 1+y^2\right) \omega ^{1}\otimes
\omega ^{1}+y(\omega ^{1}\otimes \omega ^{2}+\omega ^{2}\otimes \omega
^{1})+\left( y+y^2\right) \omega ^{2}\otimes \omega ^{2}$

with 4 QLCs, three of which are flat: 
\begin{eqnarray*}
\nabla _{F.4.1}\omega ^{1} &=&y^{2}\omega ^{1}\otimes \omega
^{1}+(y+y^{2})\omega ^{1}\otimes \omega ^{2}+\omega ^{2}\otimes \omega ^{2},
\\
\quad \nabla _{F.4.1}\omega ^{2} &=&\omega ^{1}\otimes \omega ^{1}+\omega
^{1}\otimes \omega ^{2}+y^{2}\omega ^{2}\otimes \omega ^{1}+\left(
y+y^{2}\right) \omega ^{2}\otimes \omega ^{2}, \\
R_{\nabla _{F.4.1}}\omega ^{1} &=&R_{\nabla _{F.4.1}}\omega ^{2}=0;
\end{eqnarray*}%
\begin{eqnarray*}
\nabla _{F.4.2}\omega ^{1} &=&\left( y+y^{2}\right) \omega ^{1}\otimes
\omega ^{1}+(1+y+y^{2})\left( \omega ^{1}\otimes \omega ^{2}+\omega
^{2}\otimes \omega ^{1}\right) +y^{2}\omega ^{2}\otimes \omega ^{2}, \\
\quad \nabla _{F.4.2}\omega ^{2} &=&y\omega ^{1}\otimes \omega ^{1}+\left(
y+y^{2}\right) \left( \omega ^{1}\otimes \omega ^{2}+\omega ^{2}\otimes
\omega ^{1}\right) +\left( 1+y+y^{2}\right) \omega ^{2}\otimes \omega ^{2},
\\
R_{\nabla _{F.4.2}}\omega ^{1} &=&y^2\mathrm{Vol}\otimes \omega ^{1},\quad
R_{\nabla _{F.4.2}}\omega ^{2}=y^2\mathrm{Vol}\otimes \omega ^{2};
\end{eqnarray*}%
\begin{eqnarray*}
\nabla _{F.4.3}\omega ^{1} &=&\left( 1+y+y^{2}\right) \omega ^{1}\otimes
\omega ^{1}+y\omega ^{1}\otimes \omega ^{2}+\omega ^{2}\otimes \omega ^{1},
\\
\quad \nabla _{F.4.3}\omega ^{2} &=&\omega ^{1}\otimes \omega ^{2}+\left(
1+y+y^{2}\right) \omega ^{2}\otimes \omega ^{1}+y\omega ^{2}\otimes \omega
^{2}, \\
R_{\nabla _{F.4.3}}\omega ^{1} &=&R_{\nabla _{F.4.3}}\omega ^{2}=0;
\end{eqnarray*}%
\begin{eqnarray*}
\nabla _{F.4.4}\omega ^{1} &=&\left( 1+y^{2}\right) \omega ^{1}\otimes
\omega ^{1}+\left( 1+y\right) \omega ^{1}\otimes \omega ^{2}+\left(
y+y^{2}\right) \omega ^{2}\otimes \omega ^{1}+\left( 1+y^{2}\right) \omega
^{2}\otimes \omega ^{2}, \\
\quad \nabla _{F.4.4}\omega ^{2} &=&\left( 1+y\right) \omega ^{1}\otimes
\omega ^{1}+\left( y+y^{2}\right) \omega ^{1}\otimes \omega ^{2}+\left(
1+y^{2}\right) \omega ^{2}\otimes \omega ^{1}+\left( 1+y\right) \omega
^{2}\otimes \omega ^{2}, \\
R_{\nabla _{F.4.4}}\omega ^{1} &=&R_{\nabla _{F.4.4}}\omega ^{2}=0.
\end{eqnarray*}
For this family there are 7 possible equivalent solutions for $\epsilon $ (that are invertible
and satisfy \eqref{epsrel},\eqref{dext1},\eqref{dext2},\eqref{qs} to provide a central
volume form with metric quantum symmetric). These are all nonzero, hence invertible, functions times 
$$\epsilon=\left( 
\begin{array}{cc}
y & 1 \\ 
1 & 1+y
\end{array}
\right),\quad \eps^{-1}=\frac{1}{1+y+y^2}\left(\begin{array}{cc}
1+y & 1 \\ 
1 & y
\end{array}
\right)$$ 
which is the one is used above. This has $\Vol=(\omega ^{1})^{2}$, while
the others correspond to a volume form which is the inverse of the corresponding nonzero function times this.

\section{Laplacians}

\label{seclap}

The Laplacian is defined by $\Delta=(\ ,\ )\nabla\extd$ and the quantum dimension by $\underline{\rm dim}=(\ ,\ )(\cg)$ as in Section~\ref{secpre}. 

{\bf The $m=1$ case.} Using the analysis in Section~\ref{secpre} with commutative coordinate algebra, we write $\cg=g\omega\otimes\omega$ and $\nabla\omega=\Gamma\omega\otimes
\omega$, then 
\begin{equation*}
\Delta f= {\frac{f^{\prime \prime }+f^{\prime }\Gamma}{\bar g}},\quad\quad \underline{\rm dim}={g\over \bar g}
\end{equation*}
where $\omega f=\bar f\omega$ and $\extd f=f'\omega$ in our basis, which one can compute for the $m=1$ examples above. For $n=3$ our central
metrics were not invertible and the theory does not apply. For the $n=2$
algebras in Table~1 we have $f^{\prime \prime }=0$ for all $f$ and we always
have $\Delta 1=0$ so we only need to give the value on the basis element $x$
, namely $\Delta x=\Gamma/\bar g$. Hence 
\begin{align*}
{\rm A}:&\quad g=1, \Gamma=0: & \quad \Delta=0 \\
{\rm A:}& \quad g=1, \Gamma=x: & \quad \Delta x=x \\
{\rm A:}&\quad g=1+x=\Gamma: & \quad \Delta x=1 \\
{\rm B:}&\quad g=1,\Gamma=0: & \quad \Delta=0 \\
{\rm C:}&\quad g=1, \Gamma=0: & \quad \Delta=0 \\
{\rm C:}&\quad g=1,\Gamma=x: & \quad \Delta x=x \\
{\rm C:}&\quad g=1, \Gamma=1+x: & \quad \Delta x=1+x \\
{\rm C:}&\quad g=x=\Gamma: & \quad \Delta x=1+x \\
{\rm C:}&\quad g=1+x=\Gamma: & \quad \Delta x=x
\end{align*}
remembering that $\bar x=1+x=x^{-1}$ in algebra C. Among the nonzero $\Delta$, we have ${\rm Tr}(\Delta)=1$ except for the
third case of algebra $A$, and these are precisely the ones which have an eigenvector with eigenvalue 1 (the other eigenvector is 1 with eigenvalue 0). The exception where the trace is zero has only the eigenvector 1 with zero eigenvalue and is not fully diagonalisable. Here ${\rm Tr}$ denotes the usual trace of  a linear map $A\shortrightarrow A$. All the quantum dimensions are nonzero elements of the algebra (and all except the last two are 1). 

{\bf The $m=2$ case.} We use the tensor formula (\ref{laplacian}) and we will list the
resulting operators in $M_{n}$ on our basis $\{x^{\mu }\}$ of $A$. When
listing the Laplacians, the main invariant is the dimension of the null-space
which is between 1 and $n$ and whether $\Delta $ can be diagonalised to a
basis of eigenvectors with the two eigenvalues $0,1$. 

{\em Algebra {\rm D} or ${\mathbb{F}}_{2}{\mathbb{Z}}_{3}$ with its
three metrics $g_{D.1}-g_{D.3}$.} In each case we compute the inverse metric and  the
 quantum dimension $\underline{\dim }=(\ ,\ )(\cg)$, as
\begin{align*}
g_{D.1}^{ij} &=(\omega ^{i},\omega ^{j})=\left( 
\begin{array}{cc}
0 & z^{2} \\ 
z^{2} & 1
\end{array}
\right) ,&\quad \underline{\dim }_{D.1}=1+1+1=1 \\
g_{D.2}^{ij} &=(\omega ^{i},\omega ^{j})=\left( 
\begin{array}{cc}
1 & z \\ 
z & 0
\end{array}
\right) ,&\quad   \underline{\dim }_{D.2}=1+1+1=1\\
g_{D.3}^{ij} &=(\omega ^{i},\omega ^{j})=\left( 
\begin{array}{cc}
z^{2} & 0 \\ 
0 & z
\end{array}
\right) ,&\quad \underline{\dim }_{D.3}=1+1=0.
\end{align*}
We then combine with the QLC's from Section~\ref{secDmodel} to find the following Laplacians and their traces:

$g_{D.1}$:\quad  $\Delta 1=0, \Delta z=1, \Delta z^2=z,\quad {\rm Tr}(\Delta)=0$;\\
$g_{D.2}$:\quad  $\Delta 1=0, \Delta z=z^{2}, \Delta z^{2}=1,\quad {\rm Tr}(\Delta)=0$;\\
$g_{D.3}$:\quad $\Delta =0$

independently of the four  QLCs in each case. We check the first of these  computations for $\Delta$ (for $\nabla _{D.1.1}$) by hand as a check of the implementation:

$\Delta _{D.1.1}x=(\ ,\ )_{D.1}\nabla _{D.1.1}\extd x=
(\ ,\
)_{D.1}\nabla _{D.1.1}\omega ^{1}\\
\ \qquad=(\ ,\ )_{D.1}(z^{2}\omega ^{1}\otimes \omega ^{1}+(1+z)(\omega ^{1}\otimes
\omega ^{2}+\omega ^{2}\otimes \omega ^{1})+\omega ^{2}\otimes \omega ^{2})\\ 
\ \qquad=z^{2}\cdot 0+(1+z)z^{2}+(1+z)z^{2}+1\cdot 1=1$,

$\Delta _{D.1.1}y=(\ ,\ )_{D.1}\nabla _{D.1.1}\extd y=(\ ,\
)_{D.1}\nabla _{D.1.1}\omega ^{2}\\
=
(\ ,\ )_{D.1}(z^{2}\omega ^{1}\otimes \omega ^{1}+z\omega ^{1}\otimes
\omega ^{2}+z^{2}\omega ^{2}\otimes \omega ^{1}+\omega ^{2}\otimes \omega
^{2})\\
=z^{2}\cdot 0+z\cdot z^{2}+z^{2}\cdot z^{2}+1\cdot 1=z$.

{\em Algebra {\rm B} or ${\mathbb{F}}_{2}({\mathbb{Z}}_{3})$ with its only metric $g_B$.} The inverse metric, quantum dimension and Laplacian are
\begin{equation*}
g_{B}^{ij}=(\omega ^{i},\omega ^{j})=\left( 
\begin{array}{cc}
1+x & x+y \\ 
x+y & 1+y
\end{array}
\right) ,\quad\quad \underline{\dim }_B=0,\quad \Delta=0
\end{equation*}
for all four QLC's.

{\em Algebra {\rm F} or ${\mathbb{F}}_{8}=\F_2[y]/(y^3+y^2+1)$ with its
four metrics $g_{F.1}-g_{F.4}$ which
admit QLCs.} The corresponding inverse metrics and quantum dimensions are:

\begin{align*}
g_{F.1}^{ij} &=(\omega ^{i},\omega ^{j})=\left( 
\begin{array}{cc}
0 & 1+y+y^{2} \\ 
1+y+y^{2} & 1
\end{array}
\right) ,&\quad \underline{\dim }_{F.1}=1 \\
g_{F.2}^{ij} &=(\omega ^{i},\omega ^{j})=\left( 
\begin{array}{cc}
1 & 1+y \\ 
1+y & 1
\end{array}
\right) ,&\quad \underline{\dim }_{F.2}=1\\
g_{F.3}^{ij} &=(\omega ^{i},\omega ^{j})=\left( 
\begin{array}{cc}
1 & y^{2} \\ 
y^{2} & 0
\end{array}
\right) ,&  \quad  \underline{\dim }_{F.3}=1\\
g_{F.4}^{ij} &=(\omega ^{i},\omega ^{j})=\left( 
\begin{array}{cc}
y+y^{2} & 1+y \\ 
1+y & 1+y^{2}
\end{array}
\right),&\quad \underline{\dim }_{F.4}=0.
\end{align*}

For each metric the Laplacians can be grouped into 3 cases depending on the connection:

Metric $g_{F.1}$\newline
$\nabla _{F.1.1},\nabla
_{F.1.6},\nabla _{F.1.9},\nabla _{F.1.12}$:\quad $\Delta 1 =0, \Delta y=0, \Delta y^{2}=y^{2},\quad\quad \quad{\rm Tr}(\Delta)=1$

$\nabla
_{F.1.2},\nabla _{F.1.5},\nabla _{F.1.7},\nabla _{F.1.8}$:\quad  $\Delta 1=0, \Delta y=0, \Delta y^{2}=1+y+y^{2}, \quad{\rm Tr}(\Delta)=1$ 

$\nabla _{F.1.3},\nabla
_{F.1.4},\nabla _{F.1.10},\nabla _{F.1.11}$:\quad $\Delta 1=0, \Delta y=y^{2}, \Delta y^{2}=1,\quad\quad {\rm Tr}(\Delta)=0$.

Metric $g_{F.2}$\newline
 $\nabla _{F.2.1},\nabla
_{F.2.2},\nabla _{F.2.10},\nabla _{F.2.11}$:\quad  $\Delta 1=0, \Delta y=y, \Delta y^{2}=0,\quad\quad \quad{\rm Tr}(\Delta)=1$

$\nabla
_{F.2.4},\nabla _{F.2.6},\nabla _{F.2.8},\nabla _{F.2.9}$:\quad  $\Delta 1=0, \Delta y=1+y+y^{2}, \Delta y^{2}=0,\quad {\rm Tr}(\Delta)=1$

$\nabla_{F.2.3},\nabla _{F.2.5},\nabla _{F.2.7},\nabla _{F.2.12}$:\quad $\Delta 1 =0, \Delta y=y+y^2, \Delta y^{2}=1+y+y^{2},\quad {\rm Tr}(\Delta)=0.$

Metric $g_{F.3}$\newline
$\nabla _{F.3.1},\nabla
_{F.3.4},\nabla _{F.3.5},\nabla _{F.3.8}$:\quad  $\Delta 1=0,\Delta y=\Delta y^{2}=y^{2},\quad\quad {\rm Tr}(\Delta)=1$

$\nabla _{F.3.3},\nabla
_{F.3.6},\nabla _{F.3.10},\nabla _{F.3.11}$:\quad $\Delta 1=0, \Delta y=\Delta y^{2}=y,\quad\quad {\rm Tr}(\Delta)=1$ 

 $\nabla _{F.3.2},\nabla
_{F.3.7},\nabla _{F.3.9},\nabla _{F.3.12}$:\quad $\Delta 1=0, \Delta y=1, \Delta y^{2}=1+y,\quad {\rm Tr}(\Delta)=0$.

Metric $g_{F.4}$ (for all $\nabla_{F.4.1-F.4.4}$):\quad $\Delta =0$.

We now look at eigenvalues of $\Delta$. Of course, $1$ is always an eigenvector with eigenvalue 0. 

\begin{proposition}\label{proplap} For the $n=3,m=2$ examples above where metrics and QLCs exist (namely, algebras {\rm{B,D,F}}) we have 

(i) $\Delta=0$ if and only if $\underline{\rm dim}=0$; 

(ii)  If $\Delta\ne 0$ and ${\rm Tr}(\Delta)=1$ then $\Delta$ has one eigenvector with eigenvalue 1 and two with eigenvalue 0;

(iii) If $\Delta\ne 0$ and ${\rm Tr}(\Delta)=0$ then $\Delta$ has one eigenvector (namely 1) with eigenvalue 0 and no other eigenvectors.
\end{proposition}

Hence the reasonable case for physics seems to be when ${\rm Tr}(\Delta)=1$ which for $n=3,m=2$ also entails that $\underline{\rm dim}=1$. Specifically, the massive eigenvectors  $v$ are
\[  g_{F.1}:\quad y^2,\quad 1+y+y^2;\quad g_{F.2}:\quad  y,\quad 1+y+y^2;\quad g_{F.3}:\quad y^2,\quad y\]
for the 6 relevant Laplacians in the above list for the algebra F, in order. These solve the massive Klein-Gordon equation $\Delta v+ v=0$ and in each case there are also two massless eigenvectors in the kernel of $\Delta$. 

\section{Ricci and Einstein tensors}\label{secricci}

As explained in Section~\ref{secpre}, the Ricci tensor for a general algebra $A$ with calculus $\Omega$ in the approach \cite{BegMa2}  requires the additional data of a
`lift' bimodule map $i:\Omega^2\to \Omega^1\otimes_A\Omega^1$. In our examples where $\Omega^2$ is one-dimensional with a
chosen central basis element $\Vol$,  this means  $i(\Vol)=I_{ij}\omega^i\otimes\omega^j$ for some central element $
I\in \Omega^1\otimes_A\Omega^1$ such that $\wedge(I)=\Vol$. We can then contract as explained in Section~\ref{secpre} to obtain ${\rm Ricci}\in \Omega^1\tens_A\Omega^1$ as in (\ref{ricciV}) and  adjust $I$ so that  \textrm{Ricci} has the same quantum symmetry as $\cg$ if possible. We will also be interested in when $\nabla\cdot{\rm Ricci}=0$, where $\nabla\cdot$ means to apply $\nabla$ on the element of $\Omega^1\tens_A\Omega^1$ (same as for the metric) and then contract the first two factors with $(\ ,\ )$. 

In this section, we use this method to construct Ricci and its scalar $S=(\ ,\ )({\rm Ricci})$ for our models and also  explore possible Einstein tensors. For the latter, the usual definition $\mathrm{Eins} =\mathrm{Ricci}-\frac{1}{2}{S}\cg$ makes no sense over ${\mathbb{F}}_{2}$ but in our case where the geometry is 2D we could take $\mathrm{Eins} =\mathrm{Ricci}-\frac{1}{\underline{\rm dim}}{S}\cg$ as proposed in \cite{BegMa2} for a 2D quantum geometry. Classically the quantum dimension would be 2 as per the usual Einstein tensor but for a quantum model it may have  a different value.  For the purpose of our exploration in ${\mathbb{F}}_{2}$ we actually have only two choices, 0, 1, for the coefficient of $S\cg$ and we will focus on the latter, which is consistent with $\underline{\rm dim}=1$ found to be of interest in Section~\ref{seclap}. Thus we set
\begin{equation}\label{eins}
\mathrm{Eins} :=\mathrm{Ricci}+{S}\cg
=(\mathrm{Ricci}_{\mu ij}+{S}_{\nu }g_{\rho ij}V^{\nu \rho
}{}_{\mu } )x^{\mu }\omega ^{i}\otimes \omega ^{j}\end{equation}
as the provisional definition of ${\rm Eins}$ in this section, whilst leaving open the possibility that in some models we might want ${\rm Ricci}$ alone or at the other extreme $S\cg$ alone as the Einstein tensor. We will be interested in the values of {\rm Eins} and if this is not zero (as it would be classically for a 2D manifold) then  $\nabla \cdot \mathrm{Eins}=0$. Here
\[ \nabla\cdot{\rm Eins}=\nabla\cdot{\rm Ricci}+((\ ,\ )\tens\id)(\extd S\tens \cg)=\nabla\cdot{\rm Ricci}+\extd S\]
given the properties of a connection, the inverse metric and $\nabla\cg=0$ for a QLC. In tensor component terms, this translates to
\begin{align*}
\nabla \cdot \mathrm{Eins}&=(\mathrm{Ricci}_{\mu mn}d^{\mu }{}_{\nu k}+\mathrm{Ricci}_{\mu in}\Gamma ^{i}{}_{\rho km}V^{\mu \rho }{}_{\nu}\\
&\qquad + \mathrm{Ricci}_{\mu ij}\Gamma ^{j}{}_{\beta sn}a^{i\beta }{}_{\alpha
t}V^{\mu \alpha }{}_{\lambda }V^{\lambda \sigma }{}_{\nu }\sigma ^{ts}{}_{\sigma
km})V^{\nu\zeta}{}_{\gamma}g_{\zeta}{}^{km}x^\gamma 
\omega ^{n} +S_{\mu }d^{\mu }{}_{\nu i}x^\nu\omega^i.\end{align*}

{\bf The  $m=1$ case.}  We have a unique $i(\Vol
)=\omega\otimes\omega$ as there is no quantum symmetric metric freedom that could be added to this. We also have $R_\nabla(\omega)=\rho\Vol\otimes\omega$ where $
\rho=\Gamma^{\prime }+\Gamma\bar\Gamma$ so that
\begin{equation*}
\mathrm{Ricci}={\frac{g\rho}{\bar g}}\omega\otimes\omega,\quad S={\frac{g\rho
}{\bar g^2}},\quad {\mathrm{Eins}}=0
\end{equation*}
so the Ricci tensor and scalar contain the same information as the one
component $\rho$ of the Riemann curvature, as classically in 2 dimensions.
The only case of interest is $n=2$ (to have $g$ invertible) and the algebra $A$
with $g=1$ and $\Gamma=x$ as the only case with curvature, where $\rho=1$.
Then $\mathrm{Ricci}=\omega\otimes\omega=\cg$ and $S=1$ so that both parts of {\rm Eins} are separately conserved. We also have $
\underline\dim=1$.

{\bf The $n=3$, $m=2$ case.} We generally resort to computer calculations, starting with the model where we have also done computations by hand as a check. 

{\em For the {\rm D} model in Section~\ref{secDmodel}} we already saw in (\ref{liftZ3}) what the moduli of
central elements in $\Omega^1\otimes_A\Omega^1$ looks like, from which one finds by hand that the
most general form of $I=i(\Vol)$  is 
\begin{equation*}
i(\Vol)=z^2\omega^2\otimes\omega^1+ z\omega^2\otimes\omega^2+ \gamma\cg 
\end{equation*}
 for any function $\gamma=\gamma_1+\gamma_2 z+ \gamma_3 z^2$ is any element of the algebra (so there are 8 lifts as we vary $\gamma_i\in\F_2$), where we can take $\cg=g_{D.1}$ without loss of generality. 
 
 \begin{tiny}
\begin{table}[tbp]
\begin{tabular}{|l|l|l|l|l|l|}
\hline
Metric & QLC & Ricci (central for all $\gamma _{i}$) & $S=(\ ,\ )(
\mathrm{Ricci})$ & qua. symmetric & $\nabla \cdot\mathrm{Ricci=0}$ \\ \hline
$g_{D.1}$ & $\nabla _{D.1.2}$ & $\mathrm{Ricci}=0$  (flat connection)& $S=0$ & --- & --- \\
 & $\left.\begin{array}{c}\nabla_{D.1.1}\\ \nabla_{D.1.3}\\ \nabla_{D.1.4}\end{array}\right\}$ & 
$\begin{array}{l}
\mathrm{Ricci}=\left( \gamma _{3}+\gamma _{2}z^{2}\right) \omega
^{1}\otimes \omega ^{1} \\ 
+\left( \gamma _{2}z+\gamma _{3}z^{2}\right) \omega ^{1}\otimes \omega ^{2}
\\ 
+\left( \gamma _{1}+z+\gamma _{3}z^{2}\right) \omega ^{2}\otimes \omega
^{1} \\ 
+\left( 1+\gamma _{3}z+\gamma _{1}z^{2}\right) \omega ^{2}\otimes \omega
^{2}
\end{array}$
& $S=\gamma _{2}+\gamma _{3}z$ & 
$\begin{array}{l}
\gamma _{1}=0,\gamma _{2}=1: \\ 
\mathrm{Ricci}=\\
( 1+\gamma _{3} z) z^2 \omega ^{1}\otimes
\omega ^{1} \\ 
+( 1+\gamma _{3}z)z \omega ^{1}\otimes \omega ^{2} \\ 
+( 1+\gamma _{3}z)z \omega ^{2}\otimes \omega ^{1} \\ 
+( 1+\gamma _{3}z) \omega ^{2}\otimes \omega ^{2}
\end{array}$
& 
$\begin{array}{l}
\gamma _{1}=0=\gamma _{3}: \\ 
\mathrm{Ricci}\\
= \gamma _{2}z^{2}\omega ^{1}\otimes \omega ^{1}\\
\quad  +\gamma _{2}z\omega ^{1}\otimes \omega ^{2} \\ 
\quad +z\omega ^{2}\otimes \omega ^{1}+\omega ^{2}\otimes \omega ^{2}
\end{array}$
\\  \hline
$g_{D.2}$ & $\nabla _{D.2.4}$ & $\mathrm{Ricci}=0$  (flat connection)& $S=0$ & --- & --- \\ 
& $\left.\begin{array}{c}\nabla_{D.2.1}\\ \nabla_{D.2.2}\\ \nabla_{D.2.3}\end{array}\right\}$ & 
$\begin{array}{l}
\mathrm{Ricci}=\left( 1+\gamma _{3}z+\gamma _{1}z^{2}\right) \omega
^{1}\otimes \omega ^{1} \\ 
+\left( \gamma _{3}+\gamma _{1}z+z^{2}\right) \omega ^{1}\otimes \omega ^{2}
 \\ 
+\left( \gamma _{1}z+( 1+\gamma _{2}) z^{2}\right) \omega
^{2}\otimes \omega ^{1} \\ 
+\left( \gamma _{1}+( 1+\gamma _{2})z \right) \omega ^{2}\otimes
\omega ^{2}
\end{array}$
& $S=1+\gamma _{2}+\gamma _{1} z^{2}$ & 
$\begin{array}{l}
\gamma _{2}=0=\gamma _{3}: \\ 
\mathrm{Ricci}=\\
(\gamma _{1}+z)z^2 \omega ^{1}\otimes
\omega ^{1} \\ 
+( \gamma _{1}+z)z \omega ^{1}\otimes \omega ^{2} \\ 
+( \gamma _{1}+z)z \omega ^{2}\otimes \omega ^{1} \\ 
+( \gamma _{1}+z) \omega ^{2}\otimes \omega ^{2}
\end{array}$
& 
$\begin{array}{l} \gamma _{1}=0=\gamma _{3}: \\ 
\mathrm{Ricci}\\
=\omega ^{1}\otimes \omega ^{1} +z^{2}\omega ^{1}\otimes \omega ^{2}\\ 
\quad+\left( 1+\gamma _{2}\right) z^{2}\omega ^{2}\otimes \omega ^{1} \\ 
\quad+\left( 1+\gamma _{2}\right)z \omega ^{2}\otimes \omega ^{2}
\end{array}$\\
\hline
$g_{D.3}$ & $\nabla _{D.3.1}$ & $\mathrm{Ricci}=0$ (flat connection) & $S=0$  & 
--- & --- \\ 
& $\left.\begin{array}{c}\nabla_{D.3.2}\\ \nabla_{D.3.3}\\ \nabla_{D.3.4}\end{array}\right\}$ & 
$\begin{array}{l}
\mathrm{Ricci}=\left( \gamma _{1}+( 1+\gamma _{2})
z\right) \omega ^{1}\otimes \omega ^{1} \\ 
+\left( 1+\gamma _{2}+\gamma _{1}z^{2}\right) \omega ^{1}\otimes \omega ^{2}
 \\ 
+\left( \gamma _{2}+\gamma _{3}z\right) \omega ^{2}\otimes \omega ^{1} \\ 
+\left( \gamma _{3}+\gamma _{2}z^{2}\right) \omega ^{2}\otimes \omega ^{2}
\end{array}$
& $S=1+\gamma _{3} z+\gamma _{1} z^{2}$ & never qsymm & 
$\begin{array}{l}
\gamma _{1}=0=\gamma _{3}: \\ 
\mathrm{Ricci}\\
=\left( 1+\gamma _{2}\right) z\omega ^{1}\otimes
\omega ^{1} \\ 
\quad+\left( 1+\gamma _{2}\right) \omega ^{1}\otimes \omega ^{2} \\ 
\quad+\gamma _{2}\omega ^{2}\otimes \omega ^{1} \\ 
\quad+\gamma _{2}z^{2}\omega ^{2}\otimes \omega ^{2}
\end{array}$
\\  \hline
\end{tabular}
\caption{Ricci tensor and scalar for the algebra D. For each metric one connection is Ricci flat. For metrics with $\underline{\rm dim}=1$, the other connections have two lifts making Ricci quantum symmetric. \label{t4}}
\end{table}

\begin{table}[tbp]
\begin{tabular}{|l|l|l|l|l|l|}
\hline
Metric & QLC & $\mathrm{Eins}=\mathrm{Ricci}+{S}g$
& $\mathrm{Ricci}$ qsymm & $\nabla \cdot\mathrm{Eins}=0$ \\ \hline
$g_{D.1}$ & $\nabla _{D.1.2}$ & $\mathrm{Eins}=0$   (flat connection)&--- & --- \\
 & $\left.\begin{array}{c}\nabla_{D.1.1}\\ \nabla_{D.1.3}\\ \nabla_{D.1.4}\end{array}\right\}$ & 
$\begin{array}{l}
\mathrm{Eins}=\left( \gamma _{1}+z(1+\gamma _{2})\right) \omega
^{2}\otimes \omega ^{1} \\ 
\qquad+\left( 1+\gamma _{2}+\gamma _{1}z^{2}\right) \omega ^{2}\otimes \omega ^{2}
\end{array}$
 & 
${\rm Eins}=0$& 
$\begin{array}{l}
\gamma _{1}=0: \\ 
\mathrm{Eins}=\left(1+\gamma _{2}\right)z \omega ^{2}\otimes
\omega ^{1} \\ 
\qquad+\left( 1+\gamma _{2}\right) \omega ^{2}\otimes \omega ^{2}
\end{array}$\\  \hline
$g_{D.2}$ & $\nabla _{D.2.4}$ & $\mathrm{Eins}=0$   (flat connection)& --- & --- \\ 
& $\left.\begin{array}{c}\nabla_{D.2.1}\\ \nabla_{D.2.2}\\ \nabla_{D.2.3}\end{array}\right\}$ & 
$\begin{array}{l}
\mathrm{Eins}=\left( \gamma _{2}+\gamma _{3}z)\right) \omega
^{1}\otimes \omega ^{1} \\ 
\qquad+\left( \gamma _{3}+\gamma _{2}z^{2}\right) \omega ^{1}\otimes \omega ^{2}
\end{array}$
 & 
${\rm Eins}=0$
& 
$\begin{array}{l}
\gamma _{3}=0:\\ 
\mathrm{Eins}=\gamma _{2}\omega ^{1}\otimes \omega ^{1} \\ 
\qquad+\gamma _{2}z^{2}\omega ^{1}\otimes \omega ^{2}
\end{array}$\\
\hline
$g_{D.3}$ & $\nabla _{D.3.1}$ & $\mathrm{Eins}=0$  (flat connection)& --- & --- \\ 
& $\left.\begin{array}{c}\nabla_{D.3.2}\\ \nabla_{D.3.3}\\ \nabla_{D.3.4}\end{array}\right\}$ & 
$\begin{array}{l}
\mathrm{Eins}=\left( \gamma _{2}z+\gamma _{3}z^{2}\right) \omega
^{1}\otimes \omega ^{1} \\ 
+\left( \gamma _{2}+\gamma _{3}z\right) \omega ^{1}\otimes \omega ^{2} \\ 
+\left( 1+\gamma _{2}+\gamma _{1}z^{2}\right) \omega ^{2}\otimes \omega
^{1} \\ 
+\left( \gamma _{1}z+(1+\gamma _{2})z^{2}\right) \omega ^{2}\otimes \omega
^{2}
\end{array}$
 & never qsymm & 
$\begin{array}{l}
\gamma _{1}=0=\gamma _{3}: \\ 
\mathrm{Eins}=\gamma _{2}z\omega ^{1}\otimes \omega ^{1} \\ 
\qquad+\gamma _{2}\omega ^{1}\otimes \omega ^{2} \\ 
\qquad+\left( 1+\gamma _{2}\right) \omega ^{2}\otimes \omega ^{1} \\ 
\qquad+(1+\gamma _{2})z^{2}\omega ^{2}\otimes \omega ^{2}
\end{array}$
\\  \hline
\end{tabular}
\caption{Einstein tensor for the algebra D. Metrics where $\underline{\rm dim}=1$ have zero Einstein tensor when Ricci is lifted to be quantum symmetric. The metric $g_{D.3}$ where $\underline{\rm dim}=0$ has two lifts for the non-flat connections with $\nabla\cdot{\rm Eins}=0$ and $S=1$. \label{t5}}
\end{table}
\end{tiny}

The
resulting Ricci for the first metric $g_{D.1}$ and the first connection $\nabla_{D.1.1}$ is 
\begin{equation*}
\mathrm{Ricci}=z\omega^2\tens\omega^1+\omega^2\tens\omega^2+ g_{ij}(\omega^i\rho^j{}_k\gamma g_{rs},\omega^r)\omega^s\otimes\omega^k. 
\end{equation*}
The value of this and the results for all twelve connections are shown in Table~\ref{t4}. The Einstein tensors are listed in Table~\ref{t5}. For each metric, one connection is  flat and the other three connections all have the same Ricci curvature. This is quantum symmetrizable via two lifts whenever $\underline{\rm dim}=1$, in which case Einstein=0, and never symmetrizable when $\underline{\rm dim}=0$. In all cases we can more generally chose the lift so that $\nabla\cdot {\rm Eins}=0$. In particular, we have the following natural lifts for each metric and the group of 3 connections.
\[ g_{D.1}:\quad \gamma_1=\gamma_3=0, \gamma_2=1,\quad {\rm Ricci}=g_{D.1},\quad S=1,\quad \nabla\cdot{\rm Ricci}=0,\quad {\rm Eins}=0\]
\[ g_{D.2}:\quad \gamma_1=\gamma_2=\gamma_3=0, \quad {\rm Ricci}=g_{D.2},\quad S=1,\quad \nabla\cdot{\rm Ricci}=0,\quad {\rm Eins}=0\]
\[ g_{D.3}:\quad \gamma_1=\gamma_3=0, \quad \quad S=1,\quad \nabla\cdot{\rm Ricci}=\nabla\cdot{\rm Eins}=0,\quad {\rm Eins}\ne 0\]
where the last case is unusual in that classically the Einstein tensor in 2D would vanish, but this is also the `unphysical' case where $\underline{\rm dim}=0$ and $\Delta=0$. 

\begin{tiny}
\begin{table}[tbp]
\begin{tabular}{|l|l|l|l|l|l|}
\hline  QLC & Ricci (central for all $\gamma _{i}$) & $S=(\ ,\ )(
\mathrm{Ricci})$  & $\nabla \cdot\mathrm{Ricci=0}$  \\ \hline
  $\left.\begin{array}{c}\nabla_{B.1}\\ \nabla_{B.2}\\ \nabla_{B.3}\end{array}\right\}$ & $\mathrm{Ricci}=0$ (flat connections) & $S=0$ & --- \\ \hline
  $\nabla_{B.4}$ & 
$\begin{array}{l}
\mathrm{Ricci}=\left(\gamma_2x+(\gamma_{1}+\gamma_{2})(1+y)\right)\omega ^{1}\otimes
\omega ^{1}\\
+((\gamma_{1}+\gamma_{2})(1+x)+\gamma_{1}y)\omega ^{1}\otimes \omega ^{2}\\
+((1+\gamma_{1})x+(1+\gamma_{1}+\gamma_{3})(1+y))\omega ^{2}\otimes \omega ^{1}\\
+(\gamma_{3}y+(1+\gamma_{1}+\gamma_{3})(1+x))\omega ^{2}\otimes \omega ^{2}
\end{array}$
& $\begin{array}{l}S=\gamma_{2} +\left( 1+\gamma_{2}\right) x\\ \qquad+\left( \allowbreak
1+\gamma_{3}\right) (1+y)\end{array}$ & 
$\begin{array}{l}
\gamma _{2}=1=\gamma _{3}: \\ 
\mathrm{Ricci}\\
=(\gamma_2 x+(1+\gamma_{1})(1+y))\omega ^{1}\otimes \omega
^{1}\\
+((1+\gamma_{1})(1+x)+\gamma_{1}y)\omega ^{1}\otimes \omega ^{2}\\
+((1+\gamma_{1})x+\gamma_{1}(1+y))\omega ^{2}\otimes \omega ^{1}\\
+(y+\gamma_{1}(1+x))\omega ^{2}\otimes \omega ^{2}
\end{array}$\\ \hline
\end{tabular}
\caption{Ricci tensor and scalar for the algebra B for its unique metric $g_B$. Three connections are Ricci flat and $\nabla_{B.4}$  never has Ricci quantum symmetric, so that column is omitted. \label{t6}}
\end{table}
\end{tiny}

\begin{tiny}
\begin{table}[tbp]
\begin{tabular}{|l|l|l|l|l|l|}
\hline
QLC & $\mathrm{Eins}=\mathrm{Ricci}+{S}g$
&  $\nabla \cdot\mathrm{Eins}=0$ \\ \hline
$\left.\begin{array}{c}\nabla_{B.1}\\ \nabla_{B.2}\\ \nabla_{B.3}\end{array}\right\}$  & $\mathrm{Eins}=0$ (flat connections) & --- \\ \hline
  $\nabla _{B.4}$ & 
$\begin{array}{l}
\mathrm{Eins}=\left(x+(1+\gamma_{1}+\gamma_{3})(1+y)\right)\omega ^{1}\otimes
\omega ^{1}\\
+\left((1+\gamma_{1}+\gamma_{3})(1+x)+( 1+\gamma_{1}+\gamma_{3}+\gamma_{2})y\right)\omega^{1}\otimes \omega ^{2}\\
+\left((\gamma_{1}+\gamma_{2}+\gamma_{3})x+(\gamma_{1}+\gamma_{2})(1+y)\right)\omega ^{2}\otimes \omega
^{1}\\
+\left((\gamma_{1}+\gamma_{2})(1+x)+y\right)\omega ^{2}\otimes \omega ^{2}
\end{array}$
 & 
$\begin{array}{l}
\gamma _{2}=1=\gamma _{3}: \\ 
\mathrm{Eins}=( x+\gamma_{1}(1+y))\omega ^{1}\otimes \omega
^{1}\\
+(\gamma_{1}(1+x)+( 1+\gamma_{1})y )\omega ^{1}\otimes \omega ^{2}\\
+(\gamma_{1}x+(1+\gamma_{1})(1+y))\omega ^{2}\otimes \omega
^{1}\\
+((1+\gamma_{1})(1+x)+y)\omega ^{2}\otimes \omega ^{2}
\end{array}$ \\  \hline
\end{tabular}
\caption{Einstein tensor for the algebra B for its unique metric $g_B$ showing  two lifts with for $\nabla_{B.4}$ with $\nabla\cdot{\rm Eins}=0$ and $S=1$. \label{t7}}
\end{table}
\end{tiny}

{\em For the algebra {\rm B} model in Section~\ref{secBmodel}} the most general lifting map has the form
\begin{eqnarray*}
i(\mathrm{Vol}) &=&y\omega ^{1}\otimes \omega
^{2}+\left( 1+y\right) \omega ^{2}\otimes \omega ^{1}+\left( 1+x\right) \omega ^{2}\otimes
\omega ^{2}+\gamma\, g_B
\end{eqnarray*}
where $\gamma=\gamma_1+\gamma_2x+\gamma_3y$ is any element of the algebra (so there are 8 lifts). The Ricci and Einstein tensors as functions of $\gamma$ are shown in the first columns of Tables~\ref{t6} and~\ref{t7}. Connections $\nabla_{B.1},\nabla_{B.2},\nabla_{B.3}$ are  flat. For $\nabla_{B.4}$ we can never make Ricci quantum symmetric nor Eins=0, but we do have the two lifts (the free choices of $\gamma_1$) for which $\nabla\cdot {\rm Eins}=0$. 
\[ \nabla_{B.4}:\quad \gamma_2=\gamma_3=1, \quad\quad S=1,\quad \nabla\cdot{\rm Ricci}=\nabla\cdot{\rm Eins}=0,\quad {\rm Eins}\ne 0\]
similar to the group of 3 connections for $g_{D.3}$. Indeed, this is another case where  $\underline{\rm dim}=0$ and $\Delta=0$. 

{\em For the algebra {\rm F} model in Section~\ref{secFmodel}} the most general form of the lifting map is
\[ i(\mathrm{Vol}) =(1+y^2) \omega ^{1}\otimes \omega ^{1}+ y \omega ^{1}\otimes \omega
^{2}+ \omega ^{2}\otimes \omega ^{1}+ (1+y)\omega ^{2}\otimes \omega ^{2}+\gamma \cg\]
where $\gamma=\gamma_1+\gamma_2y+\gamma_3 y^2$ is any element of the algebra (so there are 8 lifts) and $\cg=g_{F.1}$, say, without loss of generality.  The Ricci and Einstein tensors are shown in Tables~\ref{t8} -- \ref{t15} below.

\begin{tiny}
\begin{table}[H] 
\begin{tabular}{|l|l|l|l|l|}
\hline
QLC & Ricci & $S=(\ ,\ )(\mathrm{Ricci})$ & qua. symmetric & $\nabla .
\mathrm{Ricci=0}$ \\ \hline
$\left. 
\begin{array}{c}
\nabla _{F.1.1} \\ 
\nabla _{F.1.2} \\ 
\nabla _{F.1.3} \\ 
\nabla _{F.1.4} \\ 
\nabla _{F.1.11}
\end{array}
\right\} $ & $\mathrm{Ricci}=0$ (flat connections)& $S=0$ & --- & --- \\ \hline
$\nabla _{F.1.5}$ & $
\begin{array}{l}
\mathrm{Ricci}=\\
\left( 1+\left( 1+\gamma_2\right) y+\left( \gamma
_{1}+\gamma _{3}\right) y^{2}\right) \omega ^{1}\otimes \omega ^{1} \\ 
+\left( \gamma_2+\gamma _{3}+\left( \gamma_1+\gamma_2\right)
y+\left( 1+\gamma_2+\gamma _{3}\right) y^{2}\right) \omega ^{1}\otimes
\omega ^{2} \\ 
+\left( \left( \gamma_1+\gamma_2\right) y+\left( \gamma_1+\gamma
_{2}+\gamma _{3}\right) y^{2}\right) \omega ^{2}\otimes \omega ^{1} \\ 
+\left( \left( \gamma_2+\gamma _{3}\right) y+\left(
1+\gamma_1\right)(1+ y^{2})\right) \omega ^{2}\otimes \omega ^{2} \\ 
\text{central for }\gamma_1=0,\gamma_2=1,\gamma _{3}=0
\end{array}
$ & $
\begin{array}{l}
S=\left( 1+\gamma_1\right)(1+ y) \\ 
+\left( \gamma_2+\gamma _{3}\right)(1+y^{2})
\end{array}
$ & $
\begin{array}{l}
\gamma_1=\gamma_2,\gamma _{3}=0: \\ 
\mathrm{Ricci}=\\
\left( 1+\left( 1+\gamma_1\right) y+\gamma
_{1}y^{2}\right) \omega ^{1}\otimes \omega ^{1} \\ 
+\left( \gamma_1+\left( 1+\gamma_1\right) y^{2}\right) \omega
^{1}\otimes \omega ^{2} \\ 
+\left( \gamma_1 y+\left( 1+\gamma_1\right) (1+y^{2})\right)
\omega ^{2}\otimes \omega ^{2}
\end{array}
$ & no sol. \\ \hline
$\nabla _{F.1.6}$ & $
\begin{array}{l}
\mathrm{Ricci}=\\
\left( 1+\gamma_2y+\left(
\gamma_1+\gamma_2\right) (1+y^{2})\right) \omega ^{1}\otimes \omega ^{1}
\\ 
+\left( \gamma_1+\gamma _{3}+y+\left( 1+\gamma_1+\gamma_2+\gamma
_{3}\right) y^{2}\right) \omega ^{1}\otimes \omega ^{2} \\ 
+\left( \gamma_2+\gamma _{3}+\left( \gamma_1+\gamma_2\right)
y+\left( 1+\gamma_2+\gamma _{3}\right) y^{2}\right) \omega ^{2}\otimes
\omega ^{1} \\ 
+\left( \gamma_1+\gamma_2+\gamma _{3}(1+y^{2})\right) \omega
^{2}\otimes \omega ^{2} \\ 
\text{{*}not central}
\end{array}
$ & $S$ as for $\nabla_{F.1.5}$ & $
\begin{array}{l}
\gamma_1=\gamma_2,\gamma _{3}=0: \\ 
\mathrm{Ricci}=\\
\left( 1+\gamma_1y\right) \omega ^{1}\otimes
\omega ^{1} \\ 
+\left( \gamma_1+y+y^{2}\right) \omega ^{1}\otimes \omega ^{2} \\ 
+\left( \gamma_1+\left( 1+\gamma_1\right) y^{2}\right) \omega
^{2}\otimes \omega ^{1}
\end{array}
$ & no sol. \\ \hline
$\nabla _{F.1.7}$ & $
\begin{array}{l}
\mathrm{Ricci}=\\
\left( \gamma_1+\left( 1+\gamma_1+\gamma
_{2}+\gamma _{3}\right) y+\left( 1+\gamma_1+\gamma_2\right)
y^{2}\right) \omega ^{1}\otimes \omega ^{1} \\ 
+\left( \gamma _{3}+y+\gamma_2y^{2}\right) \omega ^{1}\otimes \omega ^{2}
\\ 
+\left( \gamma_1+\gamma _{3}(1+y)+\left( 1+\gamma_2+\gamma
_{3}\right) y^{2}\right) \omega ^{2}\otimes \omega ^{1} \\ 
+\left( \gamma_2+y+\left( 1+\gamma_1+\gamma_2\right) y^{2}\right)
\omega ^{2}\otimes \omega ^{2} \\ 
\text{{*}not central}
\end{array}
$ & $S$ as for $\nabla_{F.1.5}$ & $
\begin{array}{l}
\gamma_1=\gamma_2,\gamma _{3}=0: \\ 
\mathrm{Ricci}=\\
\left( \gamma_1+y+y^{2}\right) \omega ^{1}\otimes
\omega ^{1} \\ 
+\left( y+\gamma_1y^{2}\right) \omega ^{1}\otimes \omega ^{2} \\ 
+\left( \gamma_1+\left( 1+\gamma_1\right) y^{2}\right) \omega
^{2}\otimes \omega ^{1} \\ 
+\left( \gamma_1+y+y^{2}\right) \omega ^{2}\otimes \omega ^{2}
\end{array}
$ & no sol. \\ \hline
$\nabla _{F.1.8}$ & $
\begin{array}{l}
\mathrm{Ricci}=\\
\left( 1+\gamma_2+\gamma_1(1+y)+\gamma
_{3}y^{2}\right) \omega ^{1}\otimes \omega ^{1} \\ 
+\left( 1+\gamma_2+\left( 1+\gamma_1+\gamma _{3}\right) y+y^{2}\right)
\omega ^{1}\otimes \omega ^{2} \\ 
+\left( 1+\gamma_2+\gamma _{3}+(\gamma_1+\gamma _{3})y+\left( 1+\gamma
_{2}\right) y^{2}\right) \omega ^{2}\otimes \omega ^{1} \\ 
+\left( \gamma_1+\gamma_2+\left( \gamma_1+\gamma_2+\gamma
_{3}\right) y\right) \omega ^{2}\otimes \omega ^{2} \\ 
\text{central  for }\gamma_1=1,\gamma_2=0,\gamma _{3}=1
\end{array}
$ & $S$ as for $\nabla_{F.1.5}$ & $
\begin{array}{l}
\gamma_1=\gamma_2,\gamma _{3}=0: \\ 
\mathrm{Ricci}=\\
\left( 1+\gamma_1y\right) \omega ^{1}\otimes
\omega ^{1} \\ 
+\left( \left( 1+\gamma_1\right)(1+ y)+y^{2}\right) \omega
^{1}\otimes \omega ^{2} \\ 
+\left( \gamma_1y+\left( 1+\gamma_1\right) (1+y^{2})\right)
\omega ^{2}\otimes \omega ^{1}
\end{array}
$ & no sol. \\ \hline
$\nabla _{F.1.9}$ & $
\begin{array}{l}
\mathrm{Ricci}=\\
\left( 1+\left( \gamma
_{2}+\gamma _{3}\right)(1+y)+(1+\gamma_1+\gamma _{3})y^{2}\right) \omega
^{1}\otimes \omega ^{1} \\ 
+\left( \gamma_1+\left( \gamma_1+\gamma_2\right) y+(1+\gamma
_{2})y^{2}\right) \omega ^{1}\otimes \omega ^{2} \\ 
+\left( (1+\gamma_1)(1+y)+\left( 1+\gamma_1+\gamma
_{2}+\gamma _{3}\right) y^{2}\right) \omega ^{2}\otimes \omega ^{1} \\ 
+\big( 1+\gamma_1+\gamma _{3}+\left( 1+\gamma_2+\gamma _{3}\right)
y\\
+\left( 1+\gamma_1+\gamma_2\right) y^{2}\big) \omega ^{2}\otimes
\omega ^{2} \\ 
\text{central for }\gamma_1=1,\gamma_2=0,\gamma _{3}=1
\end{array}
$ & $S$ as for $\nabla_{F.1.5}$ & $
\begin{array}{l}
\gamma_1=\gamma_2,\gamma _{3}=0: \\ 
\mathrm{Ricci}=\\
\left( \gamma_1y+(1+\gamma
_{1})(1+y^{2})\right) \omega ^{1}\otimes \omega ^{1} \\ 
+\left( \gamma_1+(1+\gamma_1)y^{2}\right) \omega ^{1}\otimes \omega
^{2} \\ 
+\left( 1+\gamma_1\right)(1+y+
y^{2}) \omega ^{2}\otimes \omega ^{1} \\ 
+\left( \left( 1+\gamma_1\right) (1+y)+y^{2}\right) \omega
^{2}\otimes \omega ^{2}
\end{array}
$ & no sol. \\ \hline
$\nabla _{F.1.10}$ & $
\begin{array}{l}
\mathrm{Ricci}=\\
\left( \gamma _{3}+\left( 1+\gamma
_{2}\right) (1+y)+(1+\gamma _{3})y^{2}\right) \omega ^{1}\otimes \omega ^{1} \\ 
+\left( \gamma_2+\gamma _{3}y+(1+\gamma_2)y^{2}\right) \omega
^{1}\otimes \omega ^{2} \\ 
+\left( \gamma_1+(\gamma_1+\gamma_2)y+\left( 1+\gamma_2\right)
y^{2}\right) \omega ^{2}\otimes \omega ^{1} \\ 
+\left( 1+\gamma_2y+\left( \gamma_1+\gamma
_{2}\right) (1+y^{2})\right) \omega ^{2}\otimes \omega ^{2} \\ 
\text{central for any }\gamma _{i}
\end{array}
$ & $
\begin{array}{l}
S=1+\gamma_1\\ 
+\left( \gamma_2+\gamma _{3}\right) (1+y)
\end{array}
$ & $
\begin{array}{l}
\gamma_1=\gamma_2,\gamma _{3}=0: \\ 
\mathrm{Ricci}=\\
\left( \left( 1+\gamma_1\right)
(1+y)+y^{2}\right) \omega ^{1}\otimes \omega ^{1} \\ 
+\left( \gamma_1+(1+\gamma_1)y^{2}\right) \omega ^{1}\otimes \omega
^{2} \\ 
+\left( \gamma_1+\left( 1+\gamma_1\right) y^{2}\right) \omega
^{2}\otimes \omega ^{1} \\ 
+\left( 1+\gamma_1y\right) \omega ^{2}\otimes \omega ^{2}
\end{array}
$ & $\gamma_1=0=\gamma_2$ \\ \hline
$\nabla _{F.1.12}$ & $
\begin{array}{l}
\mathrm{Ricci}=\\
\left( \gamma_1+\gamma_2+\gamma
_{3}(1+y^{2})\right) \omega ^{1}\otimes \omega ^{1} \\ 
+\left( \left( 1+\gamma_1+\gamma _{3}\right) (1+y)+\gamma
_{1}(1+y^{2})\right) \omega ^{1}\otimes \omega ^{2} \\ 
+\left( \gamma _{3}+y+\gamma_2y^{2}\right) \omega ^{2}\otimes \omega ^{1}
\\ 
+\left( \gamma_2+(\gamma_1+\gamma_2+\gamma _{3})y+\left( 1+\gamma
_{1}+\gamma _{3}\right) y^{2}\right) \omega ^{2}\otimes \omega ^{2} \\ 
\text{central for }\gamma_1=0,\gamma_2=1,\gamma _{3}=0
\end{array}
$ & $S$ as for $\nabla_{F.1.5}$ & $
\begin{array}{l}
\gamma_1=\gamma_2,\gamma _{3}=0: \\ 
\mathrm{Ricci}=\\
\left( 1+\left( 1+\gamma_1\right) y+\gamma
_{1}y^{2}\right) \omega ^{1}\otimes \omega ^{2} \\ 
+\left( 1+\gamma_1y\right)y \omega ^{2}\otimes \omega ^{1} \\ 
+\left( \gamma_1+\left( 1+\gamma_1\right) y^{2}\right) \omega
^{2}\otimes \omega ^{2}
\end{array}
$ & no sol. \\ \hline
\end{tabular}
\caption{Ricci tensor and scalar for the algebra F, metric $g_{F.1}$. \label{t8}} 
\end{table}\end{tiny}

\begin{tiny}
\begin{table}[H] 
\begin{tabular}{|l|l|l|l|}
\hline
QLC & $\mathrm{Eins}=\mathrm{Ricci}+{S}g$ & $\mathrm{Ricci}$ qsymm & $\nabla
.\mathrm{Eins}=0$ \\ \hline
$\left. 
\begin{array}{c}
\nabla _{F.1.1} \\ 
\nabla _{F.1.2} \\ 
\nabla _{F.1.3} \\ 
\nabla _{F.1.4} \\ 
\nabla _{F.1.11}
\end{array}
\right\} $ & $\mathrm{Eins}=0$ (flat connections)& --- & --- \\ \hline
$\nabla _{F.1.5}$ & $
\begin{array}{l}
\mathrm{Eins}=\left( \gamma _{3}+(1+\gamma
_{2})(1+y)+(1+\gamma _{3})y^{2}\right) \omega ^{1}\otimes \omega ^{1} \\ 
+\left( 1+\gamma_1+\left( \gamma_1+\gamma _{3}\right) y+(1+\gamma
_{2}+\gamma _{3})y^{2}\right) \omega ^{1}\otimes \omega ^{2} \\ 
+\left( 1+\left( \gamma_1+\gamma
_{3}\right) y+(\gamma_1+\gamma_2+\gamma _{3})(1+y^{2})\right) \omega
^{2}\otimes \omega ^{1} \\ 
+\left( \gamma_2+\gamma _{3}+\left( 1+\gamma_1+\gamma_2+\gamma
_{3}\right) (y+y^{2})\right) \omega
^{2}\otimes \omega ^{2}
\end{array}
$ & $
\begin{array}{l}
\mathrm{Eins}=\left( (1+\gamma_1)(1+y)+y^{2}\right)
\omega ^{1}\otimes \omega ^{1} \\ 
+\left( \gamma_1y+(1+\gamma_1)(1+y^{2})\right) \omega
^{1}\otimes \omega ^{2} \\ 
+\left( 1+\gamma_1y\right) \omega ^{2}\otimes \omega ^{1} \\ 
+\left( \gamma_1+y+y^{2}\right) \omega ^{2}\otimes \omega ^{2}
\end{array}
$ & no sol. \\ \hline
$\nabla _{F.1.6}$ & $
\begin{array}{l}
\mathrm{Eins}=\left( 1+\gamma_1+\gamma _{3}+\gamma
_{2}y+(1+\gamma_2)y^{2}\right) \omega ^{1}\otimes \omega ^{1} \\ 
+\left( 1+\gamma_2+\left( 1+\gamma_2+\gamma _{3}\right) y+(1+\gamma
_{1}+\gamma_2+\gamma _{3})y^{2}\right) \omega ^{1}\otimes \omega ^{2} \\ 
+\left( 1+\gamma_1+\left( \gamma_1+\gamma _{3}\right) y+(1+\gamma
_{2}+\gamma _{3})y^{2}\right) \omega ^{2}\otimes \omega ^{1} \\ 
+\left( 1+\left( 1+\gamma_1\right) y+\gamma_2y^{2}\right) \omega
^{2}\otimes \omega ^{2}
\end{array}
$ & $
\begin{array}{l}
\mathrm{Eins}=\left( \gamma_1y+(1+\gamma
_{1})(1+y^{2})\right) \omega ^{1}\otimes \omega ^{1} \\ 
+\left( \left( 1+\gamma_1\right) (1+y)+y^{2}\right) \omega
^{1}\otimes \omega ^{2} \\ 
+\left( \gamma_1y+(1+\gamma_1)(1+y^{2})\right) \omega
^{2}\otimes \omega ^{1} \\ 
+\left( 1+\left( 1+\gamma_1\right) y+\gamma_1y^{2}\right) \omega
^{2}\otimes \omega ^{2}
\end{array}
$ & no sol. \\  \hline
$\nabla _{F.1.7}$ & $
\begin{array}{l}
\mathrm{Eins}=\left(1+\left(
1+\gamma_1+\gamma_2+\gamma _{3}\right) (1+y)+\gamma_2y^{2}\right)
\omega ^{1}\otimes \omega ^{1} \\ 
+\left( 1+\gamma_1+\left( 1+\gamma_2+\gamma _{3}\right)
y+\gamma_2(1+y^{2})\right) \omega ^{1}\otimes \omega ^{2} \\ 
+\left( 1+\gamma_2(1+y)+(1+\gamma_2+\gamma _{3})y^{2}\right)
\omega ^{2}\otimes \omega ^{1} \\ 
+\left(\gamma_1y+(1+\gamma_1+\gamma
_{3})(1+y^{2})\right) \omega ^{2}\otimes \omega ^{2}
\end{array}
$ & $
\begin{array}{l}
\mathrm{Eins}=\left( 1+\gamma_1y\right)y \omega ^{1}\otimes
\omega ^{1} \\ 
+\left( 1+\left( 1+\gamma_1\right) y+\gamma_2y^{2}\right) \omega
^{1}\otimes \omega ^{2} \\ 
+\left( \gamma_1y+(1+\gamma_1)(1+y^{2})\right) \omega
^{2}\otimes \omega ^{1} \\ 
+\left(\gamma_1y+(1+\gamma_1)(1+y^{2})\right) \omega
^{2}\otimes \omega ^{2}
\end{array}
$ & no sol. \\ \hline
$\nabla _{F.1.8}$ & $
\begin{array}{l}
\mathrm{Eins}=\left( \gamma_1y+\left(
1+\gamma_1+\gamma _{3}\right) (1+y^{2})\right) \omega ^{1}\otimes \omega ^{1}
\\ 
+\left( \gamma_1+\gamma _{3}+\left( 1+\gamma_1+\gamma_2\right)
y+y^{2}\right) \omega ^{1}\otimes \omega ^{2} \\ 
+\left( \gamma_1+(\gamma_1+\gamma_2)y+(1+\gamma_2)y^{2}\right)
\omega ^{2}\otimes \omega ^{1} \\ 
+\left( 1+\gamma _{3}+(1+\gamma_2+\gamma _{3})y+(\gamma_2+\gamma
_{3})y^{2}\right) \omega ^{2}\otimes \omega ^{2}
\end{array}
$ & $
\begin{array}{l}
\mathrm{Eins}=\left( \gamma_1y+\left( 1+\gamma
_{1}\right) (1+y^{2})\right) \omega ^{1}\otimes \omega ^{1} \\ 
+\left( \gamma_1+y+y^{2}\right) \omega ^{1}\otimes \omega ^{2} \\ 
+\left( \gamma_1+(1+\gamma_1)y^{2}\right) \omega ^{2}\otimes \omega
^{1} \\ 
+\left( 1+(1+\gamma_1)y+\gamma_1y^{2}\right) \omega ^{2}\otimes \omega
^{2}
\end{array}
$ & no sol. \\ \hline
$\nabla _{F.1.9}$ & $
\begin{array}{l}
\mathrm{Eins}=\left( 1+(\gamma_2+\gamma _{3})y+\gamma
_{3}y^{2}\right) \omega ^{1}\otimes \omega ^{1} \\ 
+\left( \gamma _{3}+\left( \gamma_1+\gamma _{3}\right)
y+(1+\gamma_2)(1+y^{2})\right) \omega ^{1}\otimes \omega ^{2} \\ 
+\left( \gamma_2+\gamma _{3}+(1+\gamma_1+\gamma_2+\gamma _{3})(y+y^{2})\right) \omega
^{2}\otimes \omega ^{1} \\ 
+\left( \gamma_2+(\gamma_1+\gamma_2+\gamma _{3})y+(1+\gamma
_{1}+\gamma _{3})y^{2}\right) \omega ^{2}\otimes \omega ^{2}
\end{array}
$ & $
\begin{array}{l}
\mathrm{Eins}=\left( 1+\gamma_1y\right) \omega ^{1}\otimes \omega
^{1} \\ 
+\left( \gamma_1y+(1+\gamma_1)(1+y^{2})\right) \omega
^{1}\otimes \omega ^{2} \\ 
+\left( \gamma_1+y+y^{2}\right) \omega ^{2}\otimes \omega ^{1} \\ 
+\left( \gamma_1+(1+\gamma_1)y^{2}\right) \omega ^{2}\otimes \omega
^{2}
\end{array}
$ & no sol. \\ \hline
$\nabla _{F.1.10}$ & $
\begin{array}{l}
\mathrm{Eins}=\left(\gamma _{3}+(\gamma
_{1}+\gamma_2)\left(1+ y+y^{2}\right) \right) \omega ^{1}\otimes \omega ^{1}
\\ 
+\left( \gamma _{3}(1+y)+(\gamma_1+\gamma_2)y^{2}\right)
\omega ^{1}\otimes \omega ^{2} \\ 
+\left(\gamma _{3}+(\gamma_1+\gamma
_{2})(1+y+y^{2})\right) \omega ^{2}\otimes \omega ^{1} \\ 
+\left( \gamma _{3}(1+y)+(\gamma_1+\gamma_2)y^{2}\right)
\omega ^{2}\otimes \omega ^{2}
\end{array}
$ & $\mathrm{Eins}=0$ & $\gamma_1=\gamma_2,\gamma _{3}=0$ \\ 
\hline
$\nabla _{F.1.12}$ & $
\begin{array}{l}
\mathrm{Eins}=\left( \gamma_1+(1+\gamma_1+\gamma
_{3})y^{2}\right) \omega ^{1}\otimes \omega ^{1} \\ 
+\left( \gamma_2+\left( 1+\gamma_1+\gamma_2\right)
y+\gamma_1(1+y^{2})\right) \omega ^{1}\otimes \omega ^{2} \\ 
+\left( 1+\gamma_1+\gamma_2+(1+\gamma_2+\gamma _{3})y+\gamma
_{2}y^{2}\right) \omega ^{2}\otimes \omega ^{1} \\ 
+\left( 1+\gamma_1+\gamma _{3}+\left( 1+\gamma_2+\gamma _{3}\right)
y+(1+\gamma_1+\gamma_2)y^{2}\right) \omega ^{2}\otimes \omega ^{2}
\end{array}
$ & $
\begin{array}{l}
\mathrm{Eins}=\left( \gamma_1+(1+\gamma_1)y^{2}\right) \omega
^{1}\otimes \omega ^{1} \\ 
+\left( 1+\gamma_1y\right) y \omega ^{1}\otimes \omega ^{2} \\ 
+\left( 1+(1+\gamma_1)y+\gamma_2y^{2}\right) \omega ^{2}\otimes
\omega ^{1} \\ 
+\left( \left( 1+\gamma_1\right)(1+ y)+y^{2}\right) \omega
^{2}\otimes \omega ^{2}
\end{array}
$ & no sol. \\ \hline
\end{tabular}
\caption{Einstein tensor for the algebra F, metric $g_{F.1}$ showing a unique connection $\nabla_{F.1.10}$ which is not Ricci flat but has ${\rm Eins}=0$. \label{t9}} 
\end{table}
\end{tiny}

\begin{tiny}
\begin{table}[H] 
\begin{tabular}{|l|l|l|l|l|}
\hline
QLC & Ricci & $S=(\ ,\ )(\mathrm{Ricci})$ & qua. symmetric & $\nabla .
\mathrm{Ricci=0}$ \\ \hline
$\left. 
\begin{array}{c}
\nabla _{F.2.1} \\ 
\nabla _{F.2.4} \\ 
\nabla _{F.2.5} \\ 
\nabla _{F.2.7} \\ 
\nabla _{F.2.12}
\end{array}
\right\} $ & $\mathrm{Ricci}=0$ (flat connections) & $S=0$ & --- & --- \\ \hline
$\nabla _{F.2.2}$ & $
\begin{array}{l}
\mathrm{Ricci}= \\ 
\left( \gamma_2+\left( 1+\gamma_1+\gamma_2\right) y+\left( 1+\gamma
_{3}\right) y^{2}\right) \omega ^{1}\otimes \omega ^{1} \\ 
+\left( 1+\gamma_2+\gamma_2y+\left( 1+\gamma_2+\gamma_3\right)
y^{2}\right) \omega ^{1}\otimes \omega ^{2} \\ 
+\left( 1+\gamma_1+\left( 1+\gamma_2+\gamma_3\right)
y+\gamma_2(1+y^{2})\right) \omega ^{2}\otimes \omega ^{1} \\ 
+\left( \gamma_2+\left( 1+\gamma_1+\gamma_2\right)
y+\gamma_1(1+y^{2})\right) \omega ^{2}\otimes \omega ^{2} \\ 
\text{central for }\gamma_1=0,\gamma_2=1,\gamma_3=0
\end{array}
$ & $
\begin{array}{l}
S= \\ 
1+\gamma_3 +\gamma_1y^{2}
\\ 
+\left( 1+\gamma_1+\gamma_3\right) y 
\end{array}
$ & $
\begin{array}{l}
\gamma_1=1,\gamma_2=\gamma_3: \\ 
\mathrm{Ricci}= \\ 
\left(\gamma_2(1+y)+\left( 1+\gamma_2\right) y^{2}\right)
\omega ^{1}\otimes \omega ^{1} \\ 
+\left( 1+\gamma_2(1+y)+y^{2}\right) \omega ^{1}\otimes \omega
^{2} \\ 
+\left( y+\gamma_2(1+y^{2})\right) \omega ^{2}\otimes \omega ^{1}
\\ 
+\left( 1+\gamma_2(1+y)+y^{2}\right) \omega ^{2}\otimes \omega
^{2}
\end{array}
$ & no sol. \\ \hline
$\nabla _{F.2.3}$ & $
\begin{array}{l}
\mathrm{Ricci}= \\ 
\left( 1+\gamma_3+\gamma_1y\right) \omega ^{1}\otimes \omega ^{1} \\ 
+\left( 1+\gamma_3+\left( 1+\gamma_1+\gamma_3\right) y+\gamma
_{1}y^{2}\right) \omega ^{1}\otimes \omega ^{2} \\ 
+\left( 1+\gamma_2+\left( 1+\gamma_1+\gamma_3\right) y+y^{2}\right)
\omega ^{2}\otimes \omega ^{1} \\ 
+\left( \gamma_2+\left( \gamma_1+\gamma_2+\gamma_3\right)
y+\left( 1+\gamma_1+\gamma_3\right) y^{2}\right) \omega ^{2}\otimes
\omega ^{2} \\ 
\text{central for any }\gamma _{i}
\end{array}
$ & $
\begin{array}{l}
S=\gamma_1 \\ 
+\left( 1+\gamma_1+\gamma_3\right) y^{2}
\end{array}
$ & $
\begin{array}{l}
\gamma_1=1,\gamma_2=\gamma_3: \\ 
\mathrm{Ricci}= \\ 
\left( 1+\gamma_2+y\right) \omega ^{1}\otimes \omega ^{1} \\ 
+\left( 1+\gamma_2(1+y)+y^{2}\right) \omega ^{1}\otimes \omega
^{2} \\ 
+\left( 1+\gamma_2(1+y)+y^{2}\right) \omega ^{2}\otimes \omega
^{1} \\ 
+\left( \gamma_2+y+\left( 1+\gamma_1+\gamma_2\right) y^{2}\right)
\omega ^{2}\otimes \omega ^{2}
\end{array}
$ & $\begin{array}{l}\gamma_3=0,\\
\gamma_2=0 \end{array}$ \\ \hline
$\nabla _{F.2.6}$ & $
\begin{array}{l}
\mathrm{Ricci}= \\ 
\left( 1+\gamma_3+\left( 1+\gamma_1+\gamma_3\right) y+\gamma
_{1}y^{2}\right) \omega ^{1}\otimes \omega ^{1} \\ 
+\left( \left( \gamma_1+\gamma_2\right) y+(\gamma_1+\gamma
_{2}+\gamma_3)y^{2}\right) \omega ^{1}\otimes \omega ^{2} \\ 
+\left( \gamma_2+(\gamma_1+\gamma_2+\gamma_3)y+\left( 1+\gamma
_{1}+\gamma_3\right) y^{2}\right) \omega ^{2}\otimes \omega ^{1} \\ 
+\left( 1+\left( \gamma_2+\gamma_3\right) y+\gamma_3y^{2}\right)
\omega ^{2}\otimes \omega ^{2} \\ 
\text{{*}not central}
\end{array}
$ & $S$ as for $\nabla_{F.2.2}$ & $
\begin{array}{l}
\gamma_1=1,\gamma_2=\gamma_3: \\ 
\mathrm{Ricci}= \\ 
\left( 1+\gamma_2(1+y)+\gamma_1y^{2}\right) \omega ^{1}\otimes
\omega ^{1} \\ 
+\left( \left( 1+\gamma_2\right) y+y^{2}\right) \omega ^{1}\otimes \omega
^{2} \\ 
+\left( y+\gamma_2(1+y^{2})\right) \omega ^{2}\otimes \omega ^{1}
\\ 
+\left( 1+y+\gamma_2y^{2}\right) \omega ^{2}\otimes \omega ^{2}
\end{array}
$ & no sol. \\ \hline
$\nabla _{F.2.8}$ & $
\begin{array}{l}
\mathrm{Ricci}= \\ 
\left( 1+\gamma_1+(1+\gamma_2)y+y^{2}\right) \omega ^{1}\otimes \omega
^{1} \\ 
+\left( 1+\gamma_1+\gamma_2+y +\left( \gamma_1+\gamma_3\right) (y+y^{2})\right) \omega ^{1}\otimes
\omega ^{2} \\ 
+\left( \gamma_1+\gamma_3+(1+\gamma_1+\gamma_2)y+y^{2}\right)
\omega ^{2}\otimes \omega ^{1} \\ 
+\left( \gamma_2+\gamma_3+\left( 1+\gamma_1+\gamma_2+\gamma
_{3}\right) \left( y+y^{2}\right) \right) \omega ^{2}\otimes \omega ^{2} \\ 
\text{central for }\gamma_1=0,\gamma_2=1,\gamma_3=0
\end{array}
$ & $S$ as for $\nabla_{F.2.2}$ & $
\begin{array}{l}
\gamma_1=1,\gamma_2=\gamma_3: \\ 
\mathrm{Ricci}= \\ 
\left( (1+\gamma_2)y+y^{2}\right) \omega ^{1}\otimes \omega ^{1} \\ 
+\left( \gamma_2+y +\left( \gamma
_{1}+\gamma_2\right) (y+y^{2})\right) \omega ^{1}\otimes \omega ^{2} \\ 
+\left( 1+\gamma_2(1+y)+y^{2}\right) \omega ^{2}\otimes \omega
^{1}
\end{array}
$ & no sol. \\ \hline
$\nabla _{F.2.9}$ & $
\begin{array}{l}
\mathrm{Ricci}= \\ 
\left( 1+\left( \gamma_2+\gamma_3\right) y+\gamma_3y^{2}\right)
\omega ^{1}\otimes \omega ^{1} \\ 
+\left( 1+\gamma_2y+(\gamma_1+\gamma
_{2})(1+y^{2})\right) \omega ^{1}\otimes \omega ^{2} \\ 
+\left( 1+(1+\gamma_1)y+\gamma_2y^{2}\right) \omega ^{2}\otimes \omega
^{1} \\ 
+\left( 1+\left( \gamma_1+\gamma_2\right)
y+\left( 1+\gamma_2+\gamma_3\right) (1+y^{2})\right) \omega ^{2}\otimes
\omega ^{2} \\ 
\text{central for }\gamma_1=1,\gamma_2=0,\gamma_3=1
\end{array}
$ & $S$ as for $\nabla_{F.2.2}$ & $
\begin{array}{l}
\gamma_1=1,\gamma_2=\gamma_3: \\ 
\mathrm{Ricci}= \\ 
\left( 1+\left( \gamma_2+\gamma_3\right) y+\gamma_3y^{2}\right)
\omega ^{1}\otimes \omega ^{1} \\ 
+\left( 1+\gamma_2y+(\gamma_1+\gamma
_{2})(1+y^{2})\right) \omega ^{1}\otimes \omega ^{2} \\ 
+\left( 1+(1+\gamma_1)y+\gamma_2y^{2}\right) \omega ^{2}\otimes \omega
^{1} \\ 
+(1+\left( \gamma_1+\gamma_2\right)
y\\ +
\left( 1+\gamma_2+\gamma_3\right)(1+ y^{2})) \omega ^{2}\otimes
\omega ^{2}
\end{array}
$ & no sol. \\ \hline
$\nabla _{F.2.10}$ & $
\begin{array}{l}
\mathrm{Ricci}= \\ 
\left( \gamma_1+\gamma_3+\left( \gamma_1+\gamma_2+\gamma
_{3}\right) y+(1+\gamma_1)y^{2}\right) \omega ^{1}\otimes \omega ^{1} \\ 
+\left( \gamma_3+\left( 1+\gamma_1+\gamma_2\right) y+(\gamma
_{1}+\gamma_3)y^{2}\right) \omega ^{1}\otimes \omega ^{2} \\ 
+\left( \gamma_2+(1+\gamma_3)y+\left( \gamma
_{1}+\gamma_3\right) (1+y^{2})\right) \omega ^{2}\otimes \omega ^{1} \\ 
+( 1+\gamma_2+\left( 1+\gamma_2+\gamma_3\right)( y+y^2)+\gamma_1 y^{2}) \omega ^{2}\otimes
\omega ^{2} \\ 
\text{not central}
\end{array}
$ & $S$ as for $\nabla_{F.2.2}$ & $
\begin{array}{l}
\gamma_1=1,\gamma_2=\gamma_3: \\ 
\mathrm{Ricci}= \\ 
\left( 1+\gamma_2+y\right) \omega ^{1}\otimes \omega ^{1} \\ 
+\left( \gamma_2(1+y)+(1+\gamma_2)y^{2}\right) \omega
^{1}\otimes \omega ^{2} \\ 
+\left( 1+(1+\gamma_2)(y+y^2)\right) \omega
^{2}\otimes \omega ^{1} \\ 
+\left( 1+\gamma_2+y\right) \omega ^{2}\otimes \omega ^{2}
\end{array}
$ & no sol. \\ \hline
$\nabla _{F.2.11}$ & $
\begin{array}{l}
\mathrm{Ricci}= \\ 
\left( 1+\gamma_2+\left( 1+\gamma_1+\gamma_3\right) y+y^{2}\right)
\omega ^{1}\otimes \omega ^{1} \\ 
+\left( 1+\gamma_2+\gamma_3+\left( \gamma_1+\gamma_3\right)
y+\left( 1+\gamma_2\right) y^{2}\right) \omega ^{1}\otimes \omega ^{2} \\ 
+\left( \gamma_1+\gamma_2+\left( \gamma_1+\gamma_2+\gamma
_{3}\right) y\right) \omega ^{2}\otimes \omega ^{1} \\ 
+\left( \gamma_1+\gamma_3y
+\left( 1+\gamma
_{1}+\gamma_2\right)(y+ y^{2})\right) \omega ^{2}\otimes \omega ^{2} \\ 
\text{central for }\gamma_1=1,\gamma_2=0,\gamma_3=1
\end{array}
$ & $S$ as for $\nabla_{F.2.2}$ & $
\begin{array}{l}
\gamma_1=1,\gamma_2=\gamma_3: \\ 
\mathrm{Ricci}= \\ 
\left( 1+\gamma_2(1+y)+y^{2}\right) \omega ^{1}\otimes \omega ^{1}
\\ 
+\left( 1+\left( 1+\gamma_2\right) (y+y^2)\right) \omega ^{1}\otimes \omega ^{2} \\ 
+\left( 1+\gamma_2+y\right) \omega ^{2}\otimes \omega ^{1} \\ 
+\left( 1+\gamma_2y^{2}\right) \omega ^{2}\otimes \omega ^{2}
\end{array}
$ & no sol. \\ \hline
\end{tabular}
\caption{Ricci tensor and scalar for the algebra $F$, metric $g_{F.2}$.\label{t10}}
\end{table}
\end{tiny}

\begin{tiny}
\begin{table}[H] 
\begin{tabular}{|l|l|l|l|}
\hline
QLC & $\mathrm{Eins}=\mathrm{Ricci}+{S}g$ & $\mathrm{Ricci}$ qsymm & $\nabla
.\mathrm{Eins}=0$ \\ \hline
$\left. 
\begin{array}{c}
\nabla _{F.2.1} \\ 
\nabla _{F.2.4} \\ 
\nabla _{F.2.5} \\ 
\nabla _{F.2.7} \\ 
\nabla _{F.2.12}
\end{array}
\right\} $ & $\mathrm{Eins}=0$ (flat connections) & --- & --- \\ \hline
$\nabla _{F.2.2}$ & $
\begin{array}{l}
\mathrm{Eins}= \\ 
\left( \gamma_3+(1+\gamma_2)(1+y)+\gamma
_{1}(1+y^{2})\right) \omega ^{1}\otimes \omega ^{1} \\ 
+\left( 1+\gamma_2+\gamma_1y+(1+\gamma
_{1}+\gamma_2+\gamma_3)(y+y^{2})\right) \omega ^{1}\otimes \omega ^{2} \\ 
+\left( 1+\gamma_2y+(\gamma_1+\gamma
_{2})(1+y^{2})\right) \omega ^{2}\otimes \omega ^{1} \\ 
+\left( \gamma_2+\left( \gamma_1+\gamma_2+\gamma_3\right)
y+(1+\gamma_1+\gamma_3)y^{2}\right) \omega ^{2}\otimes \omega ^{2}
\end{array}
$ & $
\begin{array}{l}
\mathrm{Eins}= \\ 
\left( (1+\gamma_2)y+y^{2}\right) \omega
^{1}\otimes \omega ^{1} \\ 
+\left( 1+\gamma_2+y\right) \omega ^{1}\otimes \omega ^{2} \\ 
+\left( \gamma_2(1+y)+(1+\gamma_2)y^{2}\right) \omega
^{2}\otimes \omega ^{1} \\ 
+\left( \gamma_2+\gamma_1y+(1+\gamma_1+\gamma_2)y^{2}\right)
\omega ^{2}\otimes \omega ^{2}
\end{array}
$ & no sol. \\ \hline
$\nabla _{F.2.3}$ & $
\begin{array}{l}
\mathrm{Eins}= \\ 
\left( \gamma_2+\gamma_3+(1+\gamma_1)y^{2}\right) \omega
^{2}\otimes \omega ^{1} \\ 
+\left( 1+\gamma_1+\left( \gamma_2+\gamma
_{3}\right)(1+ y)\right) \omega ^{2}\otimes \omega ^{2}
\end{array}
$ & $\mathrm{Eins}_{F.2.3}=0$ & $\gamma_1=1,\gamma_2=\gamma_3$ \\ \hline
$\nabla _{F.2.6}$ & $
\begin{array}{l}
\mathrm{Eins}= \\ 
\left( \gamma_1+\left( 1+\gamma_3\right) \left( y+y^{2}\right) \right)
\omega ^{1}\otimes \omega ^{1} \\ 
+\left( \left( 1+\gamma_1\right)y+(\gamma
_{2}+\gamma_3)(y+y^{2})\right) \omega ^{1}\otimes \omega ^{2} \\ 
+\left( \gamma_2+(1+\gamma_1+\gamma_2)y+(1+\gamma_3)y^{2}\right)
\omega ^{2}\otimes \omega ^{1} \\ 
+\left( 1+\gamma_1+(1+\gamma_2)y+y^{2}\right) \omega ^{2}\otimes
\omega ^{2}
\end{array}
$ & $
\begin{array}{l}
\mathrm{Eins}\\
=\left( 1+\left( 1+\gamma_2\right) \left(
y+y^{2}\right) \right) \omega ^{1}\otimes \omega ^{1} \\ 
+\left( y^2+\gamma_2(1+y+y^{2})\right) \omega
^{2}\otimes \omega ^{1} \\ 
+\left( (1+\gamma_2)y+y^{2}\right) \omega ^{2}\otimes \omega ^{2}
\end{array}
$ & no sol. \\ \hline
$\nabla _{F.2.8}$ & $
\begin{array}{l}
\mathrm{Eins}= \\ 
\left( \gamma_3+\left( 1+\gamma_1+\gamma_2\right) y+\left( \gamma
_{1}+\gamma_3\right) y^{2}\right) \omega ^{1}\otimes \omega ^{1} \\ 
+\left( 1+\gamma_1+\gamma_2+\gamma_1y+\gamma_3y^{2}\right)
\omega ^{1}\otimes \omega ^{2} \\ 
+\left( \gamma_1+\gamma_3+(\gamma_1+\gamma_2+\gamma
_{3})y+(1+\gamma_1)y^{2}\right) \omega ^{2}\otimes \omega ^{1} \\ 
+\left( \gamma_3+(\gamma_1+\gamma_2)\left(
1+y+y^{2}\right) \right) \omega ^{2}\otimes \omega ^{2}
\end{array}
$ & $
\begin{array}{l}
\mathrm{Eins}= \\ 
\left(y^2+ \gamma_2(1+y+y^{2})\right)
\omega ^{1}\otimes \omega ^{1} \\ 
+\left( y+\gamma_2(1+y^{2})\right) \omega ^{1}\otimes \omega ^{2}
\\ 
+\left( 1+\gamma_2+y\right) \omega ^{2}\otimes \omega ^{1} \\ 
+\left( 1+(1+\gamma_2)\left( y+y^{2}\right) \right) \omega ^{2}\otimes
\omega ^{2}
\end{array}
$ & no sol. \\  \hline
$\nabla _{F.2.9}$ & $
\begin{array}{l}
\mathrm{Eins}= \\ 
\left( \gamma_1+\gamma_3+(\gamma_1+\gamma_2+\gamma
_{3})y+(1+\gamma_1)y^{2}\right) \omega ^{1}\otimes \omega ^{1} \\ 
+\left( 1+\gamma_1+\left( 1+\gamma_2+\gamma_3\right)
y+\gamma_2(1+y^{2})\right) \omega ^{1}\otimes \omega ^{2} \\ 
+\left( 1+(\gamma_1+\gamma_3)y+(\gamma_1+\gamma_2)y^{2}\right)
\omega ^{2}\otimes \omega ^{1} \\ 
+\left( 1+(1+\gamma_1+\gamma
_{2}+\gamma_3)(1+y)+\gamma_2y^{2}\right) \omega ^{2}\otimes \omega ^{2}
\end{array}
$ & $
\begin{array}{l}
\mathrm{Eins}= \\ 
\left( 1+\gamma_2+y\right) \omega ^{1}\otimes \omega ^{1} \\ 
+\left( y+\gamma_2(1+y^{2})\right) \omega ^{1}\otimes \omega ^{2}
\\ 
+\left( 1+(1+\gamma_2)\left( y+y^{2}\right) \right) \omega ^{2}\otimes
\omega ^{1} \\ 
+\left( 1+\gamma_2y^{2}\right) \omega ^{2}\otimes \omega ^{2}
\end{array}
$ & no sol. \\ \hline
$\nabla _{F.2.10}$ & $
\begin{array}{l}
\mathrm{Eins}= \\ 
\left( 1+(\gamma_2+\gamma_3)y+\gamma_3y^{2}\right) \omega
^{1}\otimes \omega ^{1} \\ 
+\left((\gamma_1+\gamma_2+\gamma_3)y+\gamma
_{3}(1+y^{2})\right) \omega ^{1}\otimes \omega ^{2} \\ 
+\left( \gamma_1+\gamma_2+\gamma_3(1+y^{2})\right) \omega
^{2}\otimes \omega ^{1} \\ 
+\left( 1+\gamma_2y+(\gamma_1+\gamma
_{2})(1+y^{2})\right) \omega ^{2}\otimes \omega ^{2}
\end{array}
$ & $
\begin{array}{l}
\mathrm{Eins}= \\ 
\left( 1+\gamma_2y^{2}\right) \omega ^{1}\otimes \omega ^{1} \\ 
+\left( y+\gamma_2(1+y^{2})\right) \omega ^{1}\otimes \omega ^{2}
\\ 
+\left( 1+\gamma_2y^{2}\right) \omega ^{2}\otimes \omega ^{1} \\ 
+\left( y^2+ \gamma_2(1+y+y^{2})\right) \omega
^{2}\otimes \omega ^{2}
\end{array}
$ & no sol. \\ \hline
$\nabla _{F.2.11}$ & $
\begin{array}{l}
\mathrm{Eins}= \\ 
\left( \gamma_2+\left( 1+\gamma_3\right)
y+(\gamma_1+\gamma_3)(1+y^{2})\right) \omega ^{1}\otimes \omega ^{1} \\ 
+\left( 1+\gamma_2+\gamma_3+\left( 1+\gamma_1\right) y+(1+\gamma
_{1}+\gamma_2)y^{2}\right) \omega ^{1}\otimes \omega ^{2} \\ 
+\left( 1+(1+\gamma_1+\gamma_2)(1+y)+\gamma
_{1}y^{2})\right) \omega ^{2}\otimes \omega ^{1} \\ 
+\left( \left( \gamma_1+\gamma_2\right) y+(\gamma_1+\gamma
_{2}+\gamma_3)y^{2}\right) \omega ^{2}\otimes \omega ^{2}
\end{array}
$ & $
\begin{array}{l}
\mathrm{Eins}= \\ 
\left( 1+\left( 1+\gamma_2\right) (y+y^{2})\right) \omega
^{1}\otimes \omega ^{1} \\ 
+\left( 1+\gamma_2y^{2}\right) \omega ^{1}\otimes \omega ^{2} \\ 
+\left( 1+\gamma_2(1+y)+y^{2})\right) \omega ^{2}\otimes \omega
^{1} \\ 
+\left( \left( 1+\gamma_2\right) y+y^{2}\right) \omega ^{2}\otimes \omega
^{2}
\end{array}
$ & no sol. \\ \hline
\end{tabular}
\caption{Einstein tensor for the algebra $F$, metric $g_{F.2}$ showing a unique connection $\nabla_{F.2.3}$ which is not Ricci flat but has ${\rm Eins}=0$.\label{t11}}
\end{table}
\end{tiny}

\begin{tiny}
\begin{table}[H] 
\begin{tabular}{|l|l|l|l|l|}
\hline
QLC & Ricci & $S=(\ ,\ )(\mathrm{Ricci})$ & qua. symmetric & $\nabla .
\mathrm{Ricci=0}$ \\ \hline
$\left. 
\begin{array}{c}
\nabla _{F.3.4} \\ 
\nabla _{F.3.7} \\ 
\nabla _{F.3.9} \\ 
\nabla _{F.3.10 } \\ 
\nabla _{F.3.12}
\end{array}
\right\} $ & $\mathrm{Ricci}=0$ (flat connections)& $S=0$ & --- & --- \\ \hline
$\nabla _{F.3.1}$ & $
\begin{array}{l}
\mathrm{Ricci}= \\ 
\left( \gamma_3(1+y)+\left( \gamma_1+\gamma_2\right)
y^{2}\right) \omega ^{1}\otimes \omega ^{1} \\ 
+\left( 1+\gamma_1+(\gamma_1+\gamma_3)y+\left( 1+\gamma_2+\gamma
_{3}\right) y^{2}\right) \omega ^{1}\otimes \omega ^{2} \\ 
+\left( 1+\gamma_2(1+y)+\left( 1+\gamma_2+\gamma_3\right)
y^{2}\right) \omega ^{2}\otimes \omega ^{1} \\ 
+\left( \gamma_1+\gamma_3+\left( \gamma_1+\gamma_2+\gamma
_{3}\right) y+(1+\gamma_1)y^{2}\right) \omega ^{2}\otimes \omega ^{2} \\ 
\text{central for }\gamma_1=0,\gamma_2=1,\gamma_3=0
\end{array}
$ & $
\begin{array}{l}
S= \\ 
\gamma_1+\gamma_2+\gamma_3y \\ 
+\left( \gamma_1+\gamma_2+\gamma_3\right) y^{2}
\end{array}
$ & $
\begin{array}{l}
\gamma_1=0,\gamma_3=1: \\ 
\mathrm{Ricci}= \\ 
\left( 1+y+\gamma_2y^{2}\right) \omega ^{1}\otimes \omega ^{1} \\ 
+\left( 1+y+\gamma_2y^{2}\right) \omega ^{1}\otimes \omega ^{2} \\ 
+\left( 1+\gamma_2(1+y+y^{2})\right) \omega
^{2}\otimes \omega ^{1} \\ 
+\left( 1+\gamma_1+\left( 1+\gamma_2\right) y+y^{2}\right) \omega
^{2}\otimes \omega ^{2}
\end{array}
$ & no sol. \\ \hline
$\nabla _{F.3.2}$ & $
\begin{array}{l}
\mathrm{Ricci}= \\ 
\left( 1+\left( 1+\gamma_1+\gamma
_{2}+\gamma_3\right) (1+y)+\gamma_2y^{2}\right) \omega ^{1}\otimes \omega
^{1} \\ 
+\left( \gamma_1+\gamma_3+y+\left( 1+\gamma_1+\gamma_2+\gamma
_{3}\right) y^{2}\right) \omega ^{1}\otimes \omega ^{2} \\ 
+\left( \gamma_3(1+y)+\left( \gamma_1+\gamma_2\right)
y^{2}\right) \omega ^{2}\otimes \omega ^{1} \\ 
+\left( \gamma_1+\gamma_2+\gamma_3(1+y^{2})\right) \omega
^{2}\otimes \omega ^{2} \\ 
\text{central for any }\gamma _{i}
\end{array}
$ & $
\begin{array}{l}
S=\gamma_3 +(\gamma_1+\gamma_2)\\ \quad\quad\quad (1+y+y^{2})
\end{array}
$ & $
\begin{array}{l}
\gamma_1=0,\gamma_3=1: \\ 
\mathrm{Ricci}= \\ 
\left( 1+\gamma_2+\gamma_2y+\gamma_2y^{2}\right) \omega ^{1}\otimes
\omega ^{1} \\ 
+\left( 1+y+\gamma_2y^{2}\right) \omega ^{1}\otimes \omega ^{2} \\ 
+\left( 1+y+\gamma_2y^{2}\right) \omega ^{2}\otimes \omega ^{1} \\ 
+\left( 1+\gamma_2+y^{2}\right) \omega ^{2}\otimes \omega ^{2}
\end{array}
$ & $\begin{array}{l}\gamma_2=0,\\
\gamma_1+\gamma_3=1\end{array}$ \\ \hline
$\nabla _{F.3.3}$ & $
\begin{array}{l}
\mathrm{Ricci}= \big( 1+\gamma_1+\gamma_3 \\
+ \left( 1+\gamma_2+\gamma_3\right)
y+(1+\gamma_1+\gamma_2)y^{2}\big) \omega ^{1}\otimes \omega ^{1} \\ 
+\left( \left( 1+\gamma_1+\gamma_3\right)
(1+y)+\left( 1+\gamma_3\right) y^{2}\right) \omega ^{1}\otimes \omega ^{2} \\ 
+\left( \left( 1+\gamma_1\right)
(1+y)+\left( \gamma_2+\gamma_3\right) (1+y^{2})\right) \omega ^{2}\otimes
\omega ^{1} \\ 
+\left(1+\left( 1+\gamma_1+\gamma
_{2}+\gamma_3\right)(1+ y)+\gamma_2y^{2}\right) \omega ^{2}\otimes \omega
^{2} \\ 
\text{central for }\gamma_1=1,\gamma_2=0,\gamma_3=1
\end{array}
$ & $S$ as for $\nabla_{F.3.1}$ & $
\begin{array}{l}
\gamma_1=0,\gamma_3=1: \\ 
\mathrm{Ricci}= \\ 
\left( \gamma_2y+(1+\gamma_2)y^{2}\right) \omega ^{1}\otimes \omega
^{1} \\ 
+\left( \gamma_2+y+\left( 1+\gamma_2\right) y^{2}\right) \omega
^{2}\otimes \omega ^{1} \\ 
+\left( 1+\gamma_2\left( 1+y+y^{2}\right) \right) \omega ^{2}\otimes
\omega ^{2}
\end{array}
$ & no sol. \\ \hline
$\nabla _{F.3.5}$ & $
\begin{array}{l}
\mathrm{Ricci}= \\ 
\left( 1+\gamma_2+\left( 1+\gamma_1+\gamma_3\right) y+y^{2}\right)
\omega ^{1}\otimes \omega ^{1} \\ 
+(\gamma_1+\gamma_3+\left( \gamma_1+\gamma_2+\gamma_3\right)
y+\left( 1+\gamma_1\right) y^{2})\omega ^{1}\otimes \omega ^{2} \\ 
+\left( \gamma_1+\gamma_2+\left( \gamma_1+\gamma_2+\gamma
_{3}\right) y\right) \omega ^{2}\otimes \omega ^{1} \\ 
+\left( \gamma_2+\left( 1+\gamma_3\right)
y+\left( \gamma_1+\gamma_3\right) (1+y^{2})\right) \omega ^{2}\otimes
\omega ^{2} \\ 
\text{not central}
\end{array}
$ & $S$ as for $\nabla_{F.3.1}$ & $
\begin{array}{l}
\gamma_1=0,\gamma_3=1: \\ 
\mathrm{Ricci}= \\ 
\left( 1+\gamma_2+y^{2}\right) \omega ^{1}\otimes \omega ^{1} \\ 
+(1+\left( 1+\gamma_2\right) y+y^{2})\omega ^{1}\otimes \omega ^{2} \\ 
+\left( \gamma_2+\left( 1+\gamma_2\right) y\right) \omega ^{2}\otimes
\omega ^{1} \\ 
+\left( 1+\gamma_2+y^{2}\right) \omega ^{2}\otimes \omega ^{2} \\ 
\text{not central}
\end{array}
$ & no sol. \\ \hline
$\nabla _{F.3.6}$ & $
\begin{array}{l}
\mathrm{Ricci}= \\ 
\left( 1+\gamma_3+\gamma_1y\right) \omega ^{1}\otimes \omega ^{1} \\ 
+(\gamma_2+\left( \gamma_1+\gamma_2+\gamma_3\right) y+\left(
1+\gamma_1+\gamma_3\right) y^{2})\omega ^{1}\otimes \omega ^{2} \\ 
+\left( 1+\gamma_2+\left( 1+\gamma_1+\gamma_3\right) y+y^{2}\right)
\omega ^{2}\otimes \omega ^{1} \\ 
+\left( \gamma_3(1+y)+\left( \gamma_1+\gamma
_{2}+\gamma_3\right) (1+y^{2})\right) \omega ^{2}\otimes \omega ^{2} \\ 
\text{not central}
\end{array}
$ & $S$ as for $\nabla_{F.3.1}$ & $
\begin{array}{l}
\gamma_1=0,\gamma_3=1: \\ 
\mathrm{Ricci}= \\ 
+(\gamma_2+\left( 1+\gamma_2\right) y)\omega ^{1}\otimes \omega ^{2}
\\ 
+\left( 1+\gamma_2+y^{2}\right) \omega ^{2}\otimes \omega ^{1} \\ 
+\left( \gamma_2+y+\left( 1+\gamma_2\right) y^{2}\right) \omega
^{2}\otimes \omega ^{2}
\end{array}
$ & no sol. \\ \hline
$\nabla _{F.3.8}$ & $
\begin{array}{l}
\mathrm{Ricci}= \\ 
\left( \gamma_1+\left( 1+\gamma_3\right) \left( y+y^{2}\right) \right)
\omega ^{1}\otimes \omega ^{1} \\ 
+(\gamma_3+y+\gamma_2y^{2})\omega ^{1}\otimes \omega ^{2} \\ 
+\left( 1+\gamma_1+\gamma_3+y+\left( 1+\gamma_2\right)(y+
y^{2})\right) \omega ^{2}\otimes \omega ^{1} \\ 
+\left( \gamma_2+y+\left( 1+\gamma_1+\gamma_2\right) y^{2}\right)
\omega ^{2}\otimes \omega ^{2} \\ 
\text{central for }\gamma_1=1,\gamma_2=0,\gamma_3=1
\end{array}
$ & $S$ as for $\nabla_{F.1.5}$ & $
\begin{array}{l}
\gamma_1=0,\gamma_3=1: \\ 
\mathrm{Ricci}= \\ 
(1+y+\gamma_2y^{2})\omega ^{1}\otimes \omega ^{2} \\ 
+(\gamma_2y+\left( 1+\gamma_2\right) y^{2})\omega ^{2}\otimes \omega
^{1} \\ 
+\left( \gamma_2+y+\left( 1+\gamma_2\right) y^{2}\right) \omega
^{2}\otimes \omega ^{2}
\end{array}
$ & no sol. \\ \hline
$\nabla _{F.3.11}$ & $
\begin{array}{l}
\mathrm{Ricci}= \\ 
\left( 1+\left( 1+\gamma_2\right) (1+y)+\left(
1+\gamma_3\right) (1+y^{2})\right) \omega ^{1}\otimes \omega ^{1} \\ 
+(\gamma_1+\gamma_2+\gamma_3(1+y^{2}))\omega ^{1}\otimes
\omega ^{2} \\ 
+\left( \gamma_1+(\gamma_1+\gamma_2)y+\left( 1+\gamma_2\right)
y^{2}\right) \omega ^{2}\otimes \omega ^{1} \\ 
+\left( \gamma_3+y+\gamma_2y^{2}\right) \omega ^{2}\otimes \omega ^{2}
\\ 
\text{central for }\gamma_1=0,\gamma_2=1,\gamma_3=0
\end{array}
$ & $S$ as for $\nabla_{F.3.1}$ & $
\begin{array}{l}
\gamma_1=0,\gamma_3=1: \\ 
\mathrm{Ricci}= \\ 
\left( \gamma_2+\left( 1+\gamma_2\right) y\right) \omega ^{1}\otimes
\omega ^{1} \\ 
+(1+\gamma_2+y^{2})\omega ^{1}\otimes \omega ^{2} \\ 
+\left( \gamma_2y+\left( 1+\gamma_2\right) y^{2}\right) \omega
^{2}\otimes \omega ^{1} \\ 
+\left( 1+y+\gamma_2y^{2}\right) \omega ^{2}\otimes \omega ^{2}
\end{array}
$ & no sol. \\ \hline
\end{tabular}
\caption{Ricci tensor and scalar for the algebra $F$, metric $g_{F.3}$.\label{t12}}
\end{table}
\end{tiny}

\begin{tiny}
\begin{table}[H] 
\begin{tabular}{|l|l|l|l|}
\hline
QLC & $\mathrm{Eins}=\mathrm{Ricci}+{S}g$ & $\mathrm{Ricci}$ qsymm & $\nabla
.\mathrm{Eins}=0$ \\ \hline
$\left. 
\begin{array}{c}
\nabla _{F.3.4} \\ 
\nabla _{F.3.7} \\ 
\nabla _{F.3.9} \\ 
\nabla _{F.3.10} \\ 
\nabla _{F.3.12}
\end{array}
\right\} $ & $\mathrm{Eins}=0$ (flat connections)& --- & --- \\ \hline
$\nabla _{F.3.1}$ & $
\begin{array}{l}
\mathrm{Eins}= \left( \gamma_1+\gamma_2+\gamma_3(1+y^{2})\right) \omega
^{1}\otimes \omega ^{1} \\ 
+\left( 1+\gamma_1+\gamma_3+\gamma_2y+(1+\gamma_2)y^{2}\right)
\omega ^{1}\otimes \omega ^{2} \\ 
+\left(\gamma_3+\left( \gamma_1+\gamma_3\right)
y+\left( 1+\gamma_2\right)(1+ y^{2})\right) \omega ^{2}\otimes \omega ^{1} \\ 
+\left( \gamma_1+\gamma_3+\gamma_3y+(1+\gamma_2+\gamma
_{3})y^{2}\right) \omega ^{2}\otimes \omega ^{2}
\end{array}
$ & $
\begin{array}{l}
\mathrm{Eins}= \\ 
\left( 1+\gamma_2+\gamma_3y^{2}\right) \omega ^{1}\otimes \omega ^{1}
\\ 
+\left( \gamma_2y+(1+\gamma_2)y^{2}\right) \omega ^{1}\otimes \omega
^{2} \\ 
+\left( \gamma_2+y+\left( 1+\gamma_2\right) y^{2}\right) \omega
^{2}\otimes \omega ^{1} \\ 
+\left( 1+y+(1+\gamma_2)y^{2}\right) \omega ^{2}\otimes \omega ^{2}
\end{array}
$ & no sol. \\ \hline
$\nabla _{F.3.2}$ & $
\begin{array}{l}
\mathrm{Eins}= \left( (1+\gamma_3)y+\gamma_1y^{2}\right) \omega ^{1}\otimes \omega
^{1} \\ 
+\left( \gamma_1+(1+\gamma_3)(y+y^{2})\right) \omega
^{1}\otimes \omega ^{2}
\end{array}
$ & $\mathrm{Eins}=0$ & $\gamma_1=0,\gamma_3=1$ \\ \hline
$\nabla _{F.3.3}$ & $
\begin{array}{l}
\mathrm{Eins}= \\ 
\left( \gamma_3+\left( 1+\gamma_2\right) (1+y)+(1+\gamma
_{3})y^{2}\right) \omega ^{1}\otimes \omega ^{1} \\ 
+\left( 1+\gamma_1+(1+\gamma_2)y+y^{2}\right) \omega ^{1}\otimes
\omega ^{2} \\ 
+\left( 1+\gamma_1+\left( 1+\gamma_2+\gamma_3\right)
y+\gamma_2(1+y^{2})\right) \omega ^{2}\otimes \omega ^{1} \\ 
+\left( \gamma_2+(1+\gamma_3)y+(\gamma
_{1}+\gamma_3)(1+y^{2})\right) \omega ^{2}\otimes \omega ^{2}
\end{array}
$ & $
\begin{array}{l}
\mathrm{Eins}= \\ 
\left( \gamma_2+\left( 1+\gamma_2\right) y\right) \omega ^{1}\otimes
\omega ^{1} \\ 
+\left( 1+(1+\gamma_2)y+y^{2}\right) \omega ^{1}\otimes \omega ^{2} \\ 
+\left( 1+\gamma_2\left( 1+y+y^{2}\right) \right) \omega ^{2}\otimes
\omega ^{1} \\ 
+\left( 1+\gamma_2+y^{2}\right) \omega ^{2}\otimes \omega ^{2}
\end{array}
$ & no sol. \\ \hline
$\nabla _{F.3.5}$ & $
\begin{array}{l}
\mathrm{Eins}= \\ 
\left( \left( 1+\gamma_1\right)(1+ y)+(1+\gamma_1+\gamma
_{2}+\gamma_3)y^{2}\right) \omega ^{1}\otimes \omega ^{1} \\ 
+\left( \gamma_1+\left( 1+\gamma_1+\gamma_3\right) y^{2}\right)
\omega ^{1}\otimes \omega ^{2} \\ 
+\left( \gamma_1+\gamma_2+\gamma_3(1+y^{2})\right) \omega
^{2}\otimes \omega ^{1} \\ 
+\left(1+(1+\gamma_1+\gamma
_{2}+\gamma_3)(1+y)+\gamma_2y^{2}\right) \omega ^{2}\otimes \omega ^{2}
\end{array}
$ & $
\begin{array}{l}
\mathrm{Eins}= \\ 
\left( 1+y+\gamma_2y^{2}\right) \omega ^{1}\otimes \omega ^{1} \\ 
+\left( 1+\gamma_2+y^{2}\right) \omega ^{2}\otimes \omega ^{1} \\ 
+\left( 1+\gamma_2\left( 1+y+y^{2}\right) \right) \omega ^{2}\otimes
\omega ^{2}
\end{array}
$ & no sol. \\ \hline
$\nabla _{F.3.6}$ & $
\begin{array}{l}
\mathrm{Eins}= \\ 
\left( 1+\left( \gamma_1+\gamma
_{3}\right) y+(\gamma_1+\gamma_2+\gamma_3)(1+y^{2})\right) \omega
^{1}\otimes \omega ^{1} \\ 
+\left( \gamma_2+\gamma_3+\left( 1+\gamma_1\right) y^{2}\right)
\omega ^{1}\otimes \omega ^{2} \\ 
+\left( \gamma_3+\left( 1+\gamma_2\right) (1+y)+(1+\gamma
_{3})y^{2}\right) \omega ^{2}\otimes \omega ^{1} \\ 
+\left( \gamma_1+\gamma_2+(\gamma_1+\gamma_2+\gamma
_{3})y\right) \omega ^{2}\otimes \omega ^{2}
\end{array}
$ & $
\begin{array}{l}
\mathrm{Eins}= \\ 
\left( \gamma_2+y+(1+\gamma_2)y^{2}\right) \omega ^{1}\otimes \omega ^{1}
\\ 
+\left( 1+\gamma_2+y^{2}\right) \omega ^{1}\otimes \omega ^{2} \\ 
+\left( \gamma_2+\left( 1+\gamma_2\right) y\right) \omega ^{2}\otimes
\omega ^{1} \\ 
+\left( \gamma_2+(1+\gamma_2)y\right) \omega ^{2}\otimes \omega ^{2}
\end{array}
$ & no sol. \\ \hline
$\nabla _{F.3.8}$ & $
\begin{array}{l}
\mathrm{Eins}= \\ 
\left( \gamma_2+y+(1+\gamma_1+\gamma_2)y^{2}\right) \omega
^{1}\otimes \omega ^{1} \\ 
+\left( \left( 1+\gamma_1+\gamma_2+\gamma_3\right) y+\left(
\gamma_2+\gamma_3\right) y^{2}\right) \omega ^{1}\otimes \omega ^{2} \\ 
+\left( 1+\gamma_1+\left( \gamma_1+\gamma_3\right) y+(1+\gamma
_{2}+\gamma_3)y^{2}\right) \omega ^{2}\otimes \omega ^{1} \\ 
+\left( \gamma_2+(1+\gamma_1+\gamma_2)y+(1+\gamma_3)y^{2}\right)
\omega ^{2}\otimes \omega ^{2}
\end{array}
$ & $
\begin{array}{l}
\mathrm{Eins}= \\ 
\left( \gamma_2+y+(1+\gamma_2)y^{2}\right) \omega ^{1}\otimes \omega
^{1} \\ 
+\left( \gamma_2(y+y^{2})+y^2\right) \omega ^{1}\otimes \omega ^{2} \\ 
+\left( 1+y+\gamma_2y^{2}\right) \omega ^{2}\otimes \omega ^{1} \\ 
+\left( \gamma_2+(1+\gamma_2)y\right) \omega ^{2}\otimes \omega ^{2}
\end{array}
$ & no sol. \\ \hline
$\nabla _{F.3.11}$ & $
\begin{array}{l}
\mathrm{Eins}= \\ 
\left( 1+\gamma_1+\gamma_3+(1+\gamma_2+\gamma_3)y+(1+\gamma
_{1}+\gamma_2)y^{2}\right) \omega ^{1}\otimes \omega ^{1} \\ 
+\left( \gamma_1+\gamma_2+(\gamma_1+\gamma_2+\gamma
_{3})y\right) \omega ^{1}\otimes \omega ^{2} \\ 
+\left( \gamma_1+\gamma_3(1+y)+(1+\gamma_2+\gamma
_{3})y^{2}\right) \omega ^{2}\otimes \omega ^{1} \\ 
+\left( \gamma_3+(1+\gamma_1+\gamma_2)y+(1+\gamma_3)y^{2}\right)
\omega ^{2}\otimes \omega ^{2}
\end{array}
$ & $
\begin{array}{l}
\mathrm{Eins}= \\ 
\left( \gamma_2y+(1+\gamma_2)y^{2}\right) \omega ^{1}\otimes \omega
^{1} \\ 
+\left( \gamma_2+(1+\gamma_2)y\right) \omega ^{1}\otimes \omega ^{2}
\\ 
+\left( 1+y+\gamma_2y^{2}\right) \omega ^{2}\otimes \omega ^{1} \\ 
+\left( 1+(1+\gamma_2)y\right) \omega ^{2}\otimes \omega ^{2}
\end{array}
$ & no sol. \\ \hline
\end{tabular}
\caption{Einstein tensor for the algebra $F$, metric $g_{F.3}$ showing a unique connection $\nabla_{F.3.2}$ which is not Ricci flat but has ${\rm Eins}=0$.\label{t13}}
\end{table}
\end{tiny}
\begin{tiny}
\begin{table}[H] 
\begin{tabular}{|l|l|l|l|l|}
\hline
QLC & Ricci (central for all $\gamma _{i}$)& $S=(\ ,\ )(\mathrm{Ricci})$ &  $\nabla .
\mathrm{Ricci=0}$ \\ \hline
$\left. 
\begin{array}{c}
\nabla _{F.4.1} \\ 
\nabla _{F.4.3} \\ 
\nabla _{F.4.4}
\end{array}
\right\} $ & $\mathrm{Ricci}=0$ (flat connections)& $S=0$ & --- \\ \hline
$\nabla _{F.4.2}$ & $
\begin{array}{l}
\mathrm{Ricci}= \\ 
\left( 1+\gamma_1+\left( \gamma_2+\gamma_3\right) y+\left( 1+\gamma
_{1}\right) y^{2}\right) \omega ^{1}\otimes \omega ^{1} \\ 
+\left( \left( 1+\gamma_1+\gamma_2+\gamma_1\right) y+\left( \gamma
_{2}+\gamma_3\right) y^{2}\right) \omega ^{1}\otimes \omega ^{2} \\ 
+\left( 1+\gamma_1+\gamma_3+\gamma_1y+\left( 1+\gamma_1+\gamma
_{3}\right) y^{2}\right) \omega ^{2}\otimes \omega ^{1} \\ 
+\left( \left( 1+\gamma_3\right) y+\gamma_1y^{2}\right) \omega
^{2}\otimes \omega ^{2} 
\end{array}
$ & $
\begin{array}{l}
S= \\ 
1+\gamma_1+\gamma_3 \\ 
+\left( 1+\gamma_2+\gamma_3\right) y \\ 
+\left( 1+\gamma_1+\gamma_2\right) y^{2}
\end{array}
$  & $\begin{array}{l}
\gamma_2=\gamma_1+1,\gamma_3=\gamma_1:\\
\mathrm{Ricci}= \\ 
\left( y+\left( 1+\gamma
_{1}\right) (1+y^{2})\right) \omega ^{1}\otimes \omega ^{1} \\ 
+\left( \gamma_1 y+ y^{2}\right) \omega ^{1}\otimes \omega ^{2} \\ 
+\left( 1+\gamma_1y+ y^{2}\right) \omega ^{2}\otimes \omega ^{1} \\ 
+\left( \left( 1+\gamma_1\right) y+\gamma_1y^{2}\right) \omega
^{2}\otimes \omega ^{2} \end{array}
$ \\ 
\hline
\end{tabular}
\caption{Ricci tensor and scalar for the algebra F, metric $g_{F.4}$. Three connections are Ricci flat and $\nabla_{F.4.2}$ never has Ricci quantum symmetric, so that
column is omitted. \label{t14}}
\end{table}
\end{tiny}

\begin{tiny}
\begin{table}[H] 
\begin{tabular}{|l|l|l|l|}
\hline
QLC & $\mathrm{Eins}=\mathrm{Ricci}+{S}g$  & $\nabla
.\mathrm{Eins}=0$ \\ \hline
$\left. 
\begin{array}{c}
\nabla _{F.4.1} \\ 
\nabla _{F.4.3} \\ 
\nabla _{F.4.4}
\end{array}
\right\} $ & $\mathrm{Eins}=0$ (flat connections)&  --- \\ \hline
$\nabla _{F.4.2}$ & $
\begin{array}{l}
\mathrm{Eins}= \\ 
\left( \gamma_1+(\gamma_1+\gamma_2)y+(1+\gamma_2)y^{2}\right)
\omega ^{1}\otimes \omega ^{1} \\ 
+\left( 1+\gamma_1+\gamma_2+\gamma_2y+(\gamma_1+\gamma
_{2})y^{2}\right) \omega ^{1}\otimes \omega ^{2} \\ 
+\left( \gamma_2+\gamma_3+\left( 1+\gamma_3\right) y+y^{2}\right)
\omega ^{2}\otimes \omega ^{1} \\ 
+\left( 1+\gamma_2+\gamma_3+\left( 1+\gamma_2\right) y+(1+\gamma
_{3})y^{2}\right) \omega ^{2}\otimes \omega ^{2}
\end{array}
$ & $ \begin{array}{l}
 \gamma_2=1+\gamma_1,\gamma_3=\gamma_1:\\
\mathrm{Eins}= \left( y+\gamma_1(1+y^{2})\right)
\omega ^{1}\otimes \omega ^{1} \\ 
\qquad+\left( (1+\gamma_1)y+y^{2}\right) \omega ^{1}\otimes \omega ^{2} \\ 
\qquad+\left( 1+\left( 1+\gamma_1\right) y+y^{2}\right)
\omega ^{2}\otimes \omega ^{1} \\ 
\qquad+\left( \gamma_1 y+(1+\gamma
_{1})y^{2}\right) \omega ^{2}\otimes \omega ^{2}
\end{array}$

\\ \hline
\end{tabular}
\caption{Einstein tensor for the algebra F, metric $g_{F.4}$ showing  two lifts for $\nabla_{F.4.2}$ with $\nabla\cdot{\rm Eins}=0$ and $S=1$. \label{t15}} 
\end{table}
\end{tiny}

We see for $g_{F.1}$ -- $g_{F.3}$  with $\underline{\rm dim}=1$ that for each metric there are five flat connections. Of the other connections, we see that for each metric there are two liftings  which render all Ricci quantum symmetric (e.g. for the first metric the lift is $\gamma_1=\gamma_2\in\{0,1\}, \gamma_3=0$) and resulting in a unique connection which is not Ricci flat but  has Einstein=0, namely $\nabla_{F.1.10},\nabla_{F.2.3},\nabla_{F.3.2}$ respectively. This is also the only case for each metric where $\nabla\cdot{\rm Eins}=0$.  Indeed, the other cases have Ricci not central, which implies that it could not be any multiple of the metric. It is also striking that all the other connections in this group have the same value of $S$. By contrast, the metric $g_{F.4}$ has $\underline{\rm dim}=0$ and behaves more like $g_{D.3}$ and $g_B$ above. It has three flat connections, and the remaining connection never has Ricci quantum symmetric but has two lifts that render $\nabla\cdot{\rm Eins}=0$.

\section{Conclusions} \label{secfin}

In this paper we have mapped out the landscape of all reasonable up to 2D quantum geometries over the field $\F_2$ on unital  algebras of vector space dimension $n\le 3$. The interesting ones up to this dimension have commutative coordinate algebras, which would mean the algebra of functions on up to 3 `points' if we were working over $\C$, but over a finite field such as $\F_2$ we have more possibilities.  We used the constructive `bimodule connection' approach \cite{Ma:ltcc,DVM,Mou,BegMa1,BegMa2,Ma:gra,MaTao1}  in which the layers of geometry are added one at a time starting with a calculus $\Omega^1$ free and of dimension $m\le 2$ over the algebra.  For the exterior algebra we focussed on the case of $\Omega^2$ free and 1-dimensional with a central basis element $\Vol$, so like a 2-manifold. For $n=2$ we also covered the case of $\Omega^2=0$ as for a 1-manifold. 

The first striking conclusion is that even under this restricted set of assumptions there are a lot of such `digital' finite quantum geometries by the time we get to $n=3$. For $n=2$ there are only a few geometries. First, the calculus for $n=2$ has to have $m=1$ and there are no calculi  with $\Omega^2$ top that admit a strictly quantum-symmetric metric. If we relax that then each of the three algebras $\F_2\Z_2, \F_2(\Z_2),\F_4$ admit only the flat metric $\cg=\omega\tens\omega$ with the zero connection $\nabla \omega=0$ if we want a QLC, see Table~\ref{t1}.  If we insist on quantum symmetry of the metric, as we do elsewhere in the classification, then $n=2$ forces us to $\Omega^2=0$ (so a 1D geometry from the top form point of view), the same flat metric and now respectively 2, 1, 3 QLCs for the three algebras (the table also shows more options if we allow the weaker requirement of a WQLC).  Many more quantum geometries emerge for $n=3$. First off, there are 6 possible commutative algebras as already known from another context \cite{MaPac} and we find one further noncommutative one. But none of them meet our requirements for a well-behaved calculus $\Omega^1$ of dimension $m=1$ while still admitting a quantum metric that meets the invertibility axiom (there are some examples if we drop this, see Tables~\ref{t2} and~\ref{t3}).  These also have issues with $\Omega^2$ if we take this to be non-zero. There is also a noncommutative algebra G which for $m=1$ again does not admit a suitable calculus having a central metric.  Therefore the landscape at $n=3$ properly needs $m=2$. In this case we find that only three  of the six algebras, namely $\rm{B}=\F_2(\Z_3)$, $\rm{D}=\F_2\Z_3$, $\rm{F}=\F_8$, meet our full requirements on the calculus including $\Omega^2$ as top degree and existence of a quantum symmetric metric. For each algebra we find an essentially unique calculus and a unique quantum metric up to an invertible functional factor, giving respectively 1, 3, 4 quantum metrics  that admit QLCs. Between them there are respectively 4, 12, 40 metric and QLC pairs (or `quantum Riemannian geometries') of which 3, 3, 18 are flat in the sense of zero Riemann curvature $R_\nabla$,  see Sections~\ref{secBmodel}, \ref{secDmodel} and \ref{secFmodel} respectively. An interesting feature in the D case is that for each metric the nonzero curvatures are the same even though the connections are different.  These results suggest an even richer moduli of quantum geometries when $n\ge 4$ but  beyond reach of our current method of trying all possible $2^{24}$ Christoffel symbol values to find the QLCs.

We also used our landscape of quantum Riemannian geometries to study the canonical geometric Laplacian $\Delta$ and Ricci tensor.  For the former in Section~\ref{seclap}, a striking observation that holds across all the viable $n=3,m=2$ quantum geometries is that $\Delta=0$ if and only if the quantum dimension $\underline{\rm dim}=0$. In the $\underline{\rm dim}=1$ case we found that the trace of $\Delta$ determines if there is a massive eigenvector (i.e. eigenvalue 1) or not, see Proposition~\ref{proplap}, resulting in 6 Laplacians on  $\F_8$ that have this massive eigenvector, none for the other geometries. Another feature is that $\Delta$ always depends on the connection with a four-fold degeneracy (four connections give the same $\Delta$)  with the result that it does not depend on the connections for $\F_2(\Z_3)$ and $\F_2\Z_3$ but only on the metric, while for $\F_8$ this is also true for one of the metrics $g_{F.4}$ but for each of the other three metrics the 12 connections are divided  into groups of four. It will be interesting to see if any of these features extend as we increase the dimension. 

 For the Ricci tensor and scalar $S$ in Section~\ref{secricci}  we used an approach \cite{BegMa2} that depends on a lifting map $i$. A  corollary of our analysis of  quantum metrics on the B,D,F algebras is that the possible lifts form an affine space taking the form  $i(\Vol)=I_0+\gamma \cg$ where $\cg$ is any fixed quantum metric, $I_0$ is any fixed central 1-1 tensor with $\wedge I_0=\Vol$ and the parameter $\gamma$ is an element of the algebra (so there are 8 possible lifts) as featuring in Tables~\ref{t4} - \ref{t15}. We also tentatively proposed over $\F_2$ to take the Einstein tensor as Ricci+$S\cg$ (given that the usual factor -1/2 makes no sense). This worked as expected for $\F_2\Z_3$ with its $g_{D.1},g_{D.2}$ metrics in the sense that there exist liftings such that Ricci is quantum symmetric and then Einstein=0 independently of which lift and which connection (just as in classical geometry in dimension 2). It also worked for $\F_8$ with each of its $g_{F.1},g_{F.2},g_{F.3}$ metrics in the limited sense that Ricci could always be made quantum symmetric independently of the connection and among the QLCs there was a unique one with Einstein=0.  These are all quantum metrics with $\underline{\rm dim}=1$. By contrast the $g_{D.3}$, $g_B$ and $g_{F.4}$ metrics with $\underline{\rm dim}=0$ followed a non-classical pattern with Ricci never quantum-symmetric for any lift, but instead we found a lift existing such that $\nabla\cdot {\rm Einstein}=0$ holds.   It should be recalled that in quantum Riemannian geometry the QLC condition is linear plus quadratic in the Christoffel symbols and quite typically has classical-like solutions (sometimes unique) and non-classical ones\cite{BegMa2}. 
 
Looking forward, interesting quantum geometries over $\F_2$ for $n=4$ and higher certainly exist, for example as special cases of results known over $\C$ adapted with care over any field (this is possible in at least a few cases) and then specialised to $\F_2$. Intrinsically $\F_p$-geometries for any prime $p$ were introduced in \cite{BasMa} as the Hopf algebras $A_d=\F_p[x]$ with the relation $x^{p^d}=x$  and a natural translation-invariant differential calculus. The $n=4$ case $A_2$ over $\F_2$ was solved for its three translation-invariant quantum metrics to find in each case two translation-invariant QLCs (as natural examples, rather than a moduli of all quantum geometries on the algebra). Aside from the landscape in higher dimensions of algebra and calculus, it would also be interesting to see which of our solutions extend to higher $\F_{2^d}$ and to other $\F_{p^d}$ and $\C$. As mentioned in the introduction,  the finite field setting also allows one to test definitions and conjectures that are expected to hold over any field, even if  the main interest is over $\C$.  The $n=4, m\le 2$ case should also reveal more interesting examples of `diffeomorphisms' alluded to in the preliminaries (algebra automorphisms compatible with potentially different differential structures before and after). These are visible in Table~\ref{t2} where calculi B.5-B.8 are equivalent to B.1-B.4 in this way, but we have not explored the topic systematically since the $n=3, m=1$ case is not sufficiently interesting, while for the other cases there was only the universal calculus in the first place. 
 
As also discussed in the introduction, once we have a good handle on the moduli of classes of small $\F_{p^d}$ quantum Riemannian geometries, we can then consider quantum gravity, for example as a weighted sum over the moduli space of them much as in lattice approximations\cite{Lai}, but now finite. One may also consider how  quantum geometries could  develop by transitions much as in the dynamical poset approach \cite{Rid,Dow}, as well as finite versions of other established approaches. This is another direction for further work.

\begin{figure}
\includegraphics[scale=0.13]{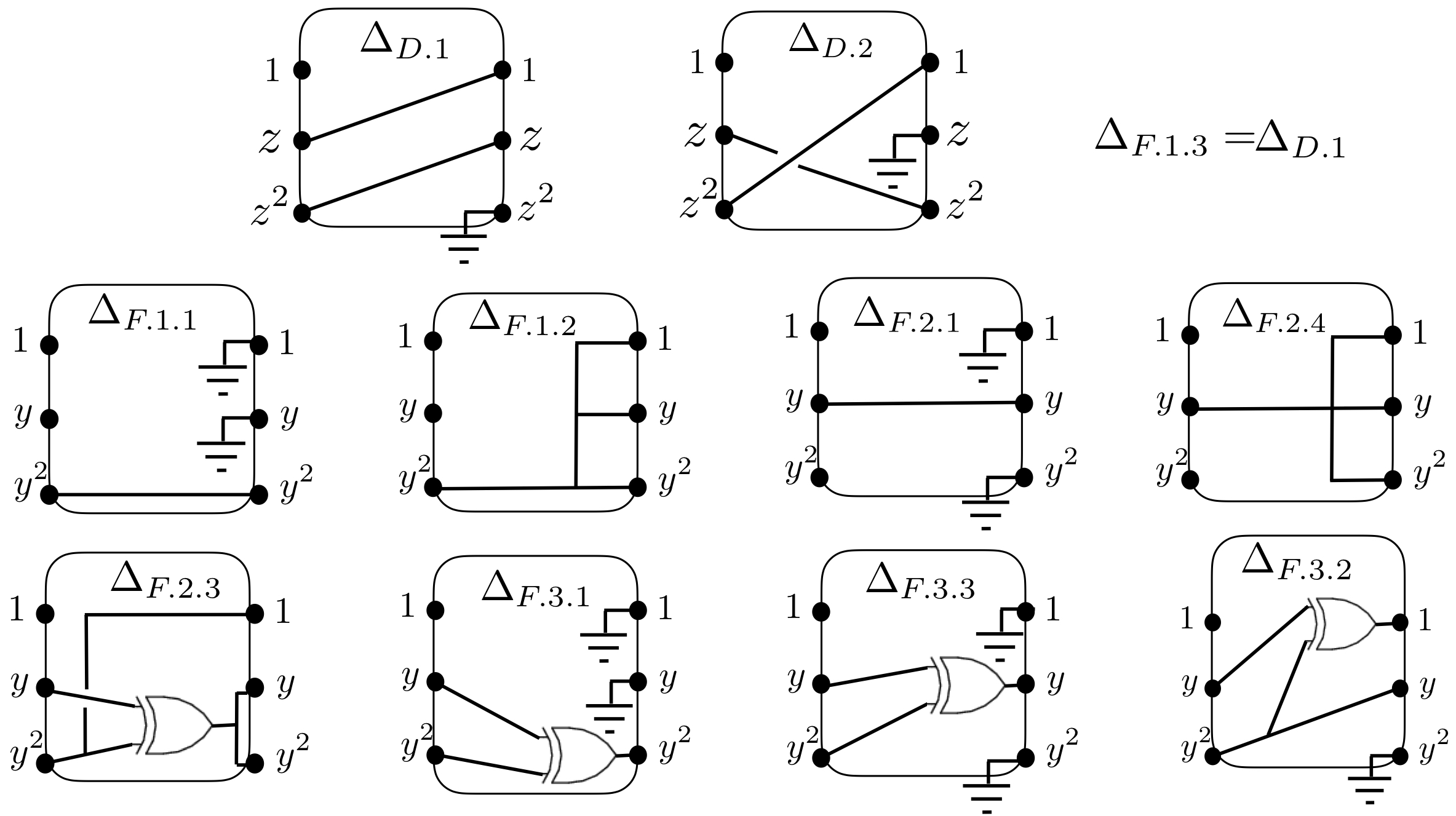}
\caption{Wiring diagrams for the nonzero Laplacians found in Section~\ref{seclap} labelled for the F algebra by a representative connection and for the D algebra by the metric. A live wire in the input at left or output at right means that the indicated algebra basis element is included. The gate shown is XOR. \label{diag}}
\end{figure}

 Finally, just as geometry has many applications, we envisage many applications of `finite' quantum geometry both over $\F_2$ and more generally over $\F_{p^d}$ (as well as over $\C$). It is not clear to what extent physics entirely over $\F_{p^d}$ makes sense  but this could be interesting to explore in terms of quantum mechanics. For quantum field theory the second quantisation can be done over $\C$ working with functions on the discrete moduli of finite solutions of the Klein-Gordon equation over $\F_2$ defined by $\Delta$. Quantum mechanics fully over $\F_2$ is unlikely to have a physical meaning but as an analogous formalism it may lead to `quantum geometric' constructions for gates in a `digital quantum computer' (as well as actual quantum geometric gates over $\C$). Discrete geometric ideas with real or complex coefficients are also used in network theory \cite{Bia} and finite versions might  be useful. Although these ideas are currently speculative, a  first step could be the Laplacian for the quantum geometry. As a map $\Delta: A\to A$, this can be realised digitally by choosing a basis of $A$. Each basis element then corresponds to an input wire with an element of $A$ specified by those basis elements where the wire is active. Similarly for the output copy of $A$. In this notation the non-zero Laplacians in Section~\ref{seclap} are shown in Figure~\ref{diag} as labelled by the metric or by a representative connection. Only two input wires are effective as $1$ is in the kernel, and clearly Laplacians of practical interest would need to be somewhat more complicated. It is explained in \cite{MaPac} how to handle tensor products and the wiring diagrams for the algebra products of $\F_2\Z_2$, $\F_2(\Z_2)$, $\F_4$ are given there. Such operations and their possible applications to engineering constitute another direction for further work.

\end{document}